\documentstyle{article}
\begin{document}
\begin{large}

\title{ On an inverse problem
for finite-difference operators of second order}
\author{  Mikhail Kudryavtsev}
\maketitle

\medskip
\begin{center}
{\it The Institute for Low Temperature Physics and Engineering\\
of the National Academy of Science of Ukraine,\\
\smallskip
Ukraine, Kharkov, 61103, Lenine Ave., 47}\\
\medskip
E-mail: kudryavtsev@ilt.kharkov.ua
\end{center}

\bigskip

\abstract{

The Jacobi matrices with bounded elements whose spectrum of multiplicity
2 is separated from its simple spectrum and contains an interval of
absolutely continuous spectrum are considered. A new type of spectral data,
which are analogous for scattering data, is introduced for this matrix.
An integral equation that allows us to reconstruct the matrix from this
spectral data is obtained. We use this equation to solve the Cauchy problem
for the Toda lattice with the initial data that are not stabilized.

}

\bigskip

\centerline{\bf Introduction}

\bigskip

    {\bf 1.} \ The Cauchy problem for the equation of the oscillation
of the doubly-infinite Toda lattice
$$ { d\,^2 \, x_k \over d \, t^2 } =
   e^{x_{k+1} - x_k} - e^{x_k - x_{k-1}} \,, \quad  k\in {\bf Z} , \eqno(0.1)
$$
$$  x_k(0)=v_k, \qquad  \dot x_k (0) = w_k\,, \eqno(0.2)
$$
with stabilized initial data
$$ \sum_{-\infty}^\infty |k| |a_k| < \infty \,,  \qquad
   \sum_{-\infty}^\infty |k| |b_k-1| < \infty \,,
$$
where
$$ a_k= w_k, \qquad b_{k-1} = {\rm exp} {v_k - v_{k-1} \over 2},  \eqno(0.3)
$$
is solved by Manakov and Flashka by the method of inverse
scattering problem (see [1],[2]). Here the finite-difference operator
of second order with the matrix
$$
J=\pmatrix{
   \ddots&\ddots&\ddots&{}&{}&{}&{}&{}\cr
   {}&b_{-2}&a_{-1}&b_{-1}&{}&0&{}&{}\cr
   {}&{}&b_{-1}&a_0&b_0&{}&{}&{}\cr
   {}&0&{}&b_0&a_1&b_1&{}&{}\cr
   {}&{}&{}&{}&b_1&a_2&b_2&{}\cr
   {}&{}&{}&{}&{}&\ddots&\ddots&\ddots\cr},
   \eqno(0.4)
$$
where the coefficients $a_k$, $b_k$ are defined by (0.3),
plays the role of the $L$-operator.

   However, we cannot apply the method of inverse scattering to solve
the Cauchy problem with non-stabilized initial data due to the fact that
the objects of with the inverse scattering problem (the Jost
solutions, the reflection coefficient, etc.) do not
exist in this case. In this connection it is important to find a more
general inverse spectral problem for the Jacobi matrices (0.4) so that
this problem can be applied to solve the equation of the oscillation of
the Toda lattice with quite a wide class of initial data.

\medskip

   Note, that in the case of the semi-infinite Toda lattices
the corresponding Cauchy problem with arbitrary (bounded) initial
data is solved by M.~Kac [3] and Yu.M.~Berezanski [4] by means of
the inverse problem of the reconstruction of the semi-infinite
Jacobi matrix by its spectral function. However, their method cannot
be extended to the case of the doubly-infinite Toda lattices.
Besides, the method of inverse problem is also used in the paper of
N.V.~Zhernakov [5], where the Cauchy problem for the doubly-infinite
Toda lattice is solved in the case when the operator defined by (0.4)
is the Hilbert-Schmidt operator. But this method cannot be used when
the spectrum of the operator is not discrete. Our aim is to solve the
inverse problem with more general type of spectrum.

\bigskip

   We denote by $P_k(\lambda)$, $ Q_k(\lambda)$ the solutions of the equation
$$b_{k-1}\omega_{k-1} + (a_k - \lambda) \omega_{k} + b_k \omega_{k+1} = 0,
\quad k\in {\bf Z},  \eqno(0.5)
$$
with initial data
$ \  P_0 (\lambda) = 1, \  P_{-1} (\lambda) = 0, \
    Q_0 (\lambda) =0, \  Q_{-1} (\lambda) =1.
$
It is evident that every solution of this equation is their linear
combination.

   It is known (see, for example, [5], [6]) that for
$\lambda \in {\bf C}\backslash{\bf R} $
the equation (0.5) has two solutions (the Weyl solutions of the matrix $J$):
$$ \varphi^R(k,\lambda) = m^R(\lambda)  P_k(\lambda)
                       -  {Q_k(\lambda) \over b_{-1}}, \quad k \in {\bf Z},
   \eqno(0.6)
$$
$$ \varphi^L(k,\lambda) = - {P_k(\lambda) \over b_{-1}}
                       + m^L(\lambda) Q_k(\lambda), \quad k \in {\bf Z},
   \eqno(0.7)
$$
such that $\sum_{k=N}^\infty  |\varphi^R(k,\lambda)|^2<\infty$, \
$\sum_{-\infty}^{k=N} |\varphi^L(k,\lambda)|^2<\infty$ for any
finite $N \in {\bf Z}$.
In this notation $m^R(\lambda)$ and $m^L(\lambda)$ are the Weyl functions
of the matrix $J$ represented in the form
$$ b_{-1} m^R(\lambda) = \int_{-\infty}^\infty
     {d\rho_R(\tau) \over \tau - \lambda } \,,  \qquad
   -{1\over b_{-1} m^L(\lambda)} = {\lambda\over b_{-1}} + \beta +
     \int_{-\infty}^\infty {d\rho_L(\tau) \over \tau - \lambda } \,,
    \eqno(0.8)
$$
where $\beta \in {\bf R}$, and $\rho_R(\lambda)$ and $\rho_L(\lambda)$
are nondecreasing functions. Here $d\rho_R(\lambda)$ is the spectral
measure of the semi-infinite (to the right) matrix $J^R$ that is contained
in $J$ and starts with the element $a_0$. The measure $d\rho_L(\lambda)$ is
spectral for the semi-infinite (to the left) matrix $J^L$ starting
with the element $a_{-2}$.

\bigskip

   {\it  In this paper we consider a larger class of matrices (0.4)
with bounded elements whose principal property is the existence of an
interval $[a,b]$ of absolutely continuous spectrum of multiplicity $2$,
which is separated from the other parts of the spectrum of the matrix $J$.
Without restriction of generality we can suppose that
$[a,b]=[-2,2]$ and $b_{-1}>0$.

   For such matrices a new inverse spectral problem is found and solved.
In this problem the functions $m^R(\lambda)$ and $m^L(\lambda)$ and the
number $b_{-1}$ play the role of spectral data from which the matrix $J$
is reconstructed. We obtain linear integral equation (3.46) of the
inverse problem and prove its unique solvability. All the parameters
contained in the equation are explicitly expressed by $b_{-1}$
and the functions $m^R(\lambda)$ and $m^L(\lambda)$.}

   In a forthcoming paper [8] we will apply this integral equation to solve
the corresponding Cauchy problem. We also note that one can find a short
and complete exposition of these results in papers [9], [10], resp.

\medskip

   In our case the functions $m^R(\lambda)$ and
$m^L(\lambda)$ have nonreal limits on the interval $[-2,2]$.
Instead of two Weyl functions and two Weyl solutions we introduce the
function
$$n(z)=  \cases{ - b_{-1} m^R(z+z^{-1}), \quad |z|<1, \cr
                 - \frac{1}{ b_{-1} m^L(z+z^{-1})}, \quad |z|>1, \cr }
  \eqno(0.9)
$$
and the solution
$$ \psi(k,z) = n(z) P_k(z+z^{-1}) + Q_k(z+z^{-1}) =
\cases{ - b_{-1} \varphi^R(k,z+z^{-1}), \quad |z|<1, \cr
     \frac{ \varphi^L(k,z+z^{-1}) }{ m^L(z+z^{-1}) }, \quad |z|>1, \cr}
   \eqno(0.10)
$$
which we will also call the Weyl function and the Weyl solution
of the matrix $J$.

   They are defined and holomorphic at such $z \in {\bf C}$, for which
both functions $m^R$ and $\frac{1}{m^L}$ are regular at the point $z+z^{-1}$.
At such points
$$\sum_N^\infty |\psi(k,z)|^2 <\infty, \ \ |z|<1, \qquad
\sum_{-\infty}^N |\psi(k,z)|^2 <\infty, \ \ |z|>1.
$$

\medskip

   Note, that in the papers of L.K.~Maslov [11],[12] the class of matrices
(0.4) was studied such that the limits of the Weyl functions $m^R(\lambda)$
and ${1\over m^L(\lambda)}$ coincide on a certain interval $(a,b)$ of the
real line:
$$ b_{-1} m^R(\lambda+i0) = {1\over b_{-1} m^L(\lambda-i0)},
   \quad a<\lambda<b.
$$
In this case $\psi(k,z)$ and $n(z)$ are holomorphic at the points
of the unit circle. In our case the limit values of functions
$\psi(k,z)$ and $n(z)$ from within and from the outside of the
circle $|\xi|=1$ are different.

\bigskip

   {\bf 2. \ Notation.} \ In the presented work we use the following
notations. Put the limit values of any function $f(z)$ on the unit circle
and the real line:
$$ f^\pm(e^{i\theta}) = \lim_{r\to 1 \mp 0} f(re^{i\theta}),
    \quad \theta \in (-\pi,\pi),
$$
and
$$ f^\pm(t) = f(t \pm i0), \quad t \in {\bf R}.
$$

   We will also denote by $\xi$ the points of the unit circle
$\{ |\xi|=1 \}$ in the complex plane.

\bigskip

   {\bf 3.} \ In order to clarify the analogy
between our case and the classical inverse scattering problem (ISP),
we will remind some essential points of ISP (see, for example, [1],[2]).
Like in the case of arbitrary Jacobi matrices, there exist two solutions
of equation (0.5), one of which is square-summable for positive $k$ and the
other one, for negative. In the case of stabilized initial data
the following asymptotic formulas are true for them:
$$\tilde\psi(k,z) = z^k+o(z^k), \quad k \to +\infty, \quad |z|<1, \eqno(0.11)
$$
$$\tilde\psi(k,z) = z^k+o(z^k), \quad k \to -\infty, \quad |z|>1, \eqno(0.12)
$$
where $z+z^{-1}=\lambda$ is the spectral parameter ($\tilde\psi(k,z)$,
$|z|<1$, are square-summable for positive $k$, and $\tilde\psi(k,z)$, $|z|>1$,
are square-summable for negative $k$). The functions $\tilde\psi(k,z)$
are called the Jost solutions.

   For a fixed point $\xi$ on the unit circle the limit values
$\tilde\psi^+(k,\xi)$, $\tilde\psi^-(k,\xi)$, $\tilde\psi^-(k,\xi^{-1})$
are solutions (with respect to $k$) of the equation (0.5) at the point
$\lambda=\xi+\xi^{-1} \in [-2,2]$. (These limit values exist because
the segment $[-2,2]$ belongs to the absolutely continuous spectrum;
here $\tilde\psi^-(k,\xi) = \overline{\tilde\psi^-(k,\xi^{-1})}$.)
The solution $\tilde\psi^+(k,\xi)$ is represented in the form of a linear
combination of two other solutions:
$$\tilde\psi^+(k,\xi) = \tilde c(\xi) \tilde\psi^-(k,\xi) +
  \tilde d(\xi) \tilde\psi^-(k,\xi^{-1}),
$$
or
$${1\over \tilde c(\xi)} \tilde \psi^+(k,\xi) =
    \tilde \psi^-(k,\xi) + {\tilde d(\xi)\over \tilde c(\xi)}
    \tilde \psi^-(k,\xi^{-1}) \,, \eqno(0.13)
$$
where
$$ {1\over \tilde c(\xi)} \equiv
   { <\tilde\psi^-(k,\xi^{-1}),\tilde\psi^-(k,\xi)> \over
     <\tilde\psi^-(k,\xi^{-1}),\tilde\psi^+(k,\xi)> } =
   { \xi - \xi^{-1} \over
     <\tilde\psi^-(k,\xi^{-1}),\tilde\psi^+(k,\xi)> }  \eqno(0.14)
$$
(the Wronskian $<\tilde\psi^-(k,\xi^{-1}),\tilde\psi^-(k,\xi)>$ is
explicitly calculated in view of the asymptotic formulae (0.11), (0.12)).
In this case the function
$ {1\over \tilde c(z)} \equiv  { z - z^{-1} \over
     <\tilde\psi(k,z^{-1}),\tilde\psi(k,z)> }$
is meromorphic inside the unit disk and its limit at the
point $\xi$ from within of the circle is equal to the coefficient
${1\over \tilde c(\xi)}$ from (0.13). Hence, this coefficient
can be analytically continued inside the circle. Therefore the function
$$ \Phi_k(z) = \cases{ {1\over \tilde c(z)} z^{-(k+1)} h_k \tilde \psi(k,z)
         , \quad |z|<1, \cr
         z^{-(k+1)} h_k \tilde \psi(k,z), \quad |z|>1, \cr }
$$
where
$$ h_k = \cases{ b_{-1} \ldots b_{k-1}, \quad   k \ge 0, \cr
                 1,                     \quad   k = - 1, \cr
                 {1\over b_k \ldots b_{-2}},     \quad   k \leq -2, \cr }
   \eqno(0.15)
$$
is meromorphic inside the unit disk. It has, according to (0.13), the jump
$$ \xi^{2(k+1)} \Bigl( \Phi_k^+(\xi)-\Phi_k^-(\xi) \Bigr)
   = r(\xi) \Phi_k^-(\xi^{-1})\,, \eqno(0.16)
$$
where $r(\xi) = {\tilde d(\xi) \over \tilde c(\xi)}$ is the reflection
coefficient, from which we reconstruct $J$. The function $\Phi_k(z)$
has inside the disk simple real poles $z_n,\ 0<|z_n|<1,\ n=1,2,\ldots,N$.
In these points $\tilde\psi(k,z_n) = \rho_n \tilde\psi(k,z_n^{-1})$
(these are the components of an eigenvector), and the residue of the function
$\Phi_k(z)$ equals ${z^{-2(k+1)}\rho_n\over\tilde c'(z_n)}\Phi_k(z_n^{-1})$.
Besides, the function $\Phi_k(z)\to 1$ as $z\to\infty$.
According to the Cauchy theorem,
$$\Phi_k(z) = 1 +
          \sum_{m=1}^N { z_m^{-2(k+1)}\rho_n \Phi_k(z_m^{-1}) \over
                                             \tilde c'(z_n) (z-z_m)}
      + {1\over 2\pi i} \int_{|\zeta|=1}
          {\zeta^{-2(k+1)}r(\zeta) \Phi_k(\zeta^{-1})
            \over \zeta-z} d\zeta \,.  \eqno(0.17)
$$
Let $z$ tend to the point $\xi^{-1}$, $|\xi^{-1}|=1$, from the outside
of the circle. According to the formulas of Plemelj--Sokhotski,
$$\Phi_k(\xi^{-1}) = 1 +
          \sum_{m=1}^N { z_m^{-2(k+1)}\rho_n \Phi_k(z_m^{-1}) \over
                                             \tilde c'(z_n) (\xi^{-1}-z_m)}
$$
$$ - {1\over2} \xi^{2(k+1)}r(\xi^{-1})\Phi_k^-(\xi)
   + {1\over 2\pi i} {\rm v.p.} \int_{|\zeta|=1}
          {\zeta^{-2(k+1)}r(\zeta) \Phi_k(\zeta^{-1})
            \over \zeta-\xi^{-1}} d\zeta \,. \eqno(0.18)
$$
(We observe that namely relation (0.16) let us replace
$\Phi_k^+(\zeta)-\Phi_k^-(\zeta)$ under the integral sign in (0.17) by
$$\zeta^{-2(k+1)}r(\zeta)\Phi_k^-(\zeta^{-1}).$$
Thus we have
one unknown function $\Phi_k^-(\xi^{-1})$ in both sides of the integral
equation (0.18).) Letting in (0.17) $z=z_n^{-1}$, we obtain $N$ equations
for $\Phi_k(z_n^{-1})$:
$$\Phi_k(z_n^{-1}) = 1 +
          \sum_{m=1}^N { z_m^{-2(k+1)}\rho_n \Phi_k(z_m^{-1}) \over
                                             \tilde c'(z_n) (z_n^{-1}-z_m)}
      + {1\over 2\pi i} \int_{|\zeta|=1}
          {\zeta^{-2(k+1)}r(\zeta) \Phi_k(\zeta^{-1})
            \over \zeta-z_n^{-1}} d\zeta \,. \eqno(0.19)
$$
Let us define the unknown function on the unit circle and at the points $z_n$
$$u_k(\beta):=\cases{ \beta^{-2(k+1)} \Phi_k^-(\beta^{-1}),
                                      \quad |\beta|=1, \cr
                      z_n^{-2(k+1)}\Phi_n^-(z_n^{-1}),
                                \quad \beta=z_n, \cr}  \eqno(0.20)
$$
and define at the points $z_n$ the point measure
$d\tilde \sigma(z_n):={\rho_n\over \tilde c'(z_n)}$. Let us also define
the function $\tilde\chi_0(\beta)$ as the indicator of the unit circle.
Then the system (0.18), (0.19) can be rewritten in the form of {\it one}
integral equation
$$\beta^{2(k+1)}u_k(\beta) +
  \int_{-\infty}^\infty {u_k(\alpha)\over \alpha-\beta^{-1}} d\tilde\sigma(\alpha) +
  {1\over2} \tilde\chi_0(\beta) r(\beta) u_k(\beta^{-1}) -
  {1\over 2\pi i} {\rm v.p.}\int_{|\xi|=1}
                  {r(\zeta)u_k(\zeta)\over\zeta-\beta^{-1}} d\zeta = 1, \eqno(0.21)
$$
where the {\it reflection coefficient} $r(\xi)$ and the {\it measure}
$d\tilde\sigma$ are the {\it scattering data}, and $u_k(\zeta)$ is the unknown
function. Further, the solution $\Phi_k(z)$ is represented in the form
$$\Phi_k(z) = 1 -
  \int_{-\infty}^\infty {u_k(\alpha)\over \alpha-z} d\tilde\sigma(\alpha) -
  {1\over 2\pi i} {\rm v.p.}\int_{|\zeta|=1}
                  {r(\zeta)u_k(\zeta)\over\zeta-z} d\zeta, \eqno(0.22)
$$
where $u_k(\beta)$ is the solution of equation (0.21). The elements
$a_k$ and $b_k$ of the matrix $J_N$ can be simply found from $\Phi_k$.

\smallskip

As we can see, the principal step of the derivation of the ISP-equation
(0.19) is to obtain the relation (0.16). To obtain this relation we
multiply the Jost solution $\tilde\psi(k,z)$ by the factor coefficient
that equals ${1\over\tilde c(z)}$ inside the disk and $1$ outside the disk
(we also multiplied the Jost solution by $z^{-2(k+1)}h_k$ in order to have
$\Phi_k(z)\to1,\ z\to\infty$).

   \medskip

   However, in our case the Jost solutions, which are defined by asymptotic
formulas (0.11), (0.12), do not exist. We consider the Weyl solutions
$\psi(k,z)$ defined by (0.10) (they exist for any Jacobi matrix).
It is easy to obtain for these solutions the following relations
on the unit circle $|\xi|=1$:
$${1\over a(\xi)} \psi^-(k,\xi) =
    \psi^+(k,\xi) + {b(\xi)\over a(\xi)} \psi^+(k,\xi^{-1}) \,, \eqno(0.23)
$$
$${1\over c(\xi)} \psi^+(k,\xi) =
    \psi^-(k,\xi) + {d(\xi)\over c(\xi)} \psi^-(k,\xi^{-1}) \, \eqno(0.24)
$$
(explicit expressions for $a(\xi), b(\xi), c(\xi),d(\xi)$ are given in
section~3). These relations are completely analogous to (0.13), but
for the coefficient
$$ {1\over c(\xi)} \equiv
   { <\psi^-(k,\xi^{-1}),\psi^-(k,\xi)> \over
     <\psi^-(k,\xi^{-1}),\psi^+(k,\xi)> }
$$
the equality (0.14) is not fulfilled, because the Weyl solutions do not
satisfy the asymptotic formulas (0.11), (0.12). From now on there is no
more coincidence between our case and classical inverse scattering, because
${1\over c(\xi)}$ does not have the necessary properties.

   Due to this reason, in order to obtain the equality, analogous to (0.16),
we have to multiply the Weyl solutions $\psi(k,z)$ by the specially chosen
function $R(z)$. The choice of the function $R(z)$ and the subsequent
reducing (0.23), (0.24) to a more symmetric look is done with the help
of the following theorem, which is a key result for the considered
inverse problem:

\medskip

{\bf Theorem 2 \ (factorization).} \ {\it The function
${(z-z^{-1})} {\Bigl( n(z)-n(z^{-1}) \Bigr) }^{-1}$
can be represented in the domain, where it is holomorphic,
in the form of the product of two functions $R(z)$, $R(z^{-1})$,
which are holomorphic in this domain
$$ - { {z-z^{-1}} \over {n(z)-n(z^{-1})} } \ = \  R(z) R(z^{-1}).  \eqno(0.25)
$$
The function $R(z)$ may only have singularities (i.e.~not to be holomorphic)
in such points $z$ that $z+z^{-1}$ belongs to the spectrum of the matrix
$J$. In addition, }
$$  R^+(\xi) R^+(\xi^{-1}) \, |{\rm Im}\, n^+(\xi)| =
    R^-(\xi) R^-(\xi^{-1}) \, |{\rm Im}\, n^-(\xi)| >0,  \quad |\xi|=1 .
$$
(The explicit definition of the function $R(z)$ is performed in section~2.)

\medskip

   Further, we define the function
$$ g(k,z) = {R(z) \over R(\infty)} z^{-(k+1)} h_k \psi(k,z),   \eqno(0.26)
$$
where $h_k$ are defined in (0.15). Theorem 2 allows us to derive from
equalities (0.23), (0.24) the following relation, which groups together
the limit values of the function $g(k,z)$ at the points $\xi$ and $\xi^{-1}$:
$$ \xi^{2(k+1)} (g^+(k,\xi)-g^-(k,\xi))
   =  - \hat r(\xi) (g^+(k,\xi^{-1})+g^-(k,\xi^{-1})) \,, \eqno(0.27)
$$
where $|\hat r(\xi)|<C<1$. (See explicit expression for $\hat r(\xi)$ in
section 3.)

   Precisely this relation is the analogue of the equality (0.16) and
makes it possible to derive the integral equation of the inverse problem,
which is appropriate for solving the Cauchy problem (0.1), (0.2).

\bigskip

{\bf 4.} \ We explain the principal steps of the
derivation of the inverse problem equation.
The transformation $z+z^{-1}=\lambda$ maps the spectrum of the matrix
$J$ to the unit circle and a certain set on the real line.
I turns out that the function $g(k,z)$ is holomorphic outside
of these sets. In view of theorem~1 (see section~1),
$g(k,z)\to1,\ z\to\infty$. Hence, according to the Cauchy theorem,
$$ g(k,z) = 1 + {1\over2\pi i}
   \int\limits_{\Gamma}{ g(k,\zeta) \over \zeta-z } d\zeta, \eqno(0.28)
$$
where $\Gamma$ is a contour enclosing all the singularities of
the function $g(k,z)$. These singularities are concentrated on the real
line and the unit circle and correspond to the different types of the spectrum
of the matrix $J$. Further we tighten the contour to the singularities
of the function $g(k,z)$ that are concentrated on the sets of each one
of these types:

1) on the set $\Omega_1\cup\Phi$, corresponding to:\\
a) disjoint parts of the spectra of the right semi-infinite matrix $J^R$
and the left semi-infinite matrix~$J^L$,\\
b) the eigenvalues of the matrix $J$ that are not contained in the spectra
of the two semi-infinite matrices;

2) on the set $\Gamma_2^s$, corresponding to the eigenvalues of the matrix
$J$ that are also eigenvalues of both matrices $J^R$ and $J^L$;

3) on the set $\Omega_2^a$, corresponding to the parts of the absolutely
continuous spectrum of multiplicity $2$ of the matrix $J$ except
the segment $[-2,2]$;

4) on the unit circle, corresponding to the segment $[-2,2]$ of
the absolutely continuous spectrum of multiplicity $2$ of the matrix $J$.

\smallskip

   Therefore the spectrum of the matrix $J$ splits into four
sets of different types. We calculate in the corresponding way the integrals
along the contours (see (3.2)) enclosing each one of these sets.
First, we remind that in the classical ISP there are parts of the spectrum
of two types: the simple spectrum (the finite set of eigenvalues) and the
segment of absolutely continuous spectrum of multiplicity $2$.
(The unit circle in $z$-plane corresponds to this segment.)
At the points of the simple spectrum we expressed the residue of the function
$\Phi_k(z)$ at the point $z_n$ by the value $\Phi_k(z^{-1})$,
substituted the obtained values of the residue in the integrals along the
contours, enclosing the eigenvalues, and obtained the sum in (0.17),
which is represented in the form of the integral by the measure
$d\tilde\sigma$. This measure is defined on the set of the
eigenvalues (i.e.~on the simple spectrum). In a resembling way,
in lemmas 3.1 and 3.2 we calculate the integrals from (0.27) along
the contours, enclosing the singularities on $\Omega_1\cup\Phi$ and
$\Gamma_2^s$, and represent these integrals in the form of integrals
along the set $\Omega_1\cup\Phi \subset {\bf R}$ and
$\Gamma_2^s \subset {\bf R}$ of the function $g_k(z^{-1})$.
In lemma 3.5 we obtain the relation (3.20) for the limit values
on the absolutely continuous spectrum (outside $[-2,2]$), which is analogous
to the relation (0.16) from ISP. Further we derive the main integral equation
(3.46) from the relations (3.20), (0.27) and lemmas 3.1 and 3.2
(see theorem~3). Formulas (3.47), (3.46) and (3.44) of the theorem 3,
are, evidently the analogs of the formulas (0.22), (0.21) and (0.20)
of the classical ISP. The functions $g(k,z)$ are simply expressed
from the solution of the integral equation (formula (3.47)).

   It is easy to conclude from the theorem 1 (formulas (1.3) and (1.4))
that the elements $a_k$ and $b_k$ of the matrix $J$ can be found
from $g(k,z)$ by formulas ($b_k$ up to its sign):
$$ b_{k-1}^2 = { g(k,0) \over g(k-1,0) },  \qquad
   a_k = \lim_{z\to\infty} { z(g(k,z)-g(k+1,z)) \over g(k,z)}. \eqno(0.29)
$$
So, the obtained inverse problem equation allows us to reconstruct
the matrix $J$ by its spectral data.

\bigskip

   {\bf 5.} \ Now we list the conditions on the spectrum for which
the solution of the inverse problem is given. Put
$$ \eta^R(\tau)= \lim_{\varepsilon \downarrow 0} \arg m^R(\tau+i\varepsilon), \qquad
   \eta^L(\tau)= \lim_{\varepsilon \downarrow 0} \arg {-1\over m^L(\tau+i\varepsilon)}.
$$

   Let $\tilde\Omega^R$ and $\tilde\Omega^L$ be the supports of measures
$d\rho_R(\tau)$ and $d\rho_L(\tau)$, resp., defined by the non-decreasing
functions $\rho_R(\tau)$ and $\rho_L(\tau)$ from (0.8),
$$ \tilde\Omega_2\ \equiv \ \tilde\Omega^R\cap \tilde\Omega^L,  \qquad
   \tilde\Omega_1\ \equiv \ (\tilde\Omega^R\backslash \tilde\Omega^L) \cup
                 (\tilde\Omega^L\backslash \tilde\Omega^R),
$$
$\tilde\Omega_2^s \subset \tilde\Omega_2$ be the set of the common poles
of the functions $m^R(\lambda)$ and ${-1\over m^L(\lambda)}$, and let
$\tilde\Omega_2^a \equiv \tilde\Omega_2 \backslash \tilde\Omega_2^s$.

\medskip

We assume that:

A) {\it All the three sets $\tilde\Omega_1$, $\tilde\Omega_2^s$,
$\tilde\Omega_2^a$ have positive mutual distances,
$[-2,2] \subset \tilde\Omega_2^a$, and the set
$\tilde\Omega_2^s$ is finite or empty.}

B) {\it For some $\varepsilon>0$ almost everywhere
(with respect to Lebesque measure) on the set
$\tilde\Omega^a_2$ }
$$ 0<\varepsilon<\eta^R(\alpha)<\pi-\varepsilon,  \quad
   0<\varepsilon<\eta^L(\alpha)<\pi-\varepsilon.
$$

C) {\it In some neighborhood of the set $\tilde\Omega^a_2$
the function $\eta^R(\alpha)-\eta^L(\alpha)$ satisfies the H\"older
condition. }

D) {\it The set $\tilde\Omega^a_2 \backslash [-2,2]$ can be covered
with mutually disjoint intervals $\delta_l$ on each of which the following
inequalities are true: }
$$
{\rm ess}\sup_{\alpha\in\delta_l}\eta^R(\alpha)-
{\rm ess}\inf_{\alpha\in\delta_l}\eta^R(\alpha)<\pi\ , \quad
{\rm ess}\sup_{\alpha\in\delta_l}\eta^L(\alpha)-
{\rm ess}\inf_{\alpha\in\delta_l}\eta^L(\alpha)<\pi\ ,
$$

E) {\it For some small $\varepsilon>0$
and $0<\alpha<\varepsilon$
$$\eta^R(2+\alpha)=\eta^L(2+\alpha)=0,   \qquad
  \eta^R(-2-\alpha)=\eta^L(-2-\alpha)=\pi,
$$
and the functions $\eta^R(\tau)$, $\eta^L(\tau)$ satisfy the
H\"older condition on the interval $[-2,2]$.
}

   From A)--E) it follows that the set  $\tilde\Omega_2^a$
is the absolutely continuous spectrum of multiplicity $2$
of the matrix $J$ and that the segment $[-2,2] \subset \tilde\Omega_2^a$
is divided from other parts of the spectrum.
We will need the conditions B)--D)
to prove that the limit values of the functions $g(k,z)$ on the
absolutely continuous spectrum of multiplicity $2$ belong to $L^2$.
(We have to add the condition E) for the segment [-2,2].)

\bigskip

   {\bf Remark.} \ The conditions A)--E) can be exposed only in the terms
of the functions $\eta^R(\tau)$ and $\eta^L(\tau)$. For this we have to
define the sets $\tilde\Omega^R$, $\tilde\Omega^L$,
$\tilde\Omega_1$, $\tilde\Omega_2$, $\tilde\Omega_2^s$
and $\tilde\Omega_2^a$ directly from $\eta^R(\tau)$ and $\eta^L(\tau)$,
without using the functions $\rho_R(\tau)$ and $\rho_L(\tau)$.
We can do it in the following way. (The equivalence of the two
definitions is a conclusion of the fact that $m^R(\lambda)$ and
${-1\over m^L(\lambda)}$ have a positive imaginary part in the upper
half-plane and are holomorphic in it.)

   Let $C\tilde\Omega^R$ ($C\tilde\Omega^L$) be the open set on the
real axis that is the union of intervals such that the function
$\eta^R(\tau)$ ($\eta^L(\tau)$) equals either $0$ or $\pi$ almost
everywhere on the interval and of the points such that
$\eta^R(\tau)$ ($\eta^L(\tau)$) is equal to $\pi$ in the left
neighborhood of the point and to $0$ in the right neighborhood.
(It means that $C\tilde\Omega^R$ ($C\tilde\Omega^L$) is the set on the real
axis where $m^R(\lambda)$ (${-1\over m^L(\lambda)}$) is holomorphic).
Let $\tilde\Omega^R$ ($\tilde\Omega^L$) be the complement of the set
$С\tilde\Omega^R$ ($С\tilde\Omega^L$) ,
$$ \tilde\Omega_2\ \equiv \ \tilde\Omega^R\cap \tilde\Omega^L,  \qquad
   \tilde\Omega_1\ \equiv \ (\tilde\Omega^R\backslash \tilde\Omega^L) \cup
                 (\tilde\Omega^L\backslash \tilde\Omega^R),
$$
$\tilde\Omega_2^s \subset \tilde\Omega_2$ be the set of points
such that $\eta^R(\tau)$, $\eta^L(\tau)$ both equal $0$ in the left
neighborhood of the point and equal $\pi$ in the right one,
$\tilde\Omega_2^a \equiv \tilde\Omega_2 \backslash \tilde\Omega_2^s$.

\medskip

    We also refer to paper [13] from which we take the idea of the
proofs of Lemmas 1.1, 1.2, 3.1--3.5, where this technique is used
for a differential operator.

\bigskip
\bigskip

\centerline {\bf 1. Asymptotics of the Weyl solutions and auxiliary results}

\bigskip

{\bf Theorem 1.\ {\rm (L.K.~Maslov)}.} {\it The Weyl solutions $\psi(k,z)$ and
the Weyl functions $n(z)$ of the matrix $J$, defined by (0.9), (0.10),
satisfy the asymptotic formulas
$${ {n(z)}\over{z} } \to b_{-1}, \quad  z \to 0, \eqno (1.1)
$$
$$\frac{n(z)}{z} \to {1\over b_{-1}}, \quad  z \to \infty, \eqno (1.2)
$$
$$ { {z^{-(k+1)}}\over{h_k} } \psi(k,z) \to 1, \quad  z \to 0, \eqno (1.3)
$$
$$ h_k z^{-(k+1)} \psi(k,z) \to 1, \quad  z \to \infty, \eqno (1.4)
$$
where $h_k$ are defined by formula (0.15).}

\medskip

We give the sketch of the proof, obtained by L.K.~Maslov in 1991.

  1. For $N>0$ we consider the matrix $J_N$ that is the matrix $J$
with changed right lower minor: instead of all $b_k,\  k\ge N-1,$ it has
$1$, and instead of $a_k,\  k\ge N,$ it has zeroes:
$$
J_N=\pmatrix{
   \ddots&\ddots&\ddots&{}&{}&{}&{}&{}&{}&{}\cr
   {}&b_{-1}&a_0&b_0&{}&{}&0&{}&{}&{}\cr
   {}&{}&b_0&a_1&b_1&{}&{}&{}&{}&{}\cr
   {}&{}&{}&\ddots&\ddots&\ddots&{}&{}&{}&{}\cr
   {}&{}&{}&{}&b_{N-3}&a_{N-2}&b_{N-2}&{}&{}&{}\cr
   {}&0&{}&{}&{}&b_{N-2}&a_{N-1}&1&{}&{}\cr
   {}&{}&{}&{}&{}&{}&1&0&1&{}\cr
   {}&{}&{}&{}&{}&{}&{}&\ddots&\ddots&\ddots\cr}.
$$
It is clear that
the orthogonal polynomials, corresponding to $J_N$, for $k\leq N-1$
coincide with $P_k(\lambda)$, $Q_k(\lambda)$. Hence for these
$k$ the Weyl solution of the matrix $J_N$, summable for positive $k$
(such that $\sum_{k=M}^\infty  |y_{N}^R(k,\lambda)|^2<\infty$),
is representable in the form
$$ y^R_{N}(k,\lambda) = m^R_N(\lambda) P_k (\lambda)
         - { Q_k (\lambda) \over b_{-1} }, \quad   k \leq N-1.
   \eqno(1.5)
$$
On the other hand, the following statement is true:

\medskip

{\bf Proposition.} {\it For $k\leq N-2$ the
solutions $ y^R_{N}(k,\lambda)$ are equal to
$$ \tilde y_k = \frac{z^k}{c_k}(1+zR_k(z)), \quad k \leq N-2, $$
to within a factor, not depending on $k$, with
$c_k=b_k b_{k+1} \ldots b_{N-2}$, and $R_k(z)$ is a polynomial
in $z$ (not of degree $k$).} (We remind that we parametrized $\lambda$
by the variable $z$, connected with $\lambda$ by relation
$z+z^{-1}=\lambda$.)

\medskip

{P r o o f} \  We use the induction in
$k=N-2, N-3, N-4, \ldots$. The function $\tilde y_{N-2}$ satisfies
the equation
$$b_{N-2} \tilde y_{N-2} + (a_{N-1} - (z+z^{-1})) \tilde y_{N-1}
   + \tilde y_N =0,
$$
from where
$$\tilde y_{N-2} = \frac{1}{b_{N-2}}((z+z^{-1}) z^{N-1} - a_N z^{N-1}
- z^N) = \frac{z^{N-2}}{c_{N-2}} (1 + z R_{N-2}(z)), $$
where the polynomial $R_{N-2}(z) = -a_N$, and $c_{N-2}=b_{N-2}$.
Thus, for $k=N-2$ the proposition is verified.

   Having the statement for $k$, let us prove it for $k-1,\  k \leq N-2$.
For $k\leq N-2$ the functions $\tilde y_k$
satisfy the equation (0.5), so
$$\tilde y_{k-1} = \frac{1}{b_{k-1}} \{
   [(z+z^{-1}) - a_k ] \frac{z^k}{c_k} (1 + z R_k(z)) \  - \
  \frac{b_k z^{k+1}}{c_{k+1}} (1 + z R_{k+1}(z)) \}  =
  \frac{z^{k-1}}{c_{k-1}} (1 + z R_{k-1}(z)) \},
$$
where $R_{k-1}(z)$ is a polynomial in $z$, as was to be proved.
\hfill\rule{0.5em}{0.5em}

\medskip

   Now we find the normallized Weyl solution $y_k = \alpha (z) \tilde y_k$,
comparing the result of the proposition with (1.5) for $k=0$ and $k=-1$:
$$ y_N^R(0,\lambda) = \alpha (z) \tilde y_0 = m_N^R P_0 - { Q_0 \over b_{-1} }, $$
$$ y_N^R(-1,\lambda) = \alpha (z) \tilde y_{-1}
   = m_N^R P_{-1} - { Q_{-1} \over b_{-1} }. $$
Taking into account the initial data for $P_k$ and $Q_k$ and the proposition
above, we rewrite it in the form
$$\alpha(z)\frac{1}{c_0}(1+zR_0(z))=m_N^R,$$
$$\alpha(z)\frac{z^{-1}}{c_{-1}}(1+zR_{-1}(z))= - {1 \over b_{-1}}.$$
Thus,
$$\alpha(z)= - { {c_{-1}z} \over { b_{-1} (1+zR_{-1}(z)) } },$$
$$m_N^R(z+z^{-1}) =
  - { {c_{-1}z}\over{b_{-1}c_0} }{ {1+zR_0(z)}\over{1+zR_{-1}(z)} }=
  - z {{1+zR_0(z)}\over{1+zR_{-1}(z)}}
$$
and $y^R_{N}(k,\lambda)$
satisfy for $k\leq N-2$ the equality
$$ y^R_{N}(k,\lambda)  = - { {c_{-1}z^{k+1}} \over {b_{-1}c_k} }
     { {1+zR_k(z)} \over {1+zR_{-1}(z)} } =
    - { h_k \over b_{-1} } z^{k+1} { {1+zR_k(z)} \over {1+zR_{-1}(z)} },
   \quad k \leq N-2,
$$
where $h_k = { {c_{-1}}\over{c_k} }$ are defined by (0.15), and
$z+z^{-1}=\lambda, \ |z|<1$.

\medskip

The functions of the form $\frac{1+zR_k(z)}{1+zR_{-1}(z)}$
are holomorphic and equal to $1$ at the point $z=0$. Thus, the last
equality implies the asymptotic formulas
$$  - \frac{b_{-1} z^{-(k+1)}}{h_k} y_{N}^R(k,z+z^{-1}) \to 1,
    \qquad z \to 0, \quad k\leq N-2.                  \eqno(1.6)
$$

\medskip

2. Let $J^R_N$ be the semi-infinite (to the right) matrix that is contained
in $J_N$ and starts with the element $a_0$.
We denote the spectral functions of the matrices $J^R$ and $J_N^R$
by $\rho^R(\tau)$ и $\rho^R_N(\tau)$. They are related with the Weyl
functions $m^R(\lambda)$, $m_N^R(\lambda)$ of the matrices $J^R$, $J_N^R$
by the transformation
$$m^R(\lambda)= \int_{-\infty}^{\infty} \frac{d\rho^R(\tau)}{\tau-\lambda},
  \qquad
  m_N^R(\lambda)= \int_{-\infty}^{\infty} \frac{d\rho_N^R(\tau)}{\tau-\lambda}.
  \eqno(1.7)
$$
Since the principal corner $N\times N$ minors of the matrices $J^R$ and
$J_N^R$ are the same, their spectral functions $\rho^R(\tau)$ and
$\rho^R_N(\tau)$ have the same first $2N$ moments (the integrals of the
functions $1,\tau,\tau^2,\ldots,\tau^{2N-1}$). Taking this into account,
after decomposing ${1\over\lambda({\tau\over\lambda}-1)}$ into the sum of
geometric progression, we estimate
$$|m^R(\lambda)-m^R_N(\lambda)|
  = | \int_{-\infty}^{\infty} \frac{d(\rho^R-\rho^R_N)(\tau)}{\tau-\lambda} |
$$
$$  = | \int_{-\infty}^{\infty} \frac{d(\rho^R-\rho^R_N)(\tau)}{\lambda(\frac{\tau}{\lambda}-1)} |
$$
$$ = | - \int_{-\infty}^{\infty} \frac{d(\rho^R-\rho^R_N)(\tau)}{\lambda}
     (1+ \frac{\tau}{\lambda} + \frac{\tau^2}{\lambda^2} + \ldots) |
$$
$$ = | - \int_{-\infty}^{\infty} \frac{d(\rho^R-\rho^R_N)(\tau)}{\lambda}
     (\frac{\tau^{2N}}{\lambda^{2N}} + \frac{\tau^{2N+1}}{\lambda^{2N+1}}  + \ldots) | $$
$$ = | \frac{1}{\lambda^{2N}} \int_{-\infty}^{\infty}
       \frac{\tau^{2N}d(\rho^R-\rho^R_N)(\tau)}{\tau-\lambda} |
$$
$$  \leq \frac{M_N(\lambda)}{|\lambda|^{2N}}, \eqno(1.8)
$$
where the function
$$ M_N(\lambda)= | \int_{-\infty}^{\infty}
        \frac{t^{2N} d(\rho^R-\rho^R_N)(\tau)}{\tau-\lambda} |
$$
is uniformly bounded outside a bounded (because of the boundedness
of $J^R$ and $J_N^R$) band along the real line.

\bigskip

3. It follows from (0.6) and (1.5) that
$$\varphi^R(k,\lambda) - y^R_{N}(k,\lambda)
   = (m^R(\lambda) - m^R_N(\lambda)) P_k(\lambda), \quad k\leq N-1.
$$
The polynomials $P_k(\lambda), \ |k|\leq N-1,$ have degrees less or equal
than $N-2$. So, $|P_k(\lambda)| \leq C_N |\lambda|^{N-2},\ |k|\leq N-1,$
where $C_N$ is a constant, depending on the coefficient of the polynomials
$P_k(\lambda), \ |k|\leq N-1$. Hence, taking into account estimate (1.8),
we have
$$ \max_{|k|\leq N-1} |\varphi^R(k,z+z^{-1}) - y^R_{N}(k,z+z^{-1})| =
   \max_{|k|\leq N-1} |(m^R(z+z^{-1}) - m^R_N(z+z^{-1})) P_k(z+z^{-1})|
$$
$$ \leq |z|^{2N}\frac{1}{|z|^{N-2}} K_N(z) = |z|^{N+2} K_N(z), \eqno(1.9)
$$
where $K_N(z)$ is a certain function, bounded in the neighborhood
of the point $z=0$.

\medskip

4. Comparing estimate (1.9) with asymptotics (1.6), we see that the functions
$\varphi^R(k,z+z^{-1})$ have the same asymptotic behavior:
$$ - \frac{b_{-1}z^{-(k+1)}}{h_k} \varphi^R(k,z+z^{-1}) \to 1,
                             \quad z \to 0. \eqno (1.10)
$$

\medskip

5. Analogously, it is easy to obtain the asymptotis of the solutions
$\varphi^L(k,\lambda)$. After using transformation (0.10), we have
for $\psi(k,z)$ asymptotic formulas (1.3) and (1.4).
Asymptotics (1.1) and (1.2) can be easily obtained, for example, from
the integral representation (0.8) of the Weyl function.

\hfill\rule{0.5em}{0.5em}

   We observe that for negative $k$ asymptotic formula (1.10)
for the increasing as $\lambda\to\infty$ solution $\varphi^R(k,\lambda)$
is a simple consequence of the well-known equalities for the leading
coefficients of the polynomials $Q_k(\lambda)$. At the same time,
it is not possible to obtain in this way formula (1.10) for positive $k$,
because the solutions $\varphi^R(k,\lambda)$, $k>0$, diminish as
$\lambda\to\infty$. The same situation takes place with
$\varphi^L(k,\lambda)$ for negative $k$. That is the reason for which
we had to prove theorem 1, which allowed us to obtain the asymptotics
(1.3), (1.4) for all the integer $k$.

    We also precise that the function $\psi(k,z)$ and $n(z)$ are
holomorphic at such points $z$ inside the unit disk for which $z+z^{-1}$
do not belong to the spectrum of the semi-infinite matrix $J^R$
(besides $z=0$) and at such points $z$ outside the unit disk for which
$z+z^{-1}$ do not belong to the spectrum of the semi-infinite
matrix $J^L$ (starting with the element $a_{-2}$). Further, the functions
$h_k z^{-(k+1)} \psi(k,z)$ and $\frac{n(z)}{z}$ are holomorphic in a certain
deleted neighborhood of the points $z=0, z=\infty$ and have a finite
limit at these points. So, they are holomorphic at the points
$z=0$ and $z=\infty$.

\bigskip

   In order to obtain factorization (0.25) with the necessary properties
of the function $R(z)$ we will need some auxiliary facts, which are
consequence of conditions A)--E). First, we remark that it follows
from the conditions C) and D) that in some neighborhood of the closed
set $\Omega^a_2\backslash [-2,2]$ we have the inequality
$$ -\pi<\eta^R(\alpha)-\eta^L(\alpha)< \pi.
$$
Hence, in condition D) we can constrict the intervals $\delta_l$ so that
this inequality holds under these intervals. Further, in condition
D) we can select a finite subcovering from the covering of the compact set
$\Omega^a_2\backslash [-2,2]$ with the open intervals $\delta_l$.
Thus, the condition D) can be replaced with a stronger one:

D$^\prime$) {\it The set $\Omega^a_2 \backslash [-2,2]$ can be covered
with a finite system of mutually disjoint intervals $\delta_l$,
under each of which the following inequalities are true: }
$$ -\pi<\eta^R(\alpha)-\eta^L(\alpha)< \pi \,,  \eqno(1.11)
$$
$$
{\rm ess}\sup_{\alpha\in\delta_l}\eta^R(\alpha)-
{\rm ess}\inf_{\alpha\in\delta_l}\eta^R(\alpha)<\pi\,, \eqno(1.12')
$$
$$
{\rm ess}\sup_{\alpha\in\delta_l}\eta^L(\alpha)-
{\rm ess}\inf_{\alpha\in\delta_l}\eta^L(\alpha)<\pi\,. \eqno(1.12'')
$$

\bigskip

   Let
$$ \tilde\Delta
  ={\bf R}\backslash (\tilde\Omega_1\cup \tilde\Omega_2^s \cup [-2,2])
  =\cup_k\tilde\Delta_k,
$$
where $\tilde\Delta_k=(\tilde\alpha_k,\tilde\beta_k)$ are mutually
disjoint intervals.

   Let us also define the function
$$ M(\lambda) = b_{-1} m^R(\lambda) - { 1 \over b_{-1} m^L(\lambda) }.
$$

   The two following lemmas clarify the behavior of the
arguments
$$\eta(\tau) := \arg M(\tau+i0),$$
$$\eta^R(\tau) := \arg m^R(\tau),$$
$$\eta^L(\tau) := \arg \Bigl( - { 1 \over b_{-1} m^L(\lambda) } \Bigr).$$
on the set
$\tilde\Delta \backslash\tilde\Omega_2^a$.

\medskip

{\bf Lemma 1.1.} {\it The arguments $\eta^R(\tau)$, $\eta(\tau)$ of the
functions $m^R(\tau+i0)$, $M(\tau+i0)$ take on the same constant value,
equal either to $0$ or $\pi$, on every set
$(\delta_l\cap\tilde\Delta)\backslash\tilde\Omega_2^a$. }

\medskip

P r o o f. \quad We see that the set $\tilde\Delta$
only include the part $\tilde\Omega_2^a$ of the set $\tilde\Omega$.
Remind, that the set $\tilde\Omega_2^a$ lies in the union of the
intervals $\delta_l$ (not to mix $\delta_l$ with $\tilde\Delta_k$!),
i.e. $\delta_l\backslash\tilde\Omega_2^a$ is not empty. Since the set
$(\delta_l\cap\tilde\Delta)\backslash\tilde\Omega_2^a$ lies in the complement
of the set $\tilde\Omega$, the functions $m^R(\lambda)$ and $M(\lambda)$
are holomorphic and real in the points of this sets, and their arguments
can only take value $0$ or $\pi$.
Condition D$'$) implies that the ocsillation of the arguments $\eta^R(\alpha)$
of $m^R(\alpha+i0)$ on the interval $\delta_l$ is less, then $\pi$;
and it follows that the range of the function $\eta^R(\alpha)$ on the set
$(\delta_l\cap\tilde\Delta)\backslash\tilde\Omega_2^a$ cannot contain both
$0$ and $\pi$. Therefore, the function $\eta^R(\alpha)$ is constant, equal
to $0$ or $\pi$ on the whole set
$(\delta_l\cap\tilde\Delta)\backslash\tilde\Omega_2^a$.
Further, since the values of $\eta^R(\alpha)-\eta^L(\alpha)$ on this set
are multiples of $\pi$ and because of the condition D$'$), the inequalities
(1.11) hold there, we have
$$ \eta^R(\alpha)-\eta^L(\alpha)=0 ,\qquad
   \alpha\in(\delta_l\cap\tilde\Delta)\backslash\tilde\Omega_2^a. \eqno(1.13)
$$
According to the definition of the function $M(\lambda)$,
$$ M(\lambda)= m^R(\lambda) \Bigl( b_{-1} +
   {-1\over b_{-1}) m^R(\lambda) m^L(\lambda) } \Bigr),
$$
So,
$$ \arg M(\tau+i0)=\arg m^R(\tau+i0)+
    \arg \Bigl( b_{-1} + {-1\over b_{-1}) m^R(\lambda) m^L(\lambda) } \Bigr),
$$
But the function ${-1\over b_{-1}) m^R(\lambda) m^L(\lambda) }$
is regular and real on the set
$(\delta_l\cap\tilde\Delta)\backslash\tilde\Omega_2^a$; its argument equals
$\eta^L(\tau) - \eta^R(\tau)$ on this set. On the other hand, we have from
(1.13) that ${-1\over b_{-1}) m^R(\lambda) m^L(\lambda) } > 0 $ for
$P(\tau,\eta^L - \eta^R) > 0$, and
$ b_{-1} + {-1\over b_{-1}) m^R(\lambda) m^L(\lambda) } >0$, too.
Hence, $\eta^R(\tau)=\eta(\tau)$ on this set.
\hfill\rule{0.5em}{0.5em}

\medskip

{\bf Lemma 1.2.} {\it Each interval
$(\tilde\alpha_k,\tilde\beta_k)=\tilde\Delta_k$
splits into two pieces
$$ \tilde\Delta_k^-=(\tilde\alpha_k,\tilde\varphi_k),
    \qquad  \tilde\Delta_k^+=(\tilde\varphi_k,\tilde\beta_k), \eqno (1.14)
$$
so that
$$ M(\alpha)\ <\ 0,\quad  \alpha\in\tilde\Delta_k^-\backslash\tilde\Omega_2^a,
    \eqno(1.15')
$$
$$  M(\alpha)\ >\ 0\quad  \alpha\in\tilde\Delta_k^+\backslash\tilde\Omega_2^a.
    \eqno(1.15'')
$$
(One of the intervals $\tilde\Delta_k^-,\tilde\Delta_k^+$ may be empty.)}

 \medskip

   Remind that $M(\lambda)$ is holomorphic on
$\tilde\Delta\backslash\tilde\Omega_2^a$; so we can write
$ M(\alpha),\ \alpha\in\tilde\Delta_k^\pm\backslash\tilde\Omega_2^a$.

\medskip

P r o o f.\quad If the interval $\tilde\Delta_k$ contains no point of
$\tilde\Omega_2^a$, then it lies in the complement of
$\tilde\Omega=\tilde\Omega^R\cup \tilde\Omega^L$ and the function
$M(\alpha)$ is real and increasing. Hence in this interval there is at most
one point, where $M(\alpha)$ changes its sign from $-$ to $+$.
Let us denote this point by $\tilde\varphi_k$ (if $M(\alpha)$ is positive,
resp., negative, on the whole interval, we set
$\tilde\varphi_k=\tilde\alpha_k$ or $\tilde\beta_k$), and
obtain the intervals (1.14), on which the inequalities (1.15) hold.

   If the interval $\tilde\Delta_k$ contains points of $\tilde\Omega_2^a$,
then, according to condition D$^\prime$), the set
$\tilde\Delta_k\cap\tilde\Omega_2^a$ is covered by a finite number of
mutually disjoint intervals $\delta_l\cap\tilde\Delta_k=(\alpha^{(l)}_k,
\beta^{(l)}_k)$, on which the inequalities (1.11), (1.12) are fulfilled.
Obviously, these intervals can be labelled so that
$$\tilde\alpha_k=\beta^{(0)}_k\leq\alpha^{(1)}_k\ < \ \beta^{(1)}_k\leq\alpha^{(2)}_k
<\beta^{(2)}_k\leq\dots\\
\dots \leq\alpha^{(n)}_k<\beta^{(n)}_k\leq\alpha^{(n+1)}_k=\tilde\beta_k
$$
(for the sake of convenience the endpoints $\tilde\alpha_k,\tilde\beta_k$
of $\tilde\Delta_k$ are denoted by $\beta^{(0)}_k,\alpha^{(n+1)}_k$).
According to Lemma 1.1 the argument $\eta(\alpha)$ of the function
$M(\alpha+i0)$ is constant, equal to $0$ or to $\pi$, on each set

$$
\Phi_p=(\alpha^{(p)}_k,\beta^{(p)}_k)\backslash\tilde\Omega^a_2\,, \qquad
 1\leq p\leq n.
$$
Therefore the function $M(\alpha+i0)=M(\alpha)$ is real, regular, of constant
sign, on each set $\Phi_p$. Let us denote by $\Phi_{p^+}$
(resp.\ $\Phi_{p^-}$) the set with the smallest (resp. greatest) number
$p^+$ (resp.\ $p^-$), at which the function $M(\alpha)$ is positive
(resp. negative) and $\Phi_0=(\beta^{(0)}_k,\alpha^{(1)}_k)=
(\tilde\alpha_k,\alpha^{(1)}_k)$ (resp.\ $\Phi_{n+1}=(\beta^{(n)}_k,
\alpha^{(n+1)}_k)=(\beta^{(n)}_k,\tilde\beta_k)$), if this function is
positive (resp.\ negative) on the set $\Phi_p\quad(1\leq p\leq n)$.

   Since the function $M(\alpha)$ is regular and grows monotonically on the
segment
$[\beta^{(p^+)}_k,\alpha^{(p^++1)}_k]$
(resp. $[\beta^{(p^--1)}_k,\alpha^{(p^-)}_k]$) and in its neighborhood,
it remains positive (negative) on this segment and on the next (previous)
set $\Phi_{p^++1}$ (resp. $\Phi_{p^--1})$. Hence, $p^+=p^-+1$
and the function $M(\alpha)$ is positive on the set
$(\alpha^{(p^+)}_k,\tilde\beta_k) \backslash\tilde\Omega_2^a$ and negative
on the set $(\tilde\alpha_k,\beta^{(p^+-1)}_k) \backslash\tilde\Omega^a_2$.
Due to the monotonicity on the segment
$[\beta^{(p^+-1)}_k,\alpha^{(p^+)}_k]$ there is only one point
$\tilde\varphi_k$, at which the function $M(\alpha)$ changes sign
from $-$ to $+$. Setting $\tilde\varphi_k=\tilde\alpha_k$, if $p^+=0$
($\tilde\varphi_k=\tilde\beta_k$, if $p^-=n+1$), we obtain
the intervals (1.14), on which (1.15) are fulfilled.
\hfill\rule{0.5em}{0.5em}

\bigskip
\bigskip

\centerline{\bf 2. Factorization and choice of the function $R(z)$}

\bigskip

   {\bf 1.} \ Let us define the function
$$N(z)=n(z)-n(z^{-1}).
$$
This function is connected with the resolvent matrix $R(m,n,\lambda)$
of the matrix $J$ in the following way (see, for example, [10]):
the two diagonal elements of this matrix are equal to
$$ R(-1,-1,\lambda) \mid_{\lambda=z+z^{-1}} = R(-1,-1,z+z^{-1})
   = { 1 \over b_{-1} N(z) },  \eqno(2.1)
$$
$$ R(0,0,\lambda) \mid_{\lambda=z+z^{-1}} = R(0,0,z+z^{-1}) =
   { n(z) n(z^{-1}) \over b_{-1} N(z) }.  \eqno(2.2)
$$

   The goal of this section is to factorize in a suitable way the function
$\frac{z-z^{-1}}{N(z)}$ in form (0.25). Such a factorization is presented
in more general case in [15].

\medskip

   The functions $b_{-1} m^R(\lambda)$ and ${-1 \over b_{-1} m^L(\lambda)}$
are holomorphic in the upper half-plane and have positive imaginary part
in it. Thus, (see, e.g., [14]), taking into account their asymptotic
behavior (see (1.1), (1.2)), we can represent these functions
in the form:
$$ b_{-1} m^R(\lambda) = \int_{-\infty}^{\infty}
   { {d\rho_R(\tau)} \over {\tau-\lambda} }, \qquad
  - { 1 \over {b_{-1} m^L(\lambda)} } = {\lambda\over b_{-1}} + \beta +
   \int_{-\infty}^{\infty} { {d\rho_L(\tau)} \over {\tau-\lambda} }, \eqno(2.3)
$$
where, according to the Stieltjes-Perron inversion formula,
$$ \rho_R (\lambda) = {1\over\pi} \lim_{\varepsilon \downarrow 0}
      \int\limits_{-\infty}^\lambda {\rm Im} \,
      \bigl( b_{-1} m^R(\tau+ i \varepsilon) \bigr)       d \tau,
   \qquad
   \rho_L (\lambda) = {1\over\pi} \lim_{\varepsilon \downarrow 0}
      \int\limits_{-\infty}^\lambda {\rm Im} \,
      { - 1 \over b_{-1} m^L(\tau+ i \varepsilon) } d \tau, \eqno(2.4)
$$
the number $\beta \in {\bf R}$, and $\rho_R(\lambda)$ and $\rho_L(\lambda)$
are constant in the infinity (the Weyl functions of the matrices $J^R$ and
$J^L$ are regular out of their spectrum, which is bounded because of the
boundedness of the matrices).

\medskip

On the parts of the real line where the two functions
$\rho_R(\tau)$ and $\rho_L(\tau)$ are absolutely continuous
and $\rho'_R(\tau)$ and $\rho'_L(\tau)$ do not vanish
(it is shown in the end of this section that it is the set $\tilde\Omega_2^a$),
we introduce
$$ \mu(\tau) :=  \sqrt{ {\rho_R}^{\prime} (\tau) \over
                        {\rho_L}^{\prime} (\tau) } \,.  \eqno(2.5)
$$

\bigskip

 {\bf 2.} \ {\it Notation 1.}  Let us define on the real line the
involutive map $V$ of symmetry with respect to the unit circle:
$$ V(t)=t^{-1}, \quad t \in {\bf R}\backslash \{0\}. \eqno(2.6)
$$
For a set $A \subset {\bf R}\backslash \{0\}$ and for a function $\rho(t)$,
defined on ${\bf R}\backslash \{0\}$, we define
$$ V(A)=\{ t \ \vert \ t^{-1}\in A \},
$$
$$   V(\rho)(t)=\rho(t^{-1}).
$$

\medskip

   {\it Notation 2.} Let us introduce
$$P(z,\gamma)\ =\ \exp \Bigl\{ { 1\over{\pi} }
   \int\limits_{-\infty}^{\infty} \gamma(t)
   ( {1\over{t-z}} - {t\over{1+t^2}} )dt \Bigr\},   \eqno(2.7)
$$
$$ \hat P(z,\hat\gamma)\ =\
   \exp\biggl\{ - { 1 \over {2\pi i} }{\int\limits_{-\pi}^\pi}
    { {e^{i\theta}+z} \over {e^{i\theta}-z} }\
    \hat\gamma(\theta)d\theta\biggr\},   \eqno(2.8)
$$
where $\gamma(t)=\overline{\gamma(t)}, \ -\infty<t<\infty,$ and
$\hat\gamma(\theta)=\overline{\hat\gamma(\theta)}, \ -\pi<t<\pi,$
are bounded and measurable functions. The function $P(z,\gamma)$ is defined
and holomorphic, at least, for nonreal $z$, and $\hat P(z,\hat\gamma)$ is
defined and holomorphic inside and outside the unit disk.
It is easy to obtain from the formulas of Plemelj-Sokhotsky the following
equalities for the limit values of the functions $P(z,\gamma)$ and
$\hat P(z,\hat\gamma)$ on the real line and the unit circle, resp.:
$$\arg P^+(t,\gamma) = - \arg P^-(t,\gamma)= \gamma(t),
   \qquad -\infty< t < \infty,   \eqno(2.9)
$$
$$ |P^+(t,\gamma)| = |P^-(t,\gamma)|,
   \qquad -\infty< t < \infty,   \eqno(2.10)
$$
$$ \arg \hat P^+(e^{i\theta},\hat\gamma)
   = - \arg \hat P^-(e^{i\theta},\hat\gamma)
   = \hat\gamma(\theta),
   \qquad -\pi<\theta<\pi,    \eqno(2.11)
$$
$$ |\hat P^+(e^{i\theta},\hat\gamma)|
        = |\hat P^-(e^{i\theta},\hat\gamma)|,
   \qquad -\pi<\theta<\pi.    \eqno(2.12)
$$

   Besides, it is easy to verify that the following lemma is true:

\medskip

{\bf Lemma 2.1.} \ {\it The functions $P(z,\gamma)$ and $\hat P(z,\hat\gamma)$
in their domain of holomorphy satisfy the equalities:
$$P(z,\gamma_1+\gamma_2)=P(z,\gamma_1) P(z,\gamma_2),    \eqno(2.13)
$$
$$ P(\overline{z},\gamma)
   = \overline{P(z,\gamma)},  \eqno(2.14)
$$
$$ P(z^{-1},\gamma)
   = P(z,-V(\gamma)),      \eqno(2.15)
$$
$$\hat P(z,\hat\gamma_1+\hat\gamma_2)
   =\hat P(z,\hat\gamma_1) \hat P(z,\hat\gamma_2),  \eqno(2.16)
$$
If the function $\hat\gamma(\theta)=-\hat\gamma(-\theta)$ is odd, then}
$$ \hat P(z^{-1},\hat\gamma) = \hat P(z,\hat\gamma).  \eqno(2.17)
$$

\medskip

   The function $M(\lambda)$,
being a function with positive imaginary part in the upper half-plane,
can be represented in the multiplicative form
$$ M(\lambda) = b_{-1} m^R(\lambda) - { 1 \over b_{-1} m^L(\lambda) }
   = C_1 P(\lambda, \eta), \eqno(2.18)
$$
where
$$ \eta(\tau)= \lim_{\varepsilon \downarrow 0} \arg M(\tau+i\varepsilon)
\eqno(2.19)
$$
and the constant $C_1 > 0$. (It is easy to obtain this representation
by considering the logarithm of these functions, which also have
positive imaginary part in the upper half-plane, and then by representating
the logarithm in the integral form of the type (2.3)\,
(see, for example, [16]).)

\bigskip

   {\bf 3.} \ We remind that we defined in the plane of the spectral
parameter $\lambda$ the sets $\tilde\Omega^R$, $\tilde\Omega^L$ as the
supports of the measures $d\rho_R(\tau)$, $d\rho_L(\tau)$. This means that
$\tilde\Omega^R$, $\tilde\Omega^L$ are the sets on which the functions
$m^R(z)$, ${-1\over m^L(z)}$ have singularities, i.e.~are not regular.
Then we introduced
$$ \tilde\Omega_2\ = \ \tilde\Omega^R\cap \tilde\Omega^L,  \qquad
   \tilde\Omega_1\ = \ (\tilde\Omega^R\backslash \tilde\Omega^L) \cup
                 (\tilde\Omega^L\backslash \tilde\Omega^R)
$$
and divided $\tilde\Omega_2$ into
$\tilde\Omega_2^s$ and $\tilde\Omega_2^a$.
Here $\tilde\Omega_2^s$ is the finite or empty set of the common poles
of the functions $m^R(z)$ and ${-1\over m^L(z)}$
($\tilde\Omega_2^s \subset \tilde\Omega_2$), and
$[-2,2] \subset \tilde\Omega_2^a  = \tilde\Omega_2 \backslash
\tilde\Omega_2^s$. It will be seen from lemma 2.3 that the set
$\tilde\Omega_2^a$ (including the segment $[-2,2]$) belongs to the
absolutely continuous spectrum of multiplicity $2$ of the semi-infinite
matrices $J^R$ and $J^L$.

   Now we are going to find the equivalent of these sets in the $z$-plane
($z+z^{-1}=\lambda$). Let us define on the real line in the $z$-plane
$$ \Omega \ = \  \{t \in {\bf R} \mid |t|<1, \ t+t^{-1} \in \tilde\Omega^R \}
          \cup \{t \in {\bf R} \mid |t|>1, \ t+t^{-1} \in \tilde\Omega^L \}\,.
$$
According to this definition, $\Omega$ is the set of the points on the real
line (except $\pm 1$) at which the function $n(z)$ is not holomorphic.
Besides the set $\Omega$, the function $n(z)$ has singularities on the
unit circle $\bf T$.

   Let us divide $\Omega$ into its symmetric and asymmetric (with respect
to the unit circle) parts:
$$ \Omega_2\ = \ \Omega\cap V(\Omega) = \{t \in {\bf R} \mid
    t\in\Omega, \ t^{-1}\in\Omega\}\,,
$$
$$ \Omega_1\ = \ \Omega\backslash V(\Omega)
   =\{t \in {\bf R} \mid t\in\Omega, \ t^{-1}\not\in\Omega\}\,,
$$
where the map $V$ is defined by formula (2.6). It may be seen that these
sets are the preimages of the sets $\tilde\Omega_2$ and
$\tilde\Omega_1$ with respect to the map $z+z^{-1}=\lambda$:
$$ \Omega_2 \ = \  \{t \in {\bf R} \mid t \ne \pm 1, \ t+t^{-1} \in \tilde\Omega_2 \},
$$
$$ \Omega_1 \ = \  \{t \in \Omega \mid t+t^{-1} \in \tilde\Omega_1 \} \,.
$$
We also define the sets $\Omega^s_2$ and $\Omega^a_2$ as the preimages
of $\tilde\Omega^s_2$ and $\tilde\Omega^a_2$ with respect to the
map $\lambda=z+z^{-1}$:
 $$ \Omega_2^s \ = \  \{t \in {\bf R} \mid t+t^{-1} \in \tilde\Omega_2^s \}, \qquad
   \Omega_2^a \ = \  \{t \in {\bf R} \mid t\ne\pm 1, \ t+t^{-1} \in \tilde\Omega_2^a \}.
$$

   Thus, the set $\Omega$ is divided into three sets:
$$ \Omega=\Omega_1\cup\Omega^s_2\cup\Omega^a_2.
$$
It follows from condition A) that the three sets have positive mutual
distances and that $\Omega^s_2$ is the finite or empty set of the common
poles of $n(z)$ and $n(z^{-1})$ (they are all simple). Since
$\Omega_1$ is the asymmetric part of the set $\Omega$, the function $n(z)$
is holomorphic on the set $V(\Omega_1)$.

Let
$$
\Delta={\bf R}\backslash
       ((-1,1)\cup \Omega_1\cup V(\Omega_1) \cup\Omega_2^s)=\cup_k\Delta_k,
   \eqno(2.20)
$$
where $\Delta_k=(\alpha_k,\beta_k)$ are mutually disjoint intervals.
(The set $\Omega_2^a$, which corresponds to the absolutely continuous spectrum
except the segment $[-2,2]$, lies in the intervals $\Delta_k$.)
According to the definition, all $|\alpha_k|\ge 1$ and $|\beta_k|\ge 1$.

   The following fact is an immediate consequence of lemma 1.2:

\medskip

{\bf Lemma 2.2.} {\it Each interval $(\alpha_k,\beta_k)=\Delta_k$
splits into two pieces:
$$ \Delta_k^{(1)}=(\alpha_k,\varphi_k),
   \qquad\Delta_k^{(2)}=(\varphi_k,\beta_k), \eqno(2.21)
$$
so that
$$ N(\alpha)\ < \ 0\quad
   \mbox{для $\alpha\in\Delta_k^{(1)}\backslash\Omega_2^a$}\,,
   \eqno(2.22')
$$
$$ N(\alpha)\ > \ 0\quad
   \mbox{для $\alpha\in\Delta_k^{(2)}\backslash\Omega_2^a$}\,.
   \eqno(2.22'')
$$
(One of the intervals $\Delta_k^{(1)},\Delta_k^{(2)}$ may be empty.) }

\medskip

   The points $\tilde\alpha_k = \alpha_k + \alpha_k^{-1}$ and
$\tilde\varphi_k = \varphi_k + \varphi_k^{-1}$
from the plane of the spectral parameter correspond to the pairs of points
$\{\alpha_k,\alpha_k^{-1}\}$ and $\{\varphi_k,\varphi_k^{-1}\}$
in $z$-plane.
To make the factorization with the properties we will need later on,
we have to choose one number $\alpha_k^*$ (resp. $\varphi_k^*$)
in each pair, according to the following rule: put
$$ \alpha^*_k=\cases{\alpha_k, {\ \rm if\ }
                              \alpha_k\not\in V(\Omega_1),\cr
                     \alpha_k^{-1},
                              {\ \rm if\ } \alpha_k\in V(\Omega_1).\cr}
   \eqno(2.23)
$$
If both endpoints of an interval $(\alpha_k,\beta_k)$ belong to the set
$\Omega_1$ or $V(\Omega_1)$, we set
$$
\varphi_k^*= \cases{\varphi_k, {\ \rm if\ }
                               \alpha_k, \beta_k \in \Omega_1,\cr
                \varphi_k^{-1},
                         {\ \rm if\ } \alpha_k, \beta_k \in V(\Omega_1),\cr}
   \eqno(2.24')
$$
and in all other cases (there are only finitely many, which is
a corollary of the disjointedness of $\Omega_1$ and $V(\Omega_1)$) we set
$$ \varphi_k^*= \cases{\varphi_k, {\ \rm if\ }
                             \varphi_k\not\in V(\Omega_1),\cr
      \varphi_k^{-1}, {\ \rm if\ }
                              \varphi_k\in V(\Omega_1).\cr}
    \eqno(2.24'')
$$

   Remark that for some $k,j$ \ $\alpha_k=-\infty$ and
(because of the condition E)) $\beta_j=\varphi_j=-1$.
We will suppose that the intervals
$\Delta_k$ are numbered so that
$\Delta_0^{(1)}=(-\infty, \varphi_0)$ and
$\Delta_1 = (\alpha_1,-1)$. (If $\beta_0=-1$, we put
$\alpha_1^*=\varphi_0^*=\varphi_1^*=-1$, and than the final formula
also cover this case.)

\medskip

Let
$$ \Delta^{(1)}=\cup_k\Delta^{(1)}_k.
$$

We denote by $\Phi=\{\varphi_k^*\}$
the set of all $\varphi_k^*$, except the point $-1$.
We also denote by $\chi_1( t),\chi_2^a( t)$ and $\chi_0( t)$
the indicators of the sets $\Omega_1$, $\Omega_2^a$
and $\Delta^{(1)}\cup V(\Delta^{(1)})$.

\medskip

The key point of the inverse problem is presented by

\bigskip

{\bf Theorem 2 \ (factorization).} \ {\it The function
${(z-z^{-1})} {N(z)}^{-1}={(z-z^{-1})} {(n(z)-n(z^{-1}))}^{-1}$
in its domain of holomorphy is represented in the form of the product of
two functions, which are holomorphic in this domain:
$$ { {z-z^{-1}} \over {N(z)} } \ = \  R(z) R(z^{-1}),
  \eqno(2.25)
$$
$$  R(z)\ = \  C R_0(z)R_1(z)R_2(z)R_3(z)R_\mu(z), \eqno(2.26)
$$
where
$$    R_0(z)\ = \ { {z-\alpha_1^*} \over {z-\varphi_0^*} }
       \prod_{k \ge 2}  { {z-\alpha_k^*} \over {z-\varphi_k^*} } \,,
                \eqno(2.27)
$$
$$  R_1(z)\ = \ P(z,-\gamma_1),\qquad R_2(z)=P(z,-\gamma_2)\,,
                \eqno(2.28)
$$
$$  \gamma_1(t)\ = \ \chi_1(t) s(t) \eta(t+t^{-1})\,,
                \eqno(2.29)
$$
$$ \gamma_2(t)\ = \ {1\over2}\chi_2^a(t) s(t) (\eta(t+t^{-1})-\chi_0(t)\pi),
                \eqno(2.30)
$$
$$ s(t) = \cases{1, \quad |t|>1, \cr
                  -1, \quad |t|<1, \cr}
    \eqno(2.31)
$$
the constant $C>0$, and the numbers $\alpha_k^*$ and $\varphi_k^*$
are defined by equalities $(2.23)$, $(2.24')$ and $(2.24'')$
(with $\eta(\tau)$, defined by formula (2.19)),
$$
  R_3(z) = \hat P(z,-{\hat\gamma_3 \over 2}), \eqno(2.32)
$$
$$\hat\gamma_3 (\theta)=
    \cases{-\eta( e^{i\theta} + e^{-i\theta} )
                               + { \pi \over 2 }, \quad 0<\theta<\pi, \cr
            \eta( e^{i\theta} + e^{-i\theta} )
                               - { \pi \over 2 }, \quad -\pi<\theta<0, \cr}
   \eqno(2.33)
$$
$$ R_\mu(z) = \exp \biggl\{ - { 1 \over 2 \pi } \int\limits_{-\pi}^\pi
    { e^{i\alpha} \over e^{i\alpha} - z }
    \ln \mu (e^{i\alpha}+e^{-i\alpha}) \  d\alpha \biggr\} \eqno(2.34)
$$
(with $\mu(\tau)$, defined by formula (2.5)).

   In addition, the limit values of the function $R(z)$ on the unit circle
${\bf T}$ satisfy the equality
$$ R^+(\xi) R^+(\xi^{-1}) \, |{\rm Im}\, n^+(\xi)| \ = \
   R^-(\xi) R^-(\xi^{-1}) \, |{\rm Im}\, n^-(\xi)| \,,
   \quad |\xi|=1.  \eqno(2.35)
$$

The function $R(z)$ may have singularities only at such points $z$, that
$z+z^{-1}$ belongs to the spectrum of the matrix $J$. If $z+z^{-1}$
belongs to $\Omega_2^a$, then $R(z)$ may have singularities at both
points $z, z^{-1}$, and if it belongs to the simple spectrum, then it only
have a singularity in one of these points.
Moreover, the functions $R_0(z)$, $R_1(z)$, $R_2(z)$ and $R_3(z) R_\mu(z)$
may have singularities {\it only} on the sets $\Phi=\{\varphi_k^*\}$,
$\Omega_1$, $\Omega_2^a$ and ${\bf T}$, resp. The sets
$$ \Omega_1\cup\Phi,\quad V(\Omega_1)\cup
V(\Phi),\quad\Omega_2^s,\quad\Omega_2^a,\quad {\bf T} $$
are at positive distances from one another. }

\medskip

   P r o o f. \
   The holomorphy of the functions $R_0(z)$, $R_1(z)$, $R_2(z)$ and
$R_3(z) R_\mu(z)$ outside the sets $\Phi=\{\varphi_k^*\}$,
$\Omega_1$, $\Omega_2^a$ and ${\bf T}$ is an immediate consequence of their
definition and properties (2.9)--(2.12). From property $(2.12)$
we have $|R_{0123}^+(e^{i\theta})| = | R_{0123}^-(e^{-i\theta})|$, where
$$  R_{0123}(z)\ := \  R_0(z)R_1(z)R_2(z)R_3(z).
$$
So, in order to show (2.35), we only have to demonstrate that
$$ R_\mu^+(\xi) R_\mu^+(\xi^{-1}) \, \rho_R'(\xi+\xi^{-1}) = \
   R_\mu^-(\xi) R_\mu^-(\xi^{-1}) \, \rho_L'(\xi+\xi^{-1}) \,,
   \quad|\xi|=1.
$$
(We remind that, according to (2.4), taking into account the absolute
continuity of $\rho_R(\tau)$, $\rho_L(\tau)$ on the set
$[-2,2]\subset\Omega_2^a$, we have
$\rho'_R(\tau)={1\over\pi} \, {\rm Im} \, m^R(\tau+i0)$,
$\rho'_L(\tau)={1\over\pi} \, {\rm Im} \, {-1\over m^L(\tau+i0)}$.)
In order to prove the last equality it is sufficient to show that
$$ { R_\mu^-(e^{i\theta}) \over R_\mu^+(e^{i\theta}) }
    = \mu(e^{i\theta}+e^{-i\theta}) =
    \sqrt{ \rho'_R (e^{i\theta}+e^{-i\theta}) \over
           \rho'_L (e^{i\theta}+e^{-i\theta}) } ,
$$
which is a simple implication of the definition of $R_\mu(z)$ and
the Plemelj-Sokhotsky formulas. So, (2.35) is proved. Further,
taking into account the evennes of the ratio
$\sqrt{ \rho'_R (e^{i\theta}+e^{-i\theta}) \over
  \rho'_L (e^{i\theta}+e^{-i\theta}) }$ \
(as a function on $\theta, \ -\pi<\theta<\pi$),
it is easy to verify from the definition of $R_\mu(z)$,
after some elementary transformations, that
$$R_\mu(z)R_\mu(z^{-1})= {\rm const} >0, \quad |z| \ne 1.
$$
Thus, to prove (2.25), it suffice to prove that
in the domain if holomorphy of the function
${ {z-z^{-1}} \over {N(z)} }$
the following equality holds
$$  { {z-z^{-1}} \over {N(z)} } \ = \  С_1 R_{0123}(z) R_{0123}(z^{-1}),
     \quad  C_1>0.
$$

   Let us now consider the functions ${ {z-z^{-1}} \over {N(z)} }$
and $R_{0123}(z) R_{0123}(z^{-1})$. They are holomorphic and do not vanish
outside the unit circle $\bf T$ and a certain set of the real line,
which is separated from the circle, zero and the infinity. (We concluded
this from the asymptotics of $n(z)$ and the definition of $R_i(z)$.)
Moreover, both functions are holomorphic and positive at the infinity
and take on positive values under the real line in a neighborhood of zero.
Thus, both functions are holomorphic and do not vanish in two
simply connected domains on the Riemann sphere, one of which
includes the nonreal points inside the unit disk and a neighborhood of zero,
and the other one includes the nonreal points outside the disk and a
neighborhood of the infinity. Hence, in each one of the two simply connected
components we can take the logarithm of these function
(so that on the real line in the neighborhood of zero and the infinity, where
these functions are positive, it is the principal branch).
Let us denote
$$ f_1(z) = \ln \biggl( { {z-z^{-1}} \over {N(z)} } \biggr),
   \qquad
   f_2(z) = \ln \biggl( R_{0123}(z) R_{0123}(z^{-1}) \biggr)
$$
and prove that the limit values of ${\rm Im}\, f_1(z)$ и
${\rm Im}\, f_2(z)$ (i.e.~the arguments of ${ {z-z^{-1}} \over {N(z)} }$
and $R_{0123}(z) R_{0123}(z^{-1})$) are the same almost everywhere on the
unit circle and the real line.

   Generally speaking, the equality of the limits of
${\rm Im} \, f_1^\pm(z)$ and ${\rm Im} \, f_2^\pm(z)$ almost everywhere on
the boundary of the domain does not imply that
$f_1(z) = f_2(z) + \alpha, \ \alpha \in {\bf R}$.
However, in our case, as it is shown in paper [15], the boundedness of
these imaginary parts will allow us to use the maximum principle
and conclude that $f_1(z)$ and $f_2(z)$ differ by a real constant.
Hence, $ { {z-z^{-1}} \over {N(z)} }$ and
$R_{012}(z) R_{012}(z^{-1})$ differ by a constant positive factor
(it is evident that this factor is the same inside and outside the circle).

   The last speculation is not completely rigorous.
In paper [15] the factorization theorem is rigorously proved in another
way and in more general form. So, now we will restrict ourselves by
proving the equality of ${\rm Im}\, f_1(z)$ and ${\rm Im}\, f_2(z)$
on the boundary of the domain of their holomorphy, omitting the consideration
of some special cases.

   Begin with the circle. First, the function
$$R_{012}(z)=R_0(z)R_1(z)R_2(z)$$
is regular and does not vanish on the circle. The expression
$$R_{012}(\xi) R_{012}(\xi^{-1})=
  R_{012}(\xi) R_{012}(\overline{\xi})=
  R_{012}(\xi) \overline{R_{012}(\xi)}, \quad |\xi|=1,
$$
is always real, its argument is precisely equal to zero (not to $2 k \pi $).
Really, for nonreal $z$ outside the unit circle
$R_{012}(\overline{z})=\overline{R_{012}(z)}$ because of $(2.14)$.
Linking with a curve an arbitrary point of the circle with, for example,
a big positive $x$, where $R_{012}(x)R_{012}(x)>0$, we make sure that
at any point of this curve $R_{012}(\overline{z}) R_{012}(z)>0$. Hence,
on the unit circle we only have to show the equality of the arguments of
the functions ${ {z-z^{-1}} \over {N(z)} }$ and $R_3(z) R_3(z^{-1})$.
But, as it may be easily seen,
$$ \lim_{r\uparrow 1}
       \arg \Bigl( { {z-z^{-1}} \over {N(z)} } \Bigr) |_{z=re^{i\theta}}
   = - \lim_{r\downarrow 1}
       \arg \Bigl( { {z-z^{-1}} \over {N(z)} } \Bigr) |_{z=re^{i\theta}}
   = - \hat\gamma_3 (\theta).
$$
On the other hand, taking into account definition
$(2.32)$, property (2.11) and the oddness of $\hat\gamma_3(\theta)$,
we see that the limit values of $\arg R_3(z)$ on the unit circle equal
$$ \arg R_3^+(e^{i\theta}) = - \arg R_3^-(e^{i\theta}) =
   - \arg R_3^+(e^{-i\theta}) = \arg R_3^-(e^{-i\theta})
   = - { \hat\gamma_3 \over 2 },
$$
that is
$$ \lim_{r\uparrow 1}
   \arg (R_3(z) R_3(z^{-1})) |_{z=re^{i\theta}}
   = \arg \bigl(R^+(e^{i\theta})R^-(e^{-i\theta})\bigr)
   = - \hat\gamma_3 (\theta),
$$
which proves the wanted equality of the arguments on the circle.

   Let us now consider the real line. Here we prove the equality
of the arguments only under the interval $(1,+\infty)$
(the proof for the intervals $(-\infty,1)$, $(-1,0)$ and $(0,1)$
is analogous). Without loss of generality, we only consider
the limit values from above of the real line. The function $R_3(z)$
is holomorphic outside the circle, and for the real $x$ we have
$R_3(x)R_3(x^{-1})>0$. That is why, instead of $R_{0123}(z)$, we
only consider $R_{012}(z)$. First, under the interval $(1,+\infty)$
$$ \arg N^+(x) = \arg M^+(x+x^{-1}) = \eta(x+x^{-1}),
$$
$$ \arg \Bigl( { {(x+i0)-(x+i0)^{-1}} \over {N(x+i0)} } \Bigr)
   = - \eta(x+x^{-1}) .  \eqno(2.36)
$$

   Let us now consider the function $R_0(z)R_0(z^{-1})$.
For the sake of simplicity we will asume that
$\alpha_1^* = \varphi_0^* = -1$, and that for $k\ge 2$ all
$\alpha_k^* > 0$, $\varphi_k^* >0$, and also either
$0< \alpha_k^*, \varphi_k^* < 1$, or
$1< \alpha_k^*, \varphi_k^*$.
(As it was already said there are finitely many cases when
$\alpha_k^*, \varphi_k^*$ lie on different sides of 1.)
In the first case $0< \varphi_k^* < \alpha_k^* <1 $ and
$$ \arg \Bigl( {x-i0 - \alpha_k^* \over x - i0 - \varphi_k^*} \Bigr)
   = \cases { 0, \quad x < \varphi_k^* \ {\rm or} \ x > \alpha_k^*, \cr
              - \pi, \quad \varphi_k^* < x < \alpha_k^*. \cr }
$$
In the second case $1< \alpha_k^* < \varphi_k^*$ and
$$ \arg \Bigl( {x+i0 - \alpha_k^* \over x + i0 - \varphi_k^*} \Bigr)
   = \cases { 0, \quad x < \alpha_k^* \ {\rm or} \ x > \varphi_k^*, \cr
              - \pi, \quad \alpha_k^* < x < \varphi_k^*. \cr }
$$
Hence, for all $x>1$, except $x=\alpha_k^*, \varphi_k^*$,
$$ \arg \Bigl( R_0 (x+i0) R_0((x+i0)^{-1}) \Bigr) = - \pi \chi_0(x).
    \eqno(2.37)
$$

   It is also easy to see (from the definition of $R_1(z) = P(z,-\gamma_1)$)
that on the set $\Omega_1 \cup V(\Omega_1)$
$$ \arg \Bigl( R_1 (x+i0) R_1((x+i0)^{-1}) \Bigr)
   = - \gamma_1 (x) + \gamma_1(x^{-1})
   = - (\chi_1 (x) + \chi_1(x^{-1}) ) \eta (x+x^{-1}), \eqno(2.38)
$$
and on the set $\Omega_2^a$
$$ \arg \Bigl( R_2 (x+i0) R_2((x+i0)^{-1}) \Bigr)
   = - \gamma_2 (x) + \gamma_2(x^{-1})
   = - \chi_2^a (x) ( \eta (x+x^{-1}) - \chi_0(t)\pi). \eqno(2.39)
$$
Thus, it follows from (2.36)--(2.39) that on the set
$\Bigl( \Omega_1 \cup V(\Omega_1) \cup \Omega_2^a \Bigr) \cap (1,+\infty)$
$$ \arg \Bigl( R_{012} (x+i0) R_{012}((x+i0)^{-1}) \Bigr) =
   \arg \Bigl( - { {(x+i0)-(x+i0)^{-1}} \over {N(x+i0)} } \Bigr).
$$

   Further, on the set
$(1,\infty)\backslash \Bigl( \Omega_1 \cup V(\Omega_1)
 \cup \Omega_2^a \cup \Omega_2^s \cup \Phi \Bigr)$, according to
$(2.22')$, $(2.22'')$,
$$ \arg \Bigl( { {(x+i0)-(x+i0)^{-1}} \over {N(x+i0)} } \Bigr)
   = - \pi \chi_0(x)
$$
(this argument may not be of the form $-\pi \chi_0(x) + 2 k \pi, \ k\ne 0$,
because in the upper half-plane except the upper half-disk
the argument of the function ${ {z-z^{-1}} \over {N(z)} }$
is contained in $(-\pi,\pi)$).

   So, we obtained that the imaginary parts of the functions
$$ f_1(z) = \ln \biggl( { {z-z^{-1}} \over {N(z)} } \biggr),
$$
$$ f_2(z) = \ln \biggl( R_{0123}(z) R_{0123}(z^{-1}) \biggr),
$$
which are holomorphic in two domains, are equal almost everywhere
on the boundaries of these domains, which completes the proof of the theorem.

\hfill\rule{0.5em}{0.5em}

\bigskip

   {\bf 4.} \ Let us establish some additional properties
of the functions $\rho_R(\tau)$ and $\rho_L(\tau)$.

\medskip

{\bf Lemma 2.3.} {\it The set $\Omega^a_2$
belongs to the absolutely continuous component of the measures
$d\rho^R(\tau)$ and $d\rho^R(\tau)$}
(i.e.~to the intersection of their absolutely continuous components).

\medskip

P r o o f. \  We only prove the lemma for $d\rho_R(\tau)$
(first on the set $\Omega^a_2 \backslash [-2,2]$).
For this we will use the following fact: if under the interval $\delta_l$
$$  \omega =
    {\rm ess}\sup_{\alpha\in\delta_l}\eta^R(\alpha)-
    {\rm ess}\inf_{\alpha\in\delta_l}\eta^R(\alpha)<\pi\,,
$$
then for every $p<\pi \omega^{-1}$ the function
$m^R(\tau \pm ih) = P(\tau \pm ih, \eta^R)$ for $h\downarrow 0$ converge
in $L^p$-norm to a limit $m^R(\tau \pm i0)$ on every compact subset of the
interval $\delta_l$. (This statement is included in lemma 3.4 in more general
form.) Applying to our case this statement and taking into account
(1.12), we have that, under a certain interval
$\tilde\delta_l \subset \delta_l$
(such that in $\delta_l \backslash \tilde\delta_l$ there is no points of
$\Omega^a_2$), the function
$$ \rho_R (\lambda) = {b_{-1}\over\pi} \lim_{\varepsilon \downarrow 0}
      \int\limits_{-\infty}^\lambda {\rm Im} \, m^R(\tau+ i \varepsilon)
      d \tau,
$$
is absolutely continuous, it has almost everywhere the derivative
$$ \rho'_R (\lambda) = {1\over\pi} {\rm Im} \, m^R(\tau+ i 0),
$$
which belongs to $L^1(\tilde\delta_l)$, and
$d\rho_R (\tau) = \rho_R'(\tau)d\tau$.
Since all the set $\Omega^a_2 \backslash [-2,2]$ can be covered with such
$\tilde\delta_l$, this means that $\rho_R(\tau)$ is absolutely continuous
on $\Omega^a_2 \backslash [-2,2]$.

\smallskip

   Let us now show that $\rho_R(\tau)$ is absolutely continuous in some
neighborhood of the set $[-2,2]$. In fact, conditions B) and E) imply
that, for some $\varepsilon>0$ the function $\eta^R(\tau)$
satisfy the condition
$$  \omega =
    {\rm ess}\sup \eta^R(\alpha)-
    {\rm ess}\inf \eta^R(\alpha)<\pi\,,
$$
individually under the intervals $(-2-\varepsilon,0]$ and
$[0,2+\varepsilon)$. This have as a conclusion that
$$ \rho'_R (\lambda) = {1\over\pi} {\rm Im} \, m^R(\tau+ i 0)
$$
belongs to $L^1$ individually on $(-2-\varepsilon,0]$ and
$[0,2+\varepsilon)$. Hence, it belongs to $L^1$ under all
$(-2-\varepsilon,2+\varepsilon)$.
\hfill\rule{0.5em}{0.5em}

\medskip

{\bf Lemma 2.4.} \ {\it The expression
$\sqrt{ \rho'_+ (\tau) \over \rho'_-(\tau) }$
satisfies the H\"older condition and is bounded away from zero
and the infinity under the segment $[-2,2]$.}

\medskip

   P r o o f.   \  Here we only prove the lemma for the right part
$[-1,2]$ of the segment $[-2,2]$ (the proof for the left part is analogous).
As it is seen from conditions B), C) and E), the functions
$\eta^+(\tau)$ и $\eta^-(\tau)$ have the same jump at the point $2$, i.e.
$$ \eta^R(\tau)= \tilde\eta^R(\tau) + \delta_2(\tau), \qquad
   \eta^L(\tau)= \tilde\eta^L(\tau) + \delta_2(\tau),
$$
where
$$ \delta_2(\tau) = \cases{ \varepsilon_1, \quad \tau<2, \cr
                            0, \quad \tau \ge 2, \cr}
$$
with a number
$\varepsilon_1: \ 0<\varepsilon_1<\pi$. Also, $\tilde\eta^R(\tau)$,
$\tilde\eta^L(\tau)$ satisfy the H\"older condition on the set
$[-1,2+\varepsilon], \ \varepsilon>0$. So,
$$ m^R(\lambda) = C P(\lambda,\tilde\eta^R(\tau))
                    P(\lambda,\delta_2(\tau))
                = C P(\lambda,\tilde\eta^R(\tau))
                    (\lambda-1)^{\varepsilon_1 \over \pi},  \eqno(2.40')
$$
$$ {-1 \over m^L(\lambda) } = C P(\lambda,\tilde\eta^L(\tau))
                                P(\lambda,\delta_2(\tau))
                            = C P(\lambda,\tilde\eta^L(\tau))
                                (\lambda-1)^{\varepsilon_1 \over \pi}  \eqno(2.40'')
$$
and $P(\tau+i0,\tilde\eta^R)$, $P(\tau+i0,\tilde\eta^L)$ are bounded
functions satisfying the H\"older condition (see, for example,
[17]), which are bounded away from the real
line\footnote{Since $\tilde\eta^R(\tau)$ and $\tilde\eta^L(\tau)$
satisfy the H\"older condition, the real part of the integral
$$ \int\limits_{-\infty}^{\infty} \tilde\eta^R(\tau)
   ( { 1\over{\tau - \lambda} } - {\tau\over{1+\tau^2}} )d\tau ,
$$
is bounded as $\lambda \to \lambda_0$, $\lambda_0 \in (-1,2+\varepsilon)$,
i.e. $P(\tau+i0,\tilde\eta^R)$ is bounded away from zero and the infinity.
Further, $\tilde\eta^R(\tau)$ is isolated from $0$ and $\pi$,
from where $P(\tau+i0,\tilde\eta^R)$ is isolated from the real line.
So is $P(\tau+i0,\tilde\eta^L)$.}.
Besides, it is seen from (2.40$'$), (2.40$''$) and the equalities
$$ \rho'_R (\lambda) = {1\over\pi} {\rm Im} \, m^R(\tau+ i 0),
$$
$$ \rho'_L (\lambda) = {1\over\pi} {\rm Im} \, { -1 \over m^L(\tau+ i 0)}
$$
that
$$ { \rho'_R (\tau) \over \rho'_L(\tau) } =
   { {\rm Im} \, m^R(\tau+ i 0) \over
     {\rm Im} \, { -1 \over m^L(\tau+ i 0)} } = C
   { {\rm Im} \,
     \Bigl( P(\tau+i0,\tilde\eta^R) (\tau-1)^{\varepsilon_1\over\pi} \Bigr)
   \over {\rm Im} \,
     \Bigl( P(\tau+i0,\tilde\eta^L) (\tau-1)^{\varepsilon_1\over\pi} \Bigr) }
  = C { {\rm Im} \,
     \Bigl( P(\tau+i0,\tilde\eta^R) e^{\varepsilon_1\over\pi} \Bigr) \over
     {\rm Im} \,
     \Bigl( P(\tau+i0,\tilde\eta^L) e^{\varepsilon_1\over\pi} \Bigr) }
$$
is a bounded function satisfying the H\"older condition.
\hfill\rule{0.5em}{0.5em}

\medskip

   In the end of this section we also observe one more property of
$R_\mu(z)$. It follows from the lemma 2.6 that the expression
$$ \ln \mu(e^{i\theta}+e^{-i\theta}) =
   \ln \sqrt{ \rho'_R (e^{i\theta}) \over \rho'_L(e^{i\theta}) }
$$
is bounded and satisfies the H\"older condition.
Hence, the following lemma is true:

\medskip

{\bf Lemma 2.5.} {\it For $\varepsilon \downarrow 0$ the function
$$ R_\mu((1\mp \varepsilon)e^{i\theta}), \quad -\pi \leq \theta \leq \pi,
$$
converges uniformly to its limit
$R_\mu^\pm(e^{i\theta}) \equiv
\lim_{\varepsilon\downarrow 0} R_\mu((1 \mp \varepsilon)e^{i\theta})$
and is bounded. }

\bigskip
\bigskip

\centerline{\bf 3. Deduction of the fundamental equation
of the inverse problem}

\bigskip

   Let us consider for every $k$
$$ g(k,z) \equiv {R(z) \over R(\infty)} z^{-(k+1)} h_k \psi(k,z)
          = {R(z) \over R(\infty)} z^{-(k+1)} h_k
            \Bigl\{ n(z) P_k (z+z^{-1}) + Q_k (z+z^{-1}) \Bigr\},  \eqno(3.1)
$$
where $h_k$ are defined in (0.15). According to (1.4), the function
$g(k,z)$ satisfies the asymptotic formula
$$ g(k,z) \to 1,  \quad z \to \infty.
$$
This function is holomorphic everywhere but the unit circle and a certain
set of the real line, which is separated from the circle, zero and the
infinity. So, according to the Cauchy theorem
$$ g(k,z) = 1 + {1\over2\pi i}
   \int\limits_{\Gamma}{ g(k,\zeta) \over \zeta-z } d\zeta,
$$
where $\Gamma$ is a contour, enclosing all the singularities of the
function $g(k,z)$, concentrated on the real line and the unit circle.
In the part of the plane, bounded by this contour, the functions
$P_k(z+z^{-1})$, $Q_k(z+z^{-1})$ are regular (we remind that they are
polynomials in $z+z^{-1}$). The singularities of the function $n(z)$
are concentrated on the set $\Omega=\Omega_1\cup\Omega_2^s\cup\Omega_2^a$
and the unit circle ${\bf T}$. The singularities of the function $R(z)$ are
concentrated on the set $\Omega_1\cup \Phi\cup\Omega_2^a$ and on the unit
circle. Hence, all the singularities of the function $g(k,z)$ in this part
of the plane lie in the union of the four sets $\Omega_1\cup\Phi$,
$\Omega_2^s$, $\Omega_2^a$, ${\bf T}$,
which lie at positive distances one from another and from the set
${V(\Omega_1)}\cup{V(\Phi)}$. This makes it possible to deform the contour
$\Gamma$ into a finite system of simple contours $\Gamma_1$, $\Gamma_2^s$,
$\Gamma_2^a$ which enclose the sets $\Omega_1\cup\Phi$, $\Omega_2^s$,
$\Omega_2^a$, respectively, and the contour $\Gamma_3$,
formed by two concentric circles bounding an annulus of small area
containing the unit circle. Here the contours $\Gamma_1$, $\Gamma_2^s$,
$\Gamma_2^a$, $\Gamma_3$, can be chosen so close to the singularities
in their interiors that they lie at positive distance one from another
and from the set ${V(\Omega_1)} \cup {V(\Phi)}$. Thus,
$$ g(k,z) = 1 + {1\over2\pi i}
   \biggl\{ \int\limits_{\Gamma_1}{ g(k,\zeta) \over \zeta-z } d\zeta
   + \int\limits_{\Gamma_2^s} { g(k,\zeta) \over \zeta-z } d\zeta
   + \int\limits_{\Gamma_2^a}{ g(k,\zeta) \over \zeta-z } d\zeta
   + \int\limits_{\Gamma_3}{ g(k,\zeta) \over \zeta-z } d\zeta \biggr\}.
   \eqno(3.2)
$$

   Let us first calculate the integrals in the right-hand side of (3.2).
We denote by $O_1$, $O_2^s$, $O_2^a$, $O_3$ the domains that lie inside the
contours $\Gamma_1$, $\Gamma_2^s$, $\Gamma_2^a$, $\Gamma_3$, which are
close to the sets $\Omega_1\cup\Phi$, $\Omega_2^s$, $\Omega_2^a$, ${\bf T}$,
so that the distances between any two of the domains
$O_1$, $O_2^s$, $O_2^a$, $O_3$ and between these domains and the set
$V({\Phi\cup\Omega_1})$ are positive.

Later in this section we will omit somewhere, for the sake of brevity,
the index $k$ of the functions $g(k,z)$, $P_k(z+z^{-1})$, $Q_k(z+z^{-1})$
(although all of them, certainly, depend on $k$).

\bigskip

\quad {\bf Calculation of the integrals along the contours belonging to
$\Gamma_1$}

\bigskip

{\bf Lemma 3.1.} {\it The integral
$\int_{\Gamma_1}$ in (3.2) is calculated by formula
$$
   {1\over2\pi i}\int\limits_{\Gamma_1}{g(\zeta)\over\zeta-z}d\zeta
   = -\int\limits_{\Phi\cup\Omega_1}
     { \alpha^{-2(k+1)} {R(\alpha^{-1})}^{-2} (\alpha-\alpha^{-1})
        \over \alpha - z }
    g(\alpha^{-1}) d\rho_1(\alpha),  \eqno (3.3)
$$
where the measure $d\rho_1(t)$ is defined by the nondecreasing function
$\rho_1(t)$ with }
$$ \rho_1 (t_2) - \rho_1 (t_1)
   = {1\over\pi} \lim_{\varepsilon \downarrow 0}
      \int\limits_{t_1}^{t_2} {\rm Im} \,
      \Bigl( - N(t + i \varepsilon)^{-1} \Bigr) d t.    \eqno(3.4)
$$

\medskip

   We remark that the measure $d\rho_1(t)$ is defined by formula (3.4) on
a larger set than $\Phi\cup\Omega_1$, and that for $\pm 1 \in [t_1,t_2]$
formula (3.4) is senseless. But, according to (3.3), we only need this
measure on $\Phi\cup\Omega_1$, so later on we suppose that $d\rho_1(t)$
is only defined on $\Phi\cup\Omega_1$.

   We also observe that the notation $R(\alpha^{-1})$, $g(\alpha^{-1})$
in formula (3.3) makes sense, because $R(z)$ and $g(z)$ are holomorphic
on the set $V(\Phi) \cup V(\Omega_1)$, i.e. $R(\alpha^{-1})$ and
$g(\alpha^{-1})$ are regular for $ \alpha \in \Phi \cup \Omega_1 $.

\medskip

P r o o f. \ According to definition (3.1),
$$ g(\zeta)= {\zeta^{-(k+1)} \over R(\infty)} h_k R(\zeta) \Bigl( Q(\zeta+\zeta^{-1}) +
                                   n(\zeta) P(\zeta+\zeta^{-1}) \Bigr),
$$
$$ g(\zeta^{-1})= {\zeta^{(k+1)} \over R(\infty)} h_k R(\zeta^{-1}) \Bigl( Q(\zeta+\zeta^{-1}) +
                    n(\zeta^{-1}) P(\zeta+\zeta^{-1}) \Bigr).
$$
Let us exclude $Q(\zeta+\zeta^{-1})$ from these relations:
$$ Q(\zeta+\zeta^{-1}) = R(\infty)
   g(\zeta^{-1}) \zeta^{-(k+1)} h_k^{-1} {R(\zeta^{-1})}^{-1}
   - n(\zeta^{-1}) P(\zeta+\zeta^{-1}),
$$
$$ g(\zeta) = g(\zeta^{-1}) \zeta^{-2(k+1)} {R(\zeta^{-1})}^{-1} R(\zeta)
   + {\zeta^{-(k+1)} \over R(\infty)} h_k R(\zeta)
           \Bigl( -n(\zeta^{-1}) P(\zeta+\zeta^{-1}) +
                   n(\zeta) P(\zeta+\zeta^{-1}) \Bigr)
$$
$$  = g(\zeta^{-1}) \zeta^{-2(k+1)} {R(\zeta^{-1})}^{-1} R(\zeta)
   + P(\zeta+\zeta^{-1}) {\zeta^{-(k+1)} \over R(\infty)} h_k R(\zeta) N(\zeta)
$$
(we remind that $N(\zeta)= n(\zeta) - n(\zeta^{-1})$). We find from
this and from the factorization
$R(\zeta)R(\zeta^{-1})N(\zeta)=-(\zeta-\zeta^{-1})$, obtained in theorem 2,
that
$$
g(\zeta)=f_1(\zeta)(-N(\zeta))^{-1}+f_2(\zeta)\,,
$$
where
$$ f_1(\zeta) = - g(\zeta^{-1}) \zeta^{-2(k+1)} {R(\zeta^{-1})}^{-2}
              (\zeta-\zeta^{-1}) \,,
$$
$$ f_2(\zeta) = P(\zeta+\zeta^{-1}) {\zeta^{-(k+1)} \over R(\infty)}
       h_k {R(\zeta^{-1})}^{-1}  (\zeta-\zeta^{-1}).
$$

    The functions $R(\zeta^{-1})$, ${R(\zeta^{-1})}^{-1}$, $n(\zeta^{-1})$,
and so as the functions $g(\zeta^{-1})$,
$f_1(\zeta)$, $f_2(\zeta)$ are holomorphic outside the unit circle and the
compact set
$$
{V(\Phi\cup\Omega_1)}\cup\Omega_2 \subset {\bf R}.
$$
In particular, all these functions are holomorphic in the domain $O_1$ and,
according to the Cauchy theorem, we have
$$ {1\over2\pi i} \int\limits_{\Gamma_1} {g(\zeta)\over\zeta-z}d\zeta
   = {1\over2\pi i}\int\limits_{\Gamma_1}{f_1(\zeta)\over\zeta-z}
              (-N(\zeta))^{-1}d\zeta\,,
$$
if $z\not\in O_1$.

   Let us consider the function
${ -1\over N(\zeta)}$. According to (1.1), (1.2),
at the zero and the infinity it vanishes and is holomorphic.
In the upper half-plane, except the arc of the unit circle,
${\rm Im}\, N(\zeta) > 0$, from where
${\rm Im}\, (-N(\zeta)^{-1}) > 0$. All the singularities of the function
are concentrated on the unit circle and the parts of the real line,
separated from $-1,0,1,\infty$. As it will be shown below,
$$ (-N(\zeta))^{-1} =
        \int\limits_{\Omega_1\cup\Phi}
        { d\rho_1(\alpha) \over \alpha-\zeta }
             + \tilde N(\zeta) \,, \eqno(3.5)
$$
where the measure $d\rho_1(t)$ is defined by (3.4),
and $\tilde N(\zeta)$ is holomorphic in the domain $O_1$.

   Taking into account the regularity of $\tilde N(\zeta)$ in the domain $O_1$,
$$ \int\limits_{\Gamma_1}{f_1(\zeta) \tilde N(\zeta) \over\zeta-z} d\zeta = 0\,,
$$
$$ {1\over2\pi i} \int\limits_{\Gamma_1} {g(\zeta)\over\zeta-z}d\zeta
   = {1\over2\pi i} \int\limits_{\Gamma_1}{f_1(\zeta)\over\zeta-z}
     \biggl[ { \int\limits_{\Phi\cup\Omega_1}
       { d\rho_1(\alpha) \over \alpha-\zeta } } \biggr]d\zeta\,.
$$

  The family of contours $\Gamma_1$ encloses the set $\Phi\cup\Omega_1$
and are at positive distance from it, which allows to change the order of
integration. We do this and use theorem of residue to calculate the inner
integral:
$$ - {1\over2\pi i} \int\limits_{\Gamma_1}
   {f_1(\zeta) d\zeta \over (\zeta-z) (\zeta-\alpha) }
   = {f_1(\alpha) \over \alpha - z }\,,
$$
because the residue
$$ {\rm Res}\vert_{\zeta=\alpha} \, { f_1(\zeta) \over (\zeta-z) (\zeta-\alpha) } =
   \lim_{\zeta \to \alpha} (\zeta-\alpha) { f_1(\zeta) \over (\zeta-z) (\zeta-\alpha) }
   = { f_1(\alpha) \over \alpha-z } \, .
$$
Thus,
$$ {1\over2\pi i} \int\limits_{\Gamma_1} {g(\zeta)\over\zeta-z}d\zeta
   = \int\limits_{\Phi\cup\Omega_1}
       { f_1(\alpha) \over \alpha - z } d\rho_1(\alpha)\,.
$$
(Here we remind that the contour $\Gamma_1$ is oriented clockwise.)
According to the definition of $f_1(\alpha)$, it follows that
$$
   {1\over2\pi i}\int\limits_{\Gamma_1}{g(\zeta)\over\zeta-z}d\zeta
   = -\int\limits_{\Phi\cup\Omega_1}
     { \alpha^{-2(k+1)} {R(\alpha^{-1})}^{-2} (\alpha-\alpha^{-1})
        \over \alpha - z }
    g(\alpha^{-1}) d\rho_1(\alpha),
$$
as was to be proved.

\medskip

{\normalsize

   Here we remind that we also have to prove the decomposition (3.4) for the
function $(-N(\zeta))^{-1}$. (Note, that this is not the Stiltjes
representation for the functions with positive imaginary part in the upper
half-plane, because the function $(-N(\zeta))^{-1}$, which has singularities
on the unit circle, is not holomorphic in the half-plane.)
Let us denote by $\Gamma_4$ the contour, enclosing the singularities
$(-N(z))^{-1}$ on the real line (remind that $\Gamma_3$ is a contour,
consisting of two concentric circles, close to the unit circle).
By the Cauchy theorem,
$$ (-N(\zeta))^{-1} =
    {1\over2\pi i} \int\limits_{\Gamma_4}
    { (-N(\alpha))^{-1} \over \alpha-\zeta } d\alpha
    + {1\over2\pi i} \int\limits_{\Gamma_3}
    { (-N(\alpha))^{-1} \over \alpha-\zeta } d\alpha.
$$
We denote the second one of these integrals by $\hat N(\zeta)$.
If the contour $\Gamma_3$ is close enough to the unit circle, the function
$N(\zeta)$ is holomorphic in the domain $O_1$. Then for little
$\varepsilon$ and great $M>0$
$$ (-N(\zeta))^{-1} =
    {1\over2\pi i} \lim_{\delta\to+0}
     \biggl( { \int\limits_{-M}^{-1-\varepsilon}
               + \int\limits_{-1+\varepsilon}^{1-\varepsilon}
               + \int\limits_{1+\varepsilon}^M } \biggr)
     \biggl( {
        { (-N(\alpha+i\delta))^{-1} \over \alpha+i\delta-\zeta } -
        { (-N(\alpha-i\delta))^{-1} \over \alpha-i\delta-\zeta }
     } \biggr) d\alpha
        + \hat N(\zeta)
$$

   It can be shown, using the methods of [15], the function $N(z)$ is
decomposed into a product of two functions: $N(z)=N_0(z)N_1(z)$, where
$N_0(z)$ is holomorphic outside the unit circle and is positive on the real line
($N_0(t)>0, \ t \in {\bf R}\backslash\{-1,1\}$), and $N_1(z)$ is holomorphic
outside the real line and has positive imaginary part in the upper half-plane.
Therefore, for the function $-N_1(z)^{-1}$ the following representation is
true:
$$ - N_1(z)^{-1} = a z + b + \int\limits_{-\infty}^\infty { d\sigma(t) \over t-z },
$$
where $a>0$, $b \in {\bf R}$ and $ d\sigma(t) $ is a certain bounded measure on
${\bf R}$, concentrated on the union of the intervals
$(-M,-1-\varepsilon)\cup(-1+\varepsilon,1-\varepsilon)\cup(1+\varepsilon,M)$.
Taking into account the holomorphy of the function $N_0(z)(az+b)$
on this union of intervals,
$$ (-N(\zeta))^{-1} = \hat N(\zeta)
   + {1\over2\pi i} \lim_{\delta\to+0}
     \biggl( { \int\limits_{-M}^{-1-\varepsilon}
              + \int\limits_{-1+\varepsilon}^{1-\varepsilon}
              + \int\limits_{1+\varepsilon}^M } \biggr)
$$
$$   \times
     \biggl( {
        { N_0(\alpha+i\delta)^{-1} N_1(\alpha+i\delta)^{-1}  \over \alpha+i\delta-\zeta }
      - { N_0(\alpha-i\delta)^{-1} N_1(\alpha-i\delta)^{-1} \over \alpha-i\delta-\zeta }
     } \biggr) d\alpha
$$
$$ =  \hat N(\zeta)
  + {1\over2\pi i} \lim_{\delta\to+0}
     \biggl( { \int\limits_{-M}^{-1-\varepsilon}
               + \int\limits_{-1+\varepsilon}^{1-\varepsilon}
               + \int\limits_{1+\varepsilon}^M } \biggr)
$$
$$ \times
   \biggl( {
        { N_0(\alpha+i\delta)^{-1}
       \int\limits_{-\infty}^\infty {  d\sigma(t) \over t-(\alpha+i\delta) }
        \over \alpha+i\delta-\zeta } -
        { N_0(\alpha-i\delta)^{-1}
       \int\limits_{-\infty}^\infty {  d\sigma(t) \over t-(\alpha-i\delta) }
       \over \alpha-i\delta-\zeta }
       } \biggr) d\alpha
$$
$$ =  \hat N(\zeta)
   + {1\over2\pi i} \lim_{\delta\to+0} \int\limits_{-\infty}^{\infty} d\sigma(t)
     \biggl( { \int\limits_{-M}^{-1-\varepsilon}
             + \int\limits_{-1+\varepsilon}^{1-\varepsilon}
             + \int\limits_{1+\varepsilon}^M } \biggr)
$$
$$ \times  \biggl( {
        { N_0(\alpha+i\delta)^{-1} \over
          (\alpha+i\delta-\zeta) (t-(\alpha+i\delta))  } -
        { N_0(\alpha-i\delta)^{-1} \over
          (\alpha-i\delta-\zeta) (t-(\alpha-i\delta)) }
       } \biggr) d\alpha.
$$
The function $ {  N_0(\alpha)^{-1} \over (\alpha-\zeta) (t-\alpha) } $
is holomorphic in the variable $\alpha$ at the points of union
$(-M,-1-\varepsilon)\cup(-1+\varepsilon,1-\varepsilon)\cup(1+\varepsilon,M)$,
except the point $\alpha=t$, at which it has a simple pole. Thus,
$$  {1\over2\pi i} \lim_{\delta\to+0}
     \biggl( { \int\limits_{-M}^{-1-\varepsilon}
             + \int\limits_{-1+\varepsilon}^{1-\varepsilon}
             + \int\limits_{1+\varepsilon}^M } \biggr)
     \biggl( {
        { N_0(\alpha+i\delta)^{-1} \over
          (\alpha+i\delta-\zeta) (t-(\alpha+i\delta))  } -
        { N_0(\alpha-i\delta)^{-1} \over
          (\alpha-i\delta-\zeta) (t-(\alpha-i\delta)) }
       } \biggr) d\alpha
$$
$$  =                           - {\rm Res}\vert_{\alpha=t}\,
          { N_0(\alpha)^{-1} \over (\alpha-\zeta) (t-\alpha) }
       =  { N_0(t)^{-1} \over t-\zeta } \,.
$$
Introducing the measure $d\rho_1(t) \equiv N_0(t)^{-1}d\sigma(t)$, we obtain
$$ (-N(\zeta))^{-1} = \hat N(\zeta)
      + {1\over2\pi i} \int\limits_{-\infty}^{\infty}
        { N_0(t)^{-1} \over t-\zeta } d\sigma(t)
      =  \hat N(\zeta)
        + \int\limits_{-\infty}^{\infty}
        { d\rho_1(t) \over t-\zeta }
      = \tilde N(\zeta)
      + {1\over2\pi i} \int\limits_{\Omega_1\cup\Phi}
        { N_0(t)^{-1} \over t-\zeta } d\sigma(t) \,,
$$
where
$$ \tilde N(\zeta) = \hat N(\zeta)
      + {1\over2\pi i} \int\limits_{{\bf R}\backslash \{\Omega_1\cup\Phi\}}
        { N_0(t)^{-1} \over t-\zeta } d\sigma(t) \,,
$$
which was our goal.

   Taking into account the definition of $\rho_1(t)$ and $\sigma(t)$
the fact that $N_0(t)>0$ for $t \in {\bf R}\backslash{\{-1,1\}}$, we have
$$ \rho_1 (t_2) - \rho_1 (t_1)
   = {1\over\pi} \lim_{\varepsilon \to + 0}
      \int\limits_{t_1}^{t_2} N_0(t)^{-1} {\rm Im} \,
      \Bigl( - N_1(t + i \varepsilon)^{-1} \Bigr) d t
   = {1\over\pi} \lim_{\varepsilon \to + 0}
      \int\limits_{t_1}^{t_2} {\rm Im} \,
      \Bigl( - N(t + i \varepsilon)^{-1} \Bigr) d t\,.
$$
It is seen that $\rho_1(t)$ is non-decreasing.
\hfill\rule{0.5em}{0.5em}

}

\bigskip

\quad {\bf Calculation of the integrals along the contours belonging to
$\Gamma_2^s$}

\bigskip

{\bf Lemma 3.2.} {\it The integral $\int_{\Gamma_2^s}$ in (3.2) is calculated
by formula
$$ {1\over2\pi i} \int\limits_{\Gamma_2^s} {g(\zeta)\over\zeta-z} d\zeta =
  \int\limits_{\Omega_2^s \cap (-1,1)}
     { \alpha^{-2(k+1)} R(\alpha)^2
       \over (\alpha-z) (\alpha^{-1}-\alpha) }
       g(\alpha^{-1}) d\rho_2(\alpha)\,, \eqno (3.6)
$$
where $d\rho_2(\alpha)$ is the measure, defined on the point set
$\Omega_2^s$ by }
$$ \rho_2(\alpha_l)=0, \quad
   \alpha_l\in \Omega_2^s\backslash [-1,1], \eqno(3.7')
$$
$$ \rho_2(\alpha_l^{-1})
    ={d\rho^R(\alpha+\alpha_l^{-1})\over d\rho^L(\alpha+\alpha_l^{-1})}
    \Bigl(d\rho^R(\alpha+\alpha_l^{-1})+d\rho^L(\alpha+\alpha_l^{-1}) \Bigr),
  \quad  \alpha_l^{-1}\in \Omega_2^s \cap (-1,1).  \eqno(3.7'')
$$

The integration is performed on the point set $\Omega_2^s\cap(-1,1)$,
in whose points the functions $R(z)$ and $g(z^{-1})$ are regular, i.e.
$R(\alpha)$ and $g(\alpha^{-1})$ make sense.

\medskip

   P r o o f. \
   According to the considerations of section 2, the part $\Omega_2^s$
of the set $\Omega$ contains of a finite number of pairs os points
$(\alpha_l,\alpha_l^{-1})$,
and the positive masses
$$\mu(\alpha_l) \equiv d\rho^L(\alpha_l+\alpha_l^{-1}),  \qquad
  \mu(\alpha_l^{-1}) \equiv d\rho^R(\alpha_l+\alpha_l^{-1})
$$
of the measures $d\rho^L(\lambda)$, $d\rho^R(\lambda)$ correspond to them
(not to mix with $\mu(t)$, introduced earlier for $t\in \Omega_2^a$).
Let us investigate the behavior of the function $g(\zeta)$ in the neighborhood
of $\alpha_l$ and $\alpha_l^{-1}$. It follows from the definition of $n(z)$
and the representation (2.3) that
$$ n(\zeta) = { 1 \over \alpha_l - \zeta} \mu(\alpha_l) +
            { 1 \over \alpha_l^{-1} - \zeta} \mu(\alpha_l^{-1}) + n_l(\zeta),
   \eqno (3.8)
$$
$$ N(\zeta) = \biggl\{ { 1 \over \alpha_l - \zeta} +
                     { 1 \over \alpha_l^{-1} - \zeta} \biggr\}
            \Bigl( \mu(\alpha_l)+\mu(\alpha_l^{-1}) \Bigr) + N_l(\zeta)\,,
   \eqno (3.9)
$$
where the functions $n_l(\zeta)$, $N_l(\zeta)$ are holomorphic in a
neighborhood of the points $\alpha_l$ and $\alpha_l^{-1}$. According to
(2.23) and theorem 2, $\alpha_l^*=\alpha_l$, the function $R(\zeta)$
is holomorphic in the points $\alpha_l$, $\alpha_l^{-1}$ and $R(\alpha_l)=0$,
$R(\alpha_l^{-1})\not=0$. Hence, both functions $R(\zeta)$, $R(\zeta)n(\zeta)$
are holomorphic in the points $\alpha_l$, and so is the function
$g(\zeta)$. The function $R(\zeta)$ is holomorphic in the point
$\alpha_l^{-1}$, and the function $R(\zeta)n(\zeta)$ has a simple pole in
it with residue $-\mu(\alpha_l^{-1}) R(\alpha_l^{-1})$. Let
$\Gamma_l$, $\Gamma_l^*$ be contours enclosing the points
$\alpha_l,\alpha_l^{-1}$ and close enough to them
(it is seen that all the contour $\Gamma_2^s$ consists of such
$\Gamma_l$, $\Gamma_l^*$). Write
$$ {1\over2\pi i}\int\limits_{\Gamma_l}{g(\zeta)\over\zeta-z}d\zeta=0\,,
$$
$$ {1\over2\pi i}\int\limits_{\Gamma_l^*}{g(\zeta)\over\zeta-z}d\zeta
   = - { \mu(\alpha_l^{-1}) R(\alpha_l^{-1})
       P(\alpha_l+\alpha_l^{-1}) \alpha_l^{k+1} h_k
              \over R(\infty) (-\alpha_l^{-1}+z) }\,.
$$
Hence,
$$ {1\over 2\pi i} \int\limits_{\Gamma_2^s} {g(\zeta)\over\zeta-z}d\zeta
   = {1\over R(\infty)} \sum_{\alpha_l\in\Omega_2^s}
       { R(\alpha_l^{-1}) \alpha_l^{k+1} h_k
         \over \alpha_l^{-1} - z } \mu(\alpha_l^{-1}) P(\alpha_l+\alpha_l^{-1})\,.
    \eqno(3.10)
$$

   Since $R(\alpha_l)=0$, we have from (3.1) and (3.8), that
$$ g(\alpha_l) = \alpha_l^{-(k+1)} {h_k \over R(\infty)}
                 P(\alpha_l+\alpha_l^{-1})  \mu(\alpha_l)
                 \lim_{\zeta\to\alpha_l} \Bigl\{(\alpha_l-\zeta)^{-1}R(\zeta)\Bigr\}.
$$
On the one hand, we have from theorem 2 that
$$ \lim_{\zeta\to\alpha_l} R(\zeta)N(\zeta) =
    { \alpha_l-\alpha_l^{-1} \over R(\alpha_l^{-1}) } \,,
$$
and, on the other hand, we have from (3.9) that
$$ \lim_{\zeta\to\alpha_l} R(\zeta)N(\zeta) =
   \Bigl( \mu(\alpha_l)+\mu(\alpha_l^{-1}) \Bigr)
   \lim_{\zeta\to\alpha_l} \Bigl\{ {(\alpha_l^{-1}-\zeta)}^{-1} R(\zeta) \Bigr\}.
$$
Therefore,
$$ \lim_{\zeta\to\alpha_l} {(\alpha_l^{-1} - \zeta)}^{-1} R(\zeta) =
    { \alpha_l-\alpha_l^{-1} \over R(\alpha_l^{-1}) }
    { \Bigl( \mu(\alpha_l)+\mu(\alpha_l^{-1}) \Bigr) }^{-1}\,,
$$
$$ g(\alpha_l)= \alpha_l^{-(k+1)} {h_k \over R(\infty)}
                                  P(\alpha_l+\alpha_l^{-1})
     \biggl({\mu(\alpha_l)\over\mu(\alpha_l)+\mu(\alpha_l^{-1})}\biggr)
     \biggl( {\alpha_l-\alpha_l^{-1}\over R(\alpha_l^{-1})}\biggr)\,,
$$
$$ P(\alpha_l+\alpha_l^{-1}) =
   \biggl( { \mu(\alpha_l)+\mu(\alpha_l^{-1}) \over \mu(\alpha_l) }\biggr)
   \biggl( { R(\alpha_l^{-1}) \over \alpha_l-\alpha_l^{-1} }\biggr)
   { \alpha_l^{k+1} R(\infty) \over h_k } g(\alpha_l).
$$
Substituting this into the right-hand side of (3.10), we obtain
$$ {1\over2\pi i} \int\limits_{\Gamma_2^s} {g(\zeta)\over\zeta-z}d\zeta =
    \sum_{\alpha_l\in\Omega_2^s}
     { R(\alpha_l^{-1})^2
      \alpha_l^{2(k+1)} \over (\alpha_l^{-1}-z) (\alpha_l-\alpha_l^{-1}) }
      \cdot {\mu(\alpha_l^{-1})\over\mu(\alpha_l)}
      \Bigl( \mu(\alpha_l)+\mu(\alpha_l^{-1}) \Bigr) g(\alpha_l).
$$
Let $d\rho_2(\alpha)$ be the measure, defined on the point set
$\Omega_2^s$ by (3.7):
$$ \rho_2(\alpha_l)=0,
$$
$$ \rho_2(\alpha_l^{-1})={\mu(\alpha_l^{-1})\over\mu(\alpha_l)}
                         \Bigl(\mu(\alpha_l)+\mu(\alpha_l^{-1}) \Bigr)
    ={d\rho^R(\alpha+\alpha_l^{-1})\over d\rho^L(\alpha+\alpha_l^{-1})}
    \Bigl(d\rho^R(\alpha+\alpha_l^{-1})+d\rho^L(\alpha+\alpha_l^{-1}) \Bigr).
$$
Then the latter integral can be rewritten as (3.7), as was to be proved.
\hfill\rule{0.5em}{0.5em}

\bigskip

\quad{\bf Calculation of the integrals along the contours belonging to
$\Gamma_2^a$}

\bigskip

{\bf Lemma 3.3.} {\it The following equality is true for the integral
along the contour $\Gamma_2^a$:
$$ {1\over2\pi i}\int\limits_{\Gamma_2^a} {g(\zeta)\over\zeta-z} d\zeta
   = {1\over2\pi i}\int\limits_{\Omega_2^a}
     { g^+(\alpha)-g^-(\alpha) \over \alpha-z } d\alpha.
   \eqno(3.11)
$$
The functions $\chi_2(\alpha)g^+(\alpha,x)$ and $\chi_2(\alpha)g^-(\alpha,x)$
belong to the Hilbert space $L^2({\bf R})$. }

P r o o f. \ The set $\Omega_2^a$ lies in a finite union of intervals
$\Delta_k^{(1)}$, $\Delta_k^{(2)}$, $ {V(\Delta_k^{(1)})}$,
${V(\Delta_k^{(2)})}$ of the real line. The contour $\Gamma_2^a$
consists of a finite number of contours each of which encloses the subsets
$\Delta_k^{(1)}\cap\Omega_2^a$, $\Delta_k^{(2)}\cap\Omega_2^a$,
$V(\Delta_k^{(1)})\cap\Omega_2^a$, $V(\Delta_k^{(2)})\cap\Omega_2^a$
of the set $\Omega_2^a$ and lies close enough to it.
The calculations of the integrals along each of these contours are completely
similar and we will examine only one of them, for example,
the contour that encloses $V(\Delta_k^{(2)})\cap\Omega_2^a$. Let
$$ V(\Delta_k^{(2)}) = (\beta_k^{-1},\varphi_k^{-1}),
   \quad \overline\beta_k = \inf(V(\Delta_k^{(2)})\cap\Omega_2^a),
   \quad \overline\varphi_k = \sup(V(\Delta_k^{(2)})\cap\Omega_2^a).
$$
Since the compact set $V(\Delta_k^{(2)})\cap\Omega_2^a$ is contained
in the interval $V(\Delta_k^{(2)})$, we have
$\beta_k^{-1}<\overline\beta_k<\overline\varphi_k<\varphi_k^{-1}$,
and, without changing the value of the integral, the corresponding contour
can be replaced by a contour $\gamma(h)$, consisting of two intercepts
$$ \gamma^{\pm}(h)=\{ \zeta \mid
              \zeta=\alpha \pm i h, \, \beta'_k<\alpha<\varphi'_k\}
$$
and of two intercepts connecting the endpoints of $\gamma^{\pm}(h)$.
Here $\alpha'_k$, $\varphi'_k$ denote numbers, chosen so that they satisfy
the inequalities
$$ \beta_k^{-1} < \beta'_k < \overline\beta_k,
   \quad \overline\varphi_k < \varphi'_k < \varphi_k^{-1},
$$
and $h>0$ is a small number. On the linear segments of the contour $\gamma(h)$
the function $g(\zeta)(\zeta-z)^{-1}$ is continuous. So, the integrals
along these intervals tend to zero when $h\downarrow 0$, and
$$ {1\over2\pi i} \int\limits_{\gamma(h)} {g(\zeta)\over\zeta-z}d\zeta=
   {1\over2\pi i} \lim_{h\downarrow 0}
   \biggl(\int\limits_{\gamma^+(h)\cup \gamma^-(h)}
   {g(\zeta)\over\zeta-z}d\zeta\biggr)
$$
(with the corresponding orientation of the integration along $\gamma^\pm(h)$).
Let us assume that under the interval $[\alpha'_k,\varphi'_k]$, the function
$g(\alpha\pm ih)$ converges in $L^p$-norm ($p\geq1$) to a finite limite
$$ \lim_{h\downarrow 0} g(\alpha \pm ih) = g^{\pm}(\alpha) .
$$
Then
$$ \lim_{h\downarrow 0}\int\limits_{\gamma^{\pm}(h)} {g(\zeta)\over\zeta-z}d\zeta
   = \int\limits_{\beta'_k}^{\varphi'_k}
     { \pm g^{\pm}(\alpha) \over \alpha-z } d\alpha
$$
and
$$ {1\over2\pi i} \int\limits_{\gamma(h)} { g(\zeta) \over\zeta-z } d\zeta
   = {1\over2\pi i} \int\limits_{\beta'_k}^{\varphi'_k}
     { g^+(\alpha)-g^-(\alpha) \over \alpha-z } d\alpha.
$$
Since the function $g(\zeta)$ is holomorphic at the points
$V(\Delta_k^{(2)}) \backslash \Omega_2^a$,
$$ g^+(\alpha)-g^-(\alpha)=0,
   \quad \alpha \in V(\Delta^{(2)}_k) \backslash \Omega_2^a.
$$
Therefore, we have
$$ {1\over2\pi i} \int\limits_{\gamma(h)} {g(\zeta)\over\zeta-z}d\zeta
   = {1\over2\pi i} \int\limits_{V(\Delta_k^{(2)})\cap\Omega_2^a}
     { g^+(\alpha)-g^-(\alpha) \over \alpha-z} d\alpha, \eqno(3.12)
$$
or
$$ \lim_{h\to0} \Vert g(\alpha \pm ih)
   - g^{\pm}(\alpha)\Vert_{L^p[\beta'_k,\varphi'_k]}=0, \quad p\ge 1. \eqno(3.13)
$$

   So, let us prove (3.13). It follows from the definition of the function
$R(\zeta)$ and the specifics of the function $n(\zeta)$ that the following
factorization is true:
$$ R(\zeta)=R^{(1)}_k(\zeta)R^{(2)}_k(\zeta),
   \quad n(\zeta)=n^{(1)}_k(\zeta)n^{(2)}_k(\zeta),
$$
where
$R^{(1)}_k(\zeta)$ and $n^{(1)}_k(\zeta)$ are functions whose singularities
are concentrated on the set $V(\Delta_k^{(2)})\cap\Omega_2^a$,
and the functions $R^{(2)}_k(\zeta)$, $n^{(2)}_k(\zeta)$ are holomorphic
on the set $V(\Delta_k^{(2)})$. The functions $R^{(1)}_k(\zeta)$,
$n^{(1)}_k(\zeta)$ are defined by formulas
$$ R^{(1)}_k(\zeta) = P\Bigl(\zeta, - J(\alpha) \gamma_2(\alpha)\Bigr), \eqno(3.14)
$$
$$ n^{(1)}_k(\zeta) = P\Bigl(\zeta,
         \pi - J(\alpha) \eta^R(\alpha+\alpha^{-1})\Bigr),
   \eqno(3.15)
$$
where $J(\alpha)$ is the indicator of the interval $V(\Delta_k^{(2)})$.

   (The representability of the function $R(\zeta)$ in the form of such product
is an immediate conclusion of its definition; the representability of
$n(\zeta)$ can be proved by decomposing $n(z)$ in a product of two functions
of the form $P(z,\gamma)$ and $\hat P(z,\hat\gamma)$ with $\gamma$ and
$\hat\gamma$, defined by the limit values of the argument of $n(z)$:
see, for example, [15].)

\medskip

   Hence,
$$ g(\zeta)=\tilde g^{(1)}_k(\zeta)R^{(1)}_k(\zeta)
   + \tilde g^{(2)}_k(\zeta)R^{(1)}_k(\zeta) n^{(1)}_k(\zeta)\,,
$$
where the functions $\tilde g_k^{(1)}(\zeta)$ and $\tilde g_k^{(2)}(\zeta)$
are holomorphic in a neighborhood of $V(\Delta^{(2)}_k)$. Therefore, the
functions $\tilde g^{(j)}_k(\alpha\pm i h)$  $(j=1, 2)$ converge uniformly
to their limits on the interval
$[\beta'_k,\varphi'_k] \subset V(\Delta^{(2)}_k)$ as $h\downarrow 0$.
So, to prove equalities (3.13) and (3.12), it is sufficient to check
that the functions $R^{(1)}_k (\alpha \pm ih )$,
$R^{(1)}_k (\alpha \pm ih ) n^{(1)}_k( \alpha \pm h)$ converge
to their limits in $L^p$-norm when $h\downarrow 0$.
This is a consequence of conditions B), C) and the following lemma.

\bigskip

{\bf Lemma 3.4.} {\it Let
$\delta(\alpha)=\delta_1(\alpha)+\delta_2(\alpha)$ and suppose that the
function $\delta_2(\alpha)$ satisfies the H\"older condition on an interval
$(\alpha_1,\alpha_2)\subset {\bf R}$ and that on the same interval
the function $\delta_1(\alpha)$ satisfies the inequality
$$
\omega={\rm ess}\sup_{\alpha_1<\alpha<\alpha_2}\delta_1(\alpha)-{\rm
ess}\inf_{ \alpha_1<\alpha<\alpha_2}\delta_1(\alpha)<\pi.
$$
Then for any $p<\pi\omega^{-1}$ the function $P(\alpha \pm ih, \delta)$
converges in $L^p$-norm to a limit $P^{\pm}(\alpha,\delta)$ on all compact
subsets of $(\alpha_1,\alpha_2)$.}

\medskip

P r o o f.\quad In our conditions the role of the interval
$(\alpha_1,\alpha_2)$ is played by
$ V(\Delta_k^{(2)}) = (\beta_k^{-1},\varphi_k^{-1})$, so in this lemma we
continue denoting the indicator of the interval $(\alpha_1,\alpha_2)$
by $J(\alpha)$. Let
$$ d_1(\alpha) = J(\alpha)\left(\delta_1(\alpha)-C\right)\,,
$$
$$ d_2(\alpha) = \delta(\alpha)-d_1(\alpha)=\left(1-J(\alpha)\right)\delta_1
   (\alpha)+\delta_2(\alpha)+J(\alpha)C\,,
$$
where
$$ C={1\over2}\left({\rm ess}\max_{\alpha_1<\alpha<\alpha_2}\delta_1(\alpha)+
   {\rm ess}\min_{\alpha_1<\alpha<\alpha_2}\delta_1(\alpha)\right).
$$
Then $\delta(\alpha)=d_1(\alpha)+d_2(\alpha)$,
$$ P(z,\delta)=P(z,d_1)P(z,d_2), \eqno(3.16)
$$
where the function $d_1(\alpha)$ satisfies the inequality
$$ {\rm ess}\max_{\alpha_1<\alpha<\alpha_2}\vert d_1(\alpha)\vert
   = {\omega\over2}<{\pi\over2}, \eqno(3.17)
$$
and the function $d_2(\alpha)$ satisfies the H\"older condition
on the interval $(\alpha_1,\alpha_2)$. Thus, $P(z,d_2)$ converges uniformly
to its limit on every segment
$[\alpha'_1,\alpha'_2]\subset(\alpha_1, \alpha_2)$. Let us introduce

\smallskip

{\bf Definition.} \ {\it Let $A$ be a set of real line, $U_\delta(A)$
be its $\delta$-neighborhood, and $f(z)$ be a holomorphic on
$U_\delta(A)\backslash A$ function. We say that the function $f(z)$
belongs locally to the Hardy class $H^p$ in a neighborhood of the set $A$,
if for a certain $\delta>0$ the functions
$f(t\pm i\varepsilon), \ t\in {\bf R} \cap U_\delta(A),$
converge as $\varepsilon\downarrow 0$ in $L^p$-norm. }

\smallskip

   It follows from the inequality (3.17) that the function
$P(z,d_1)$ belongs locally to the Hardy space $H^p$ for any
$p<\pi\omega^{-1}$ in a neighborhood of
$[\alpha'_1,\alpha'_2]\subset(\alpha_1, \alpha_2)$
(see, for example, [18], [19]).

   We explain how to demonstrate it. Let us decompose
$P(z,d_1)$ in the product of two functions one of which
($P(z,J d_1)$) is holomorphic outside the segment $[\alpha_1,\alpha_2]$,
and the other one is holomorphic on $(\alpha_1,\alpha_2)$.
Then in the domain of holomorphy of $P(z, J d_1)$
$$ \sup |\arg P(z,J d_1)| < {\omega\over2}. \eqno(3.18)
$$
The problem is reduced to the existence of the limits
of the function $P(z,J d_1)$ on the interval
$[\alpha'_1,\alpha'_2]\subset(\alpha_1, \alpha_2)$ in $L^p$-norm.
With the help of a conformal mapping, we reduce the problem to the
following: A holomorphic in the disk function, which satisfies in
the disk condition (3.18), belongs to the space $H^p$ in the unit disk.
The following statement provides us with this fact (see [19], p.100,
ex.~№13):

{\sl Smirnov's theorem.} {\it If $f(z)$ is analytic inside the disk
$|z|<1$ and ${\rm Re}\, f(z) \ge 0$, then $f \in H^p$ for all $p<1$.}

After raising the function, which satisfy inside the disk property (3.18),
to the power $\pi\omega^{-1}$, we apply the last statement.

   Thus,
$$ \lim_{h\downarrow 0} \max_{\alpha'_1<\alpha<\alpha'_2}
   \vert P(\alpha \pm ih , d_2) - P^\pm(\alpha,d_2)\vert=0\,,
$$
$$ \lim_{h\downarrow 0}\int\limits_{-\pi}^\pi
    \vert P(\alpha \pm ih, d_1) - P^\pm(\alpha,d_1)\vert^p d\alpha=0,
$$
and, together with (3.16), this means (3.13), which proves lemma 3.4.
 \hfill\rule{0.5em}{0.5em}

\medskip

   It follows from conditions B) and C) that the function
$\eta^R(\tau) - \eta^L(\tau)$ satisfies the H\"older condition.
Also, the inequality is true
$$ -\pi<\eta^R(\alpha)-\eta^L(\alpha)< \pi
$$
on the interval $\tilde\Delta_k^{(2)}$
(which is the image of $V(\Delta_k^{(2)})$ under the transformation
$\tau = \alpha + \alpha^{-1}$). So, the function
$$1 + P(\lambda+i0,\eta^L(\tau)-\eta^R(\tau))
$$
satisfies the H\"older condition on this interval (see [17]), too,
and is equal to zero almost everywhere. Hence, the function
$$ \beta(\tau)=\arg\left(1+P(\tau+i0,\eta^L -\eta^R)\right)-\pi
$$
satisfies the H\"older condition on the interval $\tilde\Delta_k^{(2)}$.
Since
$$ b_{-1} m^R(\lambda)= C_1 P(\lambda,\eta^R(\tau)); \qquad
   -{ 1 \over {b_{-1} m^L(\lambda)} } = C_2 P(\lambda,\eta^L(\tau)),
$$
then
$$ M(\lambda)= b_{-1} m^R(\lambda)
   -{ 1 \over {b_{-1} m^L(\lambda)} } =
   C_1 P(\lambda,\eta^R) + C_2 P(\lambda,\eta^L) =
$$
$$     = C_1 P(\lambda,\eta^R) (1 + C_3 P(\lambda,-\eta^R) P(\lambda,\eta^L))
       = m^R(\lambda) (1 + C_3 P(\lambda,\eta^L -\eta^R)),
$$
where the constants $C_i>0$, from where
$$ \arg M(\tau+i0)=\eta(\tau) =
   \eta^R(\tau) + \arg(1 + C_3 P(\tau+i0,\eta^L - \eta^R)),
$$
that is
$$ \eta(\alpha+\alpha^{-1})
   =\eta^R(\alpha+\alpha^{-1})+\beta(\alpha+\alpha^{-1})+\pi \,.
$$
It follows from the definition of $\gamma_2(\alpha)$ (formula (2.30))
and from what $\eta(\alpha+\alpha^{-1})=0$ when
$\alpha \in V(\Delta_k^{(2)})\backslash \Omega_2^a$ (see lemma 2.2)
that we have the following equalities on the interval $V(\Delta_k^{(2)})$:
$$ -\gamma_2(t) = \chi_2^a(t) \Bigl(
           {1\over2}\eta^R(t+t^{-1}) + {1\over2} \beta(t+t^{-1}) \Bigr)
    = {1\over2}\eta^R(t+t^{-1}) + {1\over2} \beta(t+t^{-1}) \,,
$$
$$ - \eta^R(t+t^{-1}) - \gamma_2(t)
   = - {1\over2} \eta^R(t+t^{-1})+ {1\over2}\beta(t+t^{-1})\,.
$$
Taking into account (3.14), (3.15), we conclude that
$$ R_k^{(1)}(\zeta)=P(\zeta,\delta_1+\delta_2); \qquad
   R_k^{(1)}(\zeta)n_k^{(1)}(\zeta)=-P(\zeta,-\delta_1+\delta_2), \eqno(3.19)
$$
where
$$ \delta_1(\alpha)={J(\alpha)\over2}\eta^R(\alpha+\alpha^{-1}),\qquad
   \delta_2(\alpha)={J(\alpha)\over2}\beta(\alpha+\alpha^{-1})\,,
$$
and $J(\alpha)$ is the indicator of $V(\Delta_k^{(2)})$. The function
$\delta_2(\alpha)$, evidently, satisfies the H\"older condition
on the interval $V(\Delta_k^{(2)})$ and, according to condition B), the
functions $\pm\delta_1(\alpha)$ satisfy the inequality
$$ {\rm ess} \max_{\alpha\in V(\Delta_k^{(2)})} (\pm\delta_1(\alpha)) -
   {\rm ess} \min_{\alpha\in V(\Delta_k^{(2)})} (\pm\delta_1(\alpha))=
   { \omega \over 2 } < {\pi\over2} \,.
$$

   In view of lemma 3.4 this means that functions (3.19)
converge in $L^p$-norm to their limits on the segment
$[\beta'_k,\varphi'_k]\subset V(\Delta_k^{(2)})$ for all
$p<2\pi\omega^{-1}$. Since $\omega<\pi$, they also converge in
$L^2(\beta'_k,\varphi'_k)$.

Thus, we proved equality (3.13) and formula (3.12).
Analogous speculations can be applied to all the other components
of the contour $\Gamma_2^a$, which completes the proof of lemma 3.3.
 \hfill\rule{0.5em}{0.5em}

\medskip

{\bf Remark.} {\it It follows from the proof above that the functions
$\chi_2(\alpha)g^+(\alpha,x)$ and $\chi_2(\alpha)g^-(\alpha,x)$
belong to the Hilbert space $L^2({\bf R})$. }

\bigskip

{\bf Lemma 3.5.} {\it The following equalities hold almost everywhere
on the set $\Omega_2^a$:
$$ {\alpha-\alpha^{-1} \over \vert \alpha -\alpha^{-1} \vert}
   \left(g^+(\alpha^{-1})+g^-(\alpha^{-1})\right) =
$$
$$ ={ 2i \alpha^{2(k+1)} q(\alpha) \left(g^+(\alpha)-g^-(\alpha) \right)
      + i m (\alpha) {\alpha \over \vert \alpha \vert}
        \left(g^+(\alpha^{-1})-g^-(\alpha^{-1})\right)
      \over p(\alpha) },   \eqno(3.20)
$$
with
$$ p(\alpha) = \left( \mu(\alpha+\alpha^{-1}) +
                      \left( \mu(\alpha+\alpha^{-1}) \right) ^{-1} \right)
               \tan\vert\gamma_2(\alpha)\vert >0, \eqno(3.21)
$$
$$ q(\alpha) = {\vert R(\alpha^{-1}) \vert \over \vert R(\alpha) \vert }
     \left( \mu(\alpha + \alpha^{-1})  \right)^{s(\alpha)} > 0,
    \eqno(3.22)
$$
$$ m(\alpha) = s(\alpha)  \left( \mu(\alpha+\alpha^{-1}) -
           \left( \mu(\alpha+\alpha^{-1}) \right) ^{-1} \right)
               \cos\chi_0(\alpha )\pi.   \eqno(3.23)
$$
(The indicator $\chi_0(\alpha)$ is defined above in theorem 2,
and the functions $\mu(\alpha)$, $\gamma_2(\alpha)$ and
$s(\alpha)$ are defined by (2.5), (2.30) and (2.31), resp.)}

\medskip

P r o o f. \ It can be easily seen from the definition (2.5) of the function
$\mu(t)$ that (3.21), (3.22), (3.23)
can be rewritten in the following way:

$$ p(\alpha) = \left(\sqrt{ \sigma(\alpha^{-1})\over\sigma(\alpha) } +
                     \sqrt{ \sigma(\alpha)\over\sigma(\alpha^{-1}) } \right)
               \tan\vert\varphi_2(\alpha)\vert >0, \eqno(3.21')
$$
$$ q(\alpha) = {\vert R(\alpha^{-1}) \vert \over \vert R(\alpha) \vert }
    \sqrt { \sigma(\alpha^{-1}) \over\sigma(\alpha) } > 0,   \eqno(3.22')
$$
$$ m(\alpha) = \left(\sqrt{\sigma(\alpha^{-1}) \over \sigma(\alpha)} -
                     \sqrt{\sigma(\alpha)\over\sigma(\alpha^{-1})}\right)
               \cos\chi_0(\alpha )\pi,   \eqno(3.23')
$$
where
$$ \sigma(\alpha)\equiv \cases { \rho'_R(\alpha+\alpha^{-1}), \quad |\alpha|<1, \cr
                          \rho'_L(\alpha+\alpha^{-1}), \quad |\alpha|>1. \cr }
$$

Thus, we will prove the lemma with this coefficients.

\medskip

Since the functions $P(z+z^{-1})$, $Q(z+z^{-1})$
are holomorphic everywhere but $0$ and $\infty$, the limit values of the
function $g(z)$ on the set $\Omega_2^a$ equal
$$
g^\pm(\alpha)=\alpha^{-(k+1)} {h_k \over R(\infty)}
  \left( n^\pm(\alpha) R^\pm(\alpha) P(\alpha+\alpha^{-1}) +
         R^\pm(\alpha) Q(\alpha+\alpha^{-1}) \right).
$$
Therefore
$$ g^+(\alpha)-g^-(\alpha)
  = \alpha^{-(k+1)} {h_k \over R(\infty)}
    \biggl\{ \left( n^+(\alpha) R^+(\alpha) -
                n^-(\alpha) R^-(\alpha) \right) P(\alpha+\alpha^{-1})=
$$
$$ + \left( R^+(\alpha) - R^-(\alpha) \right) Q(\alpha+\alpha^{-1})\biggr\} .
    \eqno(3.24)
$$

Let us now consider the limit values of the functions $g(z)$ at the point
$\alpha^{-1}$:
$$ \alpha^{-(k+1)} h_k^{-1} R(\infty) g^\pm(\alpha^{-1})
    = \Bigl\{ n^\pm(\alpha^{-1}) R^\pm(\alpha^{-1}) P(\alpha+\alpha^{-1}) +
              R^\pm(\alpha^{-1}) Q(\alpha^{-1}+\alpha^{-1}) \Bigr\} .
$$
The two latter equalities (for $g^+(\alpha^{-1})$ and $g^-(\alpha^{-1})$)
provides us with a system of equations for the components of the functions
$P(\alpha+\alpha^{-1})$, $Q(\alpha+\alpha^{-1})$ with the natrix
$$ \pmatrix{ R^+(\alpha^{-1}) n^+(\alpha^{-1}) & R^+(\alpha^{-1}) \cr
             R^-(\alpha^{-1}) n^-(\alpha^{-1}) & R^-(\alpha^{-1}) \cr },
$$
whose deteminant is
$$ D(\alpha)=R^+(\alpha^{-1}) R^-(\alpha^{-1})
    \left(n^+(\alpha^{-1})-n^-(\alpha^{-1})\right)\,.
$$
Taking into account the definition of $n(z)$, formulas
(2.3), (2.4) and the absolute continuity of the functions
$\rho_R(\tau+\tau^{-1})$, $\rho_L(\tau+\tau^{-1})$ on $\Omega_2^a$,
we have
$$ D(\alpha)= 2\pi i { \vert R(\alpha^{-1})\vert }^2 \sigma(\alpha^{-1})\not=0
$$
almost everywhere on the set $\Omega_2^a$ (with $\sigma(\alpha)$, defined
as above). Thus, the linear system is uniqually solvable and
$$ P(\alpha+\alpha^{-1}) = { \alpha^{-(k+1)} \over h_k } R(\infty) D(\alpha)^{-1}
    \Bigl\{ R^-(\alpha^{-1}) g^+(\alpha^{-1}) -
       R^+(\alpha^{-1}) g^-(\alpha^{-1}) \Bigr\} \,,
$$
$$ Q(\alpha+\alpha^{-1}) = { \alpha^{-(k+1)} \over h_k } R(\infty) D(\alpha)^{-1}
   \Bigl\{ - R^-(\alpha^{-1}) n^-(\alpha^{-1}) g^+(\alpha^{-1})
      +R^+(\alpha^{-1}) n^+(\alpha^{-1}) g^-(\alpha^{-1}) \Bigr\} \,.
$$
Substituting these expressions into (3.24), we obtain
$$ g^+(\alpha)-g^-(\alpha) =
   \alpha^{-2(k+1)} D(\alpha)^{-1}
   \biggl\{ \Bigl( { n^+(\alpha) R^+(\alpha)
                   - n^-(\alpha) R^-(\alpha) } \Bigr)
$$
$$ \times \Bigl( { R^-(\alpha^{-1}) g^+(\alpha^{-1}) -
                   R^+(\alpha^{-1}) g^-(\alpha^{-1}) } \Bigr)
$$
$$
    + \Bigl( { R^+(\alpha) - R^-(\alpha)  } \Bigr)
       \Bigl( { - R^-(\alpha^{-1}) n^-(\alpha^{-1}) g^+(\alpha^{-1})
                + R^+(\alpha^{-1}) n^+(\alpha^{-1}) g^-(\alpha^{-1}) }
       \Bigr) \biggr\}\,
$$

$$
    = \alpha^{-2(k+1)}
       D(\alpha)^{-1}
     \biggl\{ \Bigl[  n^+(\alpha) R^+(\alpha) R^-(\alpha^{-1})
               - n^-(\alpha) R^-(\alpha) R^-(\alpha^{-1})
$$
$$             - R^+(\alpha) R^-(\alpha^{-1}) n^-(\alpha^{-1})
               + R^-(\alpha) R^-(\alpha^{-1}) n^-(\alpha^{-1})  \Bigr]
         g^+(\alpha^{-1}) +
$$
$$    +  \Bigl[   - n^+(\alpha) R^+(\alpha) R^+(\alpha^{-1})
                  + n^-(\alpha) R^-(\alpha) R^+(\alpha^{-1})
$$
$$                + R^+(\alpha) R^+(\alpha^{-1}) n^+(\alpha^{-1})
                  - R^-(\alpha) R^+(\alpha^{-1}) n^+(\alpha^{-1})  \Bigr]
         g^-(\alpha^{-1}) \biggr\}
$$

$$
    = \alpha^{-2(k+1)} D(\alpha)^{-1} \
       \biggl\{ R^-(\alpha) R^-(\alpha^{-1})
       \Bigl[ { { R^+(\alpha) \over R^-(\alpha) } n^+(\alpha) - n^-(\alpha)
               - { R^+(\alpha) \over R^-(\alpha) } n^-(\alpha^{-1})
                   + n^-(\alpha^{-1}) } \Bigr] g^+(\alpha^{-1})
$$
$$  +  R^+(\alpha) R^+(\alpha^{-1})
       \Bigl[ { - n^+(\alpha) + { R^-(\alpha) \over R^+(\alpha) } n^-(\alpha)
               + n^-(\alpha^{-1}) -
               { R^-(\alpha) \over R^+(\alpha) } n^-(\alpha^{-1})
       }\Bigr] g^-(\alpha^{-1}) \biggr\}
$$

$$
    = \alpha^{-2(k+1)} D(\alpha)^{-1}
      \biggl\{ A(\alpha)g^+(\alpha^{-1}) + B(\alpha)g^-(\alpha^{-1}) \biggr\}\, ,
    \eqno(3.25)
$$
where
$$ A(\alpha) = R^-(\alpha)R^-(\alpha^{-1}) \times
   \biggl\{ {R^+(\alpha)\over R^-(\alpha)}
            \Bigl( n^+(\alpha)-n^-(\alpha^{-1}) \Bigr)
              + n^-(\alpha^{-1})-n^-(\alpha)
   \biggr\}
$$
$$ = R^-(\alpha) R^-(\alpha^{-1}) \times
   \biggl\{ \Bigl( n^+(\alpha)-n^-(\alpha^{-1})\Bigr)
                      {R^+(\alpha)\over R^-(\alpha)}
            \Bigl( 1 - {R^-(\alpha)\over R^+(\alpha)} \Bigr)
          + \Bigl(n^+(\alpha)-n^-(\alpha)\Bigr)
   \biggr\},
$$
$$ B(\alpha) = R^+(\alpha)R^+(\alpha^{-1}) \times
   \biggl\{ {R^-(\alpha)\over R^+(\alpha)}
            \Bigl( n^-(\alpha)-n^+(\alpha^{-1}) \Bigr)
         + n^+(\alpha^{-1})-n^+(\alpha)
   \biggr\}
$$
$$ = R^+(\alpha)R^+(\alpha^{-1}) \times
    \biggl\{ \Bigl( n^-(\alpha)-n^+(\alpha^{-1}) \Bigr)
                {R^-(\alpha)\over R^+(\alpha)}
             \Bigl( 1 - {R^+(\alpha)\over R^-(\alpha)} \Bigr)
           + \Bigl( n^-(\alpha)-n^+(\alpha) \Bigr)
    \biggr\}=\overline {A(\alpha)}\,. \eqno(2.3.26)
$$
In view of theorem 2 and formulas (2.28), (2.30) we have the following
equalities on $\Omega_2^a$:
$$ { R^-(\alpha) \over R^+(\alpha) } = e^{2i\varphi_2(\alpha)} ; \quad
   R^+(\alpha)R^-(\alpha^{-1})
   = { \alpha-\alpha^{-1} \over N^+(\alpha) }
   = { \alpha-\alpha^{-1} \over \vert N^(\alpha)\vert } e^{-i\nu(\alpha)},
$$
with
$$\nu(\alpha) = \arg N^+(\alpha) =
               \cases{\eta(\alpha+\alpha^{-1}), \quad |\alpha|>1, \cr
                  \pi-\eta(\alpha+\alpha^{-1}), \quad |\alpha|<1. \cr }
$$
(We have taken into account that $|N^+(\alpha)| = |N^-(\alpha)| \equiv |N(\alpha)|$
for $\alpha \in {\bf R}$; further we also used
$|R^+(\alpha)| = |R^-(\alpha)| \equiv |R(\alpha)|$.)
It follows from these equalities that
$$
  R^-(\alpha)R^-(\alpha^{-1}) =
  R^+(\alpha)R^-(\alpha^{-1}) {R^-(\alpha) \over R^+(\alpha)}
  = { \alpha-\alpha^{-1} \over \vert N^(\alpha)\vert }
    e^{-i \Bigl( \nu(\alpha) - 2\varphi_2(\alpha) \Bigr) },
$$
$$ \Bigl(n^+(\alpha)-n^-(\alpha^{-1})\Bigr)
   {R^+(\alpha) \over R^-(\alpha)}
   = N^+(\alpha) e^{-2i\varphi_2(\alpha)}
   = \vert N(\alpha)\vert e^{i \Bigl( \nu(\alpha)-2\varphi_2(\alpha) \Bigr) },
$$
$$ {R^-(\alpha)\over R^+(\alpha)} - 1 = e^{2i\varphi_2(\alpha)}-1
   = \sin 2\varphi_2(\alpha) \Bigl( i-\tan\varphi_2(\alpha) \Bigr) \,.
$$

   Since
$\nu(\alpha)-2\varphi_2(\alpha) = \Bigl({1-s(\alpha)}\Bigr)
{\pi\over2} + s(\alpha) \chi_0(\alpha)\pi$ \footnote{In fact,
by the definition
$$\nu(t) = \cases{\eta(\alpha+\alpha^{-1}), \quad |\alpha|>1, \cr
                  \pi-\eta(\alpha+\alpha^{-1}), \quad |\alpha|<1, \cr }
$$
$$ 2 \varphi_2(\alpha)
   = s(\alpha) \Bigl(\eta(\alpha+\alpha^{-1})-\chi_0(\alpha)\pi\Bigr).
$$
for $|\alpha|>1$ we have
$\nu(\alpha)-2\varphi_2(\alpha) = \chi_0(\alpha)\pi$, i.e.~the formula is
true, and for $|\alpha|<1$ we have
$\nu(\alpha)-2\varphi_2(\alpha) = \pi - \chi_0(\alpha)\pi$, i.e.~the fornula
is also true.}, we have
$$ e^{i\Bigl(\nu(\alpha)-2\varphi_2(\alpha)\Bigr)}=
   e^{i\Bigl(1-s(\alpha)\Bigr) {\pi\over2} }
   e^{is(\alpha) \chi_0(\alpha) \pi}
   = s(\alpha) \cos\chi_0(\alpha)\pi;
$$
$$ \sin2\varphi_2(\alpha)
   = \sin \nu(\alpha) \cos \biggl\{ \Bigl({1-s(\alpha)}\Bigr) {\pi\over2}
                                       + s(\alpha) \chi_0(\alpha)\pi \biggr\}
   - \cos \nu(\alpha) \sin \biggl\{ \Bigl({1-s(\alpha)}\Bigr) {\pi\over2}
                                       + s(\alpha) \chi_0(\alpha)\pi \biggr\}
$$
$$ = \sin \nu(\alpha) \biggl\{ \cos \Bigl( { (1-s(\alpha)) {\pi\over2} }\Bigr)
                         \cos \Bigl( { s(\alpha) \chi_0(\alpha)\pi } \Bigr)
              - \sin \Bigl( { (1-s(\alpha)) {\pi\over2} }\Bigr)
                \sin \Bigl( { s(\alpha) \chi_0(\alpha)\pi } \Bigr) \biggr\}
   + 0
$$
$$ = \sin \nu(\alpha) \cos \Bigl( { (1-s(\alpha)) {\pi\over2} }\Bigr)
                      \cos \Bigl( { s(\alpha) \chi_0(\alpha)\pi } \Bigr)
    + 0 + 0
$$
$$ = \sin \nu(\alpha) \cdot s(\alpha) \cos \Bigl( { \chi_0(\alpha)\pi } \Bigr)
  = s(\alpha) \cos\chi_0(\alpha)\pi \sin\nu(\alpha).
$$
Hence,
$$ R^-(\alpha)R^-(\alpha^{-1}) =
   { \alpha - \alpha^{-1} \over \vert N(\alpha)\vert } s(\alpha)
   \cos\chi_0(\alpha)\pi\,,
$$
$$ \Bigl(n^+(\alpha)-n^-(\alpha^{-1})\Bigr) {R^+(\alpha)\over R^-(\alpha)}
   = \vert N(\alpha)\vert s(\alpha) \cos\chi_0(\alpha)\pi\,,
$$
$$ {R^-(\alpha)\over R^+(\alpha)}-1
   = s(\alpha) \cos\chi_0(\alpha)\pi \sin\nu(\alpha)
     \Bigl( i-\tan\varphi_2(\alpha) \Bigr).
$$
Therefore
$$ R^-(\alpha)R^-(\alpha^{-1}) =
   \vert R(\alpha)\vert \vert R(\alpha^{-1})\vert
   {\alpha-\alpha^{-1} \over \vert \alpha -\alpha^{-1} \vert}
   s(\alpha) \cos\chi_0(\alpha)\pi,
$$

$$ \Bigl(n^+(\alpha)-n^-(\alpha^{-1})\Bigr) {R^+(\alpha)\over R^-(\alpha)}
   \Bigl({R^-(\alpha)\over R^+(\alpha)}-1\Bigr)
$$
$$ = \vert N(\alpha)\vert s(\alpha) \cos \Bigl( \chi_0(\alpha)\pi \Bigr)
                          s(\alpha) \cos \Bigl( \chi_0(\alpha)\pi \Bigr)
     \sin\nu(\alpha) \Bigl( i-\tan\varphi_2(\alpha) \Bigr)
$$
$$ = \vert N(\alpha)\vert
       \sin\nu(\alpha)\Bigl(i-\tan\varphi_2(\alpha)\Bigr)
   = \Bigl(i-\tan\varphi_2(\alpha)\Bigr)\, {\rm Im}\, N^+(\alpha).
$$
Substituting these expressions into (3.26), we obtain
$$ A(\alpha)=\vert R(\alpha)\vert\vert R(-\alpha)\vert
     {\alpha-\alpha^{-1} \over \vert \alpha -\alpha^{-1} \vert}
     s(\alpha) \cos\chi_0(\alpha) \pi
  \biggl\{ -\Bigl(i-\tan\varphi_2(\alpha)\Bigr)\, {\rm Im}\, N^+(\alpha) +
      2i\,{\rm Im}\,n^+(\alpha) \biggr\}.
$$
Finally, from the definition of $N(z)$, formulas (2.3), (2.4) and
the absolute continuity of the functions $\rho_R(\tau+\tau^{-1})$,
$\rho_L(\tau+\tau^{-1})$ on $\Omega_2^a$ it follows that
$$ {\rm Im}\, N^+(\alpha)= \pi\Bigl(\sigma(\alpha)+\sigma(\alpha^{-1})\Bigr),
   \qquad {\rm Im}\,n^+(\alpha)=\pi\sigma(\alpha), $$
and, since
$$ \tan\varphi_2(\alpha) = s(\alpha)
   \tan{1\over2} \Bigl(\eta(\alpha+\alpha^{-1})-\chi_0(\alpha)\pi \Bigr)=
   s(\alpha) \cos\chi_0(\alpha)\pi\tan\vert\varphi_2(\alpha)\vert , \
\footnote{In fact, if $\chi_0(\alpha)=0$, then the argument of
${1\over2} (\eta(\alpha+\alpha^{-1})-\chi_0(\alpha)\pi)$ lies in $(0,\pi)$
and
$$ {1\over2} (\eta(\alpha+\alpha^{-1})-\chi_0(\alpha)\pi) =
   \vert {1\over2} (\eta(\alpha+\alpha^{-1})-\chi_0(\alpha)\pi) \vert =
   \vert  {1\over2} (\eta(\alpha+\alpha^{-1}) \vert .
$$
And if $\chi_0(\alpha)=1$, then the argument of
${1\over2} (\eta(\alpha+\alpha^{-1})-\chi_0(\alpha)\pi)$ lies in $(-\pi,0)$
and
$$ \tan {1\over2} (\eta(\alpha+\alpha^{-1})-\chi_0(\alpha)\pi) =
   - \tan \vert{1\over2} (\eta(\alpha+\alpha^{-1})-\chi_0(\alpha)\pi) \vert\,.
$$
}
$$
we obtain the equality
$$ A(\alpha) = \vert R(\alpha)\vert \vert R(\alpha^{-1})\vert
      {\alpha-\alpha^{-1} \over \vert \alpha -\alpha^{-1} \vert}
      s(\alpha) \cos\Bigl(\chi_0(\alpha)\pi\Bigr)
   \biggl\{ - \Bigl(i-\tan\varphi_2(\alpha)\Bigr) \pi
          (\sigma(\alpha) + \sigma(\alpha^{-1})) + 2i \pi \sigma(\alpha) \biggr\}
$$

$$ = \vert R(\alpha)\vert \vert R(\alpha^{-1})\vert
   {\alpha-\alpha^{-1} \over \vert \alpha -\alpha^{-1} \vert}
   s(\alpha) \pi \cos\Bigl(\chi_0(\alpha)\pi\Bigr)
$$
$$ \times \biggl\{ s(\alpha) \cos\Bigl(\chi_0(\alpha)\pi\Bigr)
                                 \tan|\varphi_2(\alpha)|
     \Bigl( \sigma(\alpha) + \sigma(\alpha^{-1}) \Bigr) -
     i \Bigl( \sigma(\alpha^{-1}) - \sigma(\alpha) \Bigr) \biggr\}
$$

$$ = \vert R(\alpha)\vert \vert R(\alpha^{-1})\vert
   {\alpha-\alpha^{-1} \over \vert \alpha -\alpha^{-1} \vert}
   s(\alpha) \pi
$$
$$ \times
   \biggl\{ s(\alpha) \tan|\varphi_2(\alpha)| (\sigma(\alpha)
                   + \sigma(\alpha^{-1}))
     - i \Bigl( \sigma(\alpha^{-1}) - \sigma(\alpha) \Bigr)
              \cos\Bigl(\chi_0(\alpha)\pi\Bigr) \biggr\}
$$

$$ = \vert R(\alpha)\vert \vert R(\alpha^{-1})\vert
   {\alpha-\alpha^{-1} \over \vert \alpha -\alpha^{-1} \vert}
   s(\alpha) \pi \sqrt{\sigma(\alpha)\sigma(\alpha^{-1})}
$$
$$ \times \biggl\{ s(\alpha) \tan|\varphi_2(\alpha)|
   \Bigl( \sqrt{\sigma(\alpha) \over \sigma(\alpha^{-1})}
        + \sqrt{\sigma(\alpha^{-1}) \over \sigma(\alpha)} \Bigr)
     - i \Bigl( \sqrt{\sigma(\alpha^{-1}) \over \sigma(\alpha)}
        - \sqrt{\sigma(\alpha) \over \sigma(\alpha^{-1})} \Bigr)
           \cos\Bigl(\chi_0(\alpha)\pi\Bigr)  \biggr\}
$$
$$ = \vert R(\alpha)\vert\vert R(\alpha^{-1})\vert \pi
     \sqrt{\sigma(\alpha)\sigma(\alpha^{-1})}
    \Bigl\{ { {\alpha-\alpha^{-1} \over \vert \alpha -\alpha^{-1} \vert}
         p(\alpha) - i {\alpha \over \vert \alpha \vert} m(\alpha)}\Bigr\}\,,
$$
where the functions $p(\alpha)$ and $m(\alpha)$ are defined by
(3.21$^\prime$), (3.23$^\prime$).

    Substitute the obtained expressions of $D(\alpha)$, $A(\alpha)$ and
$B(\alpha)=\overline{A(\alpha)}$ into (3.25):
$$ g^+(\alpha)-g^-(\alpha) =
   \alpha^{-2(k+1)} \cdot
   { 1 \over 2 \pi i { \vert R(\alpha^{-1})\vert }^2 \sigma(\alpha^{-1}) }
$$
$$ \times \biggl\{ \vert R(\alpha)\vert \vert R(\alpha^{-1})\vert \pi
     \sqrt{\sigma(\alpha)\sigma(\alpha^{-1})}
    \Bigl[ { {\alpha-\alpha^{-1} \over \vert \alpha -\alpha^{-1} \vert}
           p(\alpha) - i {\alpha \over \vert \alpha \vert} m(\alpha)}\Bigr]
        g^+(\alpha^{-1}) +
$$
$$  + \vert R(\alpha)\vert \vert R(\alpha^{-1})\vert \pi
     \sqrt{\sigma(\alpha)\sigma(\alpha^{-1})}
    \Bigl[ { {\alpha-\alpha^{-1} \over \vert \alpha -\alpha^{-1} \vert}
           p(\alpha) + i {\alpha \over \vert \alpha \vert} m(\alpha)}\Bigr]
        g^-(\alpha^{-1}) \biggr\}
$$

$$ = - { \alpha^{-2(k+1)} \over 2i } \cdot
   { \vert R(\alpha)\vert \over \vert R(\alpha^{-1}) }
     \sqrt{ \sigma(\alpha) \over \sigma(\alpha^{-1}) }
$$
$$ \times \biggl\{
   { - {\alpha-\alpha^{-1} \over \vert \alpha -\alpha^{-1} \vert}
      p(\alpha) \Bigl( g^+(\alpha^{-1}) + g^-(\alpha^{-1}) \Bigr)
     + i {\alpha \over \vert \alpha \vert} m(\alpha)
      \Bigl( g^+(\alpha^{-1}) - g^-(\alpha^{-1}) \Bigr) }
   \biggr\} \,,
$$
or
$$ - 2i q(\alpha) \alpha^{2(k+1)}
    \Bigl( g^+(\alpha) - g^-(\alpha) \Bigr)
$$
$$
  = - p(\alpha) {\alpha-\alpha^{-1} \over \vert \alpha -\alpha^{-1} \vert}
      \Bigl( g^+(\alpha^{-1}) + g^-(\alpha^{-1}) \Bigr)
     + i {\alpha \over \vert \alpha \vert} m(\alpha)
      \Bigl( g^+(\alpha^{-1}) - g^-(\alpha^{-1}) \Bigr),
$$
which is equivalent to (3.20). \hfill\rule{0.5em}{0.5em}

\smallskip

   Note that the function $\varphi_2(\alpha)$
on the set $\Omega_2^a$ is devided of $0$, $\pi\over2$, $\pi$, $3\pi\over2$,
from which it follows that $\tan \varphi_2(\alpha)$ is bounded and bounded
away from zero. Since $\mu(\alpha+\alpha^{-1})$ is also bounded and bounded
away from zero, on this set, the function $p(\alpha)$ is also bounded and
bounded away from zero.

\bigskip

{\bf Calculation of the integrals along the contour $\Gamma_3$}

\medskip

{\bf Lemma 3.6.} {\it For $\theta \in [-\pi,\pi]$ the functions
$ g \Bigl(k, (1\mp h)e^{i\theta} \Bigr) $ have limits
$g^\pm(k,e^{i\theta}) \in L^2[-\pi,\pi]$ as $h \downarrow 0$.
The following equality is true for the integral along the contour
$\Gamma_3$ in (3.2):}
$$ {1 \over 2\pi i} \int\limits_{\Gamma_3}{ g(\zeta) \over \zeta-z } d\zeta =
   {1 \over 2\pi i} \int\limits_{|\xi|=1}
                    { g^+(\xi) - g^-(\xi) \over \xi-z } d\xi
   = {1 \over 2\pi } \int\limits_{-\pi}^\pi { g^+(e^{i\theta}) - g^-(e^{i\theta})
              \over e^{i\theta}-z } e^{i\theta} d\theta.
   \eqno (3.27)
$$

P r o o f.\  We remind that the contour $\Gamma_3$ consists of two
concentric circles of radius $(1-\delta)$ and $(1+\delta)$ with the center
at the point $0$. Let us deform the contour $\Gamma_3$ in the following way:
we deform its part that lies in the $\varepsilon$-neighborhoods
($\varepsilon>\delta$) of the points $-1,1$ by the arcs of radius
$\varepsilon$ with the center at these points. Tending $\delta\downarrow 0$
and using the methods analogous to those of the previous subsection
we obtain that
$$ {1 \over 2\pi i} \int\limits_{\Gamma_3}{ g(\zeta) \over \zeta-z } d\zeta =
   {1 \over 2\pi i} \int\limits_{ \small \matrix{|\xi|=1, \cr
                                 |\xi-1|> \varepsilon, \cr
                                 |\xi+1|> \varepsilon \cr } }
                                { g^+(\xi) - g^-(\xi) \over \xi-z } d\xi
  + {1 \over 2\pi i} \int\limits_{|\zeta-1| = \varepsilon }
                     { g(\zeta) \over \zeta-z } d\zeta
  + {1 \over 2\pi i} \int\limits_{|\zeta+1| = \varepsilon }
                               { g(\zeta) \over \zeta-z } d\zeta,
$$
if
$$ \lim_{h\downarrow 0} \Vert g \Bigl( (1\mp h)e^{i\theta} \Bigr) - g^\pm(e^{i\theta})
   \Vert_{L^p[\alpha(\varepsilon),\pi-\alpha(\varepsilon)]} = 0,
   \eqno(3.28)
$$
$$ \lim_{h\downarrow 0} \Vert g \Bigl((1\mp h)e^{i\theta} \Bigr) - g^\pm(e^{i\theta})
   \Vert_{L^p[-\pi + \alpha(\varepsilon),-\alpha(\varepsilon)]} = 0,
   \eqno(3.29)
$$
with $p\ge 1$, where $\alpha(\varepsilon)$ is such a small number that
$e^{i\alpha(\varepsilon)}$ lies in the $\varepsilon$-neighborhood of the
point $1$. Tending $\varepsilon \downarrow 0$ in the last formula
we obtain that
$$ {1 \over 2\pi i} \int\limits_{\Gamma_3}{ g(\zeta) \over \zeta-z } d\zeta =
   {1 \over 2\pi i} \int\limits_{|\xi|=1}
                    { g^+(\xi) - g^-(\xi) \over \xi-z } d\xi
   = {1 \over 2\pi } \int\limits_{-\pi}^\pi { g^+(e^{i\theta}) - g^-(e^{i\theta})
              \over e^{i\theta}-z } e^{i\theta} d\theta,
$$
if
$$ \lim_{\varepsilon \downarrow 0} \int\limits_{|\zeta \pm 1| = \varepsilon }
                      { g(\zeta) \over \zeta-z } d\zeta = 0 \,,
$$
$$ g^\pm(e^{i\theta}) \in L^p[-\pi,\pi], \quad p\ge 1,
$$
and also if for all small $\varepsilon > 0$ the conditions (3.28) and (3.29)
are fulfilled.

   Without restriction of generality, we will prove these facts
for the interior part of the disk, i.e. we will prove that
$$ \lim_{\varepsilon\downarrow 0} \int\limits_{ |\zeta \pm 1| = \varepsilon, \
                                            |\zeta| < 1  }
                      { g(\zeta) \over \zeta-z } d\zeta = 0 \,,
$$
$$ g^+(e^{i\theta}) \in L^2[-\pi,\pi],
$$
$$ \lim_{h\downarrow 0} \Vert g \Bigl( (1 - h) e^{i\theta} \Bigr) - g^+(e^{i\theta})
   \Vert_{L^p[\alpha(\varepsilon),\pi-\alpha(\varepsilon)]} = 0,
   \eqno(3.28')
$$
$$ \lim_{h\downarrow 0} \Vert g \Bigl((1 - h)e^{i\theta} \Bigr) - g^+(e^{i\theta})
   \Vert_{L^p[-\pi + \alpha(\varepsilon),-\alpha(\varepsilon)]} = 0.
   \eqno(3.29')
$$

   As it is known from theorem 2, for the function $R(z)$
the representation
$$ R(z) = C \, R_\mu(z) R_3(z) \, R_0(z) R_1(z) R_2(z)
          = R_\mu(z) R^{(1)}(z) R^{(2)}(z) \,,
$$
is true, where $ R^{(1)}(z) = R_3(z) =
      \hat P\Bigl(z, - {\hat\gamma_3\over2} \Bigr)$, and
$R^{(2)}(z) = C \, R_0(z) R_1(z) R_2(z)$ is a holomorphic function on the
unit circle. Let us derive an analogous decomposition for the function $n(z)$.
It can be shown (see this method, for example, in [15]) that inside the
unit disk
$$ n(z) = n_0(z) n_1(z)\,, \quad |z|<1,
$$
where
$$ n_0(z) = \hat P (z, \hat\eta^R)\,,
$$
$$ \hat\eta^R(\theta) =
   \cases{ - \eta^R(e^{i\theta}+e^{-i\theta}) + {\pi \over 2}, \quad 0<\theta<\pi, \cr
            \eta^R(e^{i\theta}+e^{-i\theta}) - {\pi \over 2}, \quad -\pi<\theta<0, \cr}
$$
and $n_1(z)$ is a holomorphic on the circle function, which vanish
at the points $-1,1$, i.e.
$$ n(z) = (z-1) (z+1) n^{(1)}(z) n^{(2)}(z)\,,
$$
where $ n^{(1)}(z) = n_0(z) = \hat P(z, \hat\eta^R)$, and
$n^{(2)}(z)$ is holomorphic on the circle. Thus, in view of (3.1),
$$ g(\zeta) = g^{(1)}(\zeta) R_\mu(\zeta) R^{(1)}(\zeta) +
            g^{(2)}(\zeta) R_\mu(\zeta) R^{(1)}(\zeta) n^{(1)}(\zeta) (\zeta-1)(\zeta+1)\,,
   \eqno(3.30)
$$
where the functions $ g^{(1)}(\zeta)$, $g^{(2)}(\zeta)$ are holomorphic on
the unit circle. Since the functions $\eta(\tau)$ and $\eta^R(\tau)$ satisfy
the H\"older condition under $(-2,2)$, the functions
$R^{(1)}((1-h)e^{i\theta}) = \hat P((1-h)e^{i\theta},-{\hat\gamma_3\over2})$,
$R^{(1)}((1-h)e^{i\theta}) n^{(1)}((1-h)e^{i\theta}) (e^{i\theta}-1)
  = \hat P((1-h)e^{i\theta},-{\hat\gamma_3\over2}+\hat\eta^R) (e^{i\theta}-1)$,
(as functions of $\theta$), converge uniformly when $h\to \pm 0$ to their
limits on the segments $[\alpha(\varepsilon),\pi-\alpha(\varepsilon)]$,
$[-\pi + \alpha(\varepsilon),-\alpha(\varepsilon)]$. All the more,
they converge in $L^2$-norm, and their limits also belong to $L^2$. Further,
according to lemma 2.5, the function $R_\mu((1-h)e^{i\theta})$ converges
uniformly to its limit $R_\mu^+(e^{i\theta})$ as $h\downarrow 0$.
That is why all the function $g((1-h)e^{i\theta})$ converges uniformly
as $h\downarrow 0$ to its limit $g^+(e^{i\theta})$ on segments
$[\alpha(\varepsilon),\pi-\alpha(\varepsilon)]$,
$[-\pi + \alpha(\varepsilon),-\alpha(\varepsilon)]$, i.e.~the
properties (3.28$^\prime$) and (3.29$^\prime$) are proven.
Without loss of generality, we will prove the two other properties
only in a neighborhood of the point $1$, i.e.~we will prove that
$$ \lim_{\varepsilon\downarrow 0}
    \int\limits_{ |\zeta - 1| = \varepsilon, \ |\zeta|<1 }
                      { g(\zeta) \over \zeta-z } d\zeta  = 0 \,,
    \qquad
   g^+(e^{i\theta}) \in L^2[-{\pi\over2},{\pi\over2}]\,.
$$

   We remind that the functions $\hat\gamma_3(\theta)$, $\hat\eta^R(\theta)$
are odd and have a jump at the point $\theta=0$. Let
$$ \hat\eta_1^R(\theta) =
   \cases{ \hat\eta^R(+0), \quad 0<\theta<\pi, \cr
           \hat\eta^R(-0), \quad -\pi<\theta<0, \cr}
$$
$$ \hat\eta_2^R(\theta) = \hat\eta^R(\theta) - \hat\eta_1^R(\theta);
$$
$$ \hat\nu_1(\theta) =
   \cases{ -{\hat\gamma_3(+0)\over2}, \quad 0<\theta<\pi, \cr
           -{\hat\gamma_3(-0)\over2}, \quad -\pi<\theta<0, \cr}
$$
$$ \hat\nu_2(\theta) = -{\hat\gamma_3(\theta)\over2} - \hat\nu_1(\theta).
$$

    The functions $\hat\nu_2(\theta)$, $\hat\eta_2^R(\theta)$ satisfy the
H\"older condition under $(-{\pi\over2},{\pi\over2})$.
That is why the functions
$\hat P \Bigl( (1-h)e^{i\theta}), \hat\nu_2 \Bigr)$,
$\hat P \Bigl( (1-h)e^{i\theta}),  \hat\eta_2^R \Bigr)$
converge uniformly to their limits
$\hat P^+ (e^{i\theta}, \hat\nu_2)$,
$\hat P^+ (e^{i\theta},  \hat\eta_2^R)$ for
$ -{\pi\over2} < \theta < {\pi\over2} $,
in particular, they are uniformly bounded in a neighborhood of $\theta=0$.
Further,
$$ \hat P (z,  \hat\eta_1^R) = C
   \Bigl( {z-1 \over z+1} \Bigr)^{ 2\hat\eta^R(+0) \over \pi }\,, \qquad
  \hat P (z, \hat\nu_1) = C
   \Bigl( {z-1 \over z+1} \Bigr)^{ - {\hat\gamma_3(+0) \over \pi }}\,.
$$

   (In fact, looking at the limit values of the arguments of the function
${z-1 \over z+1}$ on the circle, it can be seen that, to within a constant
factor,
$$\hat P (z, {\theta \over |\theta|} \cdot {\pi \over 2} )
  = {z-1 \over z+1} ,
$$
and, since
$$ \hat\eta^R_1(\theta)= {2 \hat\eta^R(+0) \over \pi}
   \Bigl( {\theta \over |\theta|} \cdot {\pi \over 2} \Bigr) \,,  \qquad
   \hat\nu_1(\theta)= - {\hat\gamma_3(+0) \over \pi}
   \Bigl( {\theta \over |\theta|} \cdot {\pi \over 2} \Bigr) \,,
$$
then
$$   \hat P (z,  \hat\eta_1^R) =
   \hat P (z, {2 \hat\eta^R(+0) \over \pi}
   \Bigl( {\theta \over |\theta|} \cdot {\pi \over 2} \Bigr)),  \qquad
  \hat P (z, \hat\nu_1) =
   \hat P (z, - {\hat\gamma_3(+0) \over \pi}
   \Bigl( {\theta \over |\theta|} \cdot {\pi \over 2} \Bigr) ) ,
$$
from where we have the necessary identities.)

   Thus,
$$ n^{(1)} \Bigl( (1-h)e^{i\theta} \Bigr) \Bigl( (1-h)e^{i\theta} - 1 \Bigr) =
   \hat P ((1-h)e^{i\theta}, \hat\eta_1^R)
{   { ((1-h)e^{i\theta} - 1) }^{ 1 + {2\hat\eta^R(+0) \over \pi} } \over
   { ((1-h)e^{i\theta} + 1) }^{2\hat\eta^R(+0) \over \pi}
}\,.
$$
Since $ -{\pi \over 2} < \eta^R(+0) < {\pi \over 2} $, the function
$ ((1-h)e^{i\theta} - 1)^{1 + {2\hat\eta^R(+0) \over \pi} } $
is uniformly bounded in a neighborhood of $\theta=0$ and converges uniformly
as $h\downarrow 0$ to its limit
$(e^{i\theta} - 1)^{ 1 + {2\hat\eta^R(+0) \over \pi} }$. It means that
$ n^{(1)} \Bigl( (1-h)e^{i\theta} \Bigr) \Bigl( (1-h)e^{i\theta} - 1 \Bigr) $
for $ -{\pi \over 2} < \theta < {\pi \over 2} $ is also uniformly bounded and
uniformly converges to its limit (which belongs to $L^2$). According to lemma
2.5, the function $ R_\mu((1-h)e^{i\theta})$ has the same property.
That is why, according to decomposition (3.30), to prove that
$ g^+(e^{i\theta}) \in L^2[-{\pi\over2},{\pi\over2}]$ we have to prove
that ${R^{(1)}}^+(e^{i\theta}) \in L^2[-{\pi\over2},{\pi\over2}]$.
But this fact follows from what
$$ {R^{(1)}}^+((1-h)e^{i\theta})
    = \hat P((1-h)e^{i\theta},\hat\nu_1) \hat P((1-h)e^{i\theta},\hat\nu_2)
$$
$$  = \Bigl( {(1-h)e^{i\theta}-1 \over (1-h)e^{i\theta}+1} \Bigr)^{ - {\hat\gamma_3(+0) \over \pi }}
        \hat P((1-h)e^{i\theta},\hat\nu_2) \,.
$$

   As it was already said, the second one of the functions of the last product,
$\hat P((1-h)e^{i\theta},\hat\nu_2)$, is uniformly bounded and converges
uniformly as $h\downarrow 0$. As for the first function,
$ \Bigl( {z-1 \over z+1} \Bigr)^{ - {\hat\gamma_3(+0) \over \pi }}$,
in view of condition B) we have $0<\eta(2-0)<\pi$, from where
$ -{\pi \over 2} < \hat\gamma_3(+0) < {\pi \over 2} $. Hence, the function
$ \Bigl( {e^{i\theta}-1 \over e^{i\theta}+1}
  \Bigr)^{ - {\hat\gamma_3(+0) \over \pi }}$ is square-summable for
$ -{\pi \over 2} < \theta < {\pi \over 2} $, from where we deduce
the belonging ${R^{(1)}}^+(e^{i\theta}) \in L^2[-{\pi\over2},{\pi\over2}]$.
Therefore,
$ g^+(e^{i\theta}) \in L^2[-{\pi\over2},{\pi\over2}]$.

\bigskip

   Let us now pass to the property
$$ \lim_{\varepsilon\downarrow 0}
    \int\limits_{ |\zeta - 1| = \varepsilon, \ |\zeta|<1 }
                      { g(\zeta) \over \zeta-z } d\zeta  = 0 \,.
$$
As it may be seen,
$$ g(\zeta) = g^{(3)}(\zeta)
   \Bigl( { \zeta-1 \over \zeta+1 } \Bigr)^{ -\hat\gamma_3(+0)\over\pi} \,,
$$
where the function
$$ g^{(3)}(\zeta) = \Bigl( g^{(1)}(\zeta) +
   g^{(2)}(\zeta) n^{(1)}(\zeta) (\zeta-1) (\zeta+1) \Bigr)
   R_\mu(\zeta) \hat P(\zeta, \hat\nu_2)
$$
is uniformly bounded in a neighborhood of $1$ (we have just proved it
individually for $n^{(1)}(\zeta) (\zeta-1)$, $R_\mu(\zeta)$ and
$\hat P(\zeta, \hat\nu_2)$). Hence,
$$ \vert \int\limits_{ |\zeta - 1| = \varepsilon, \ |\zeta|<1 }
          { g(\zeta) \over \zeta-z } d\zeta \vert           \leq
  {\max_\zeta} \Bigl| { g^{(3)}(\zeta) \over
                    (\zeta-z)(\zeta+1)^{ - \hat\gamma_3(+0) \over \pi } } \Bigr| \
  \int\limits_{ |\zeta - 1| = \varepsilon, \ |\zeta|<1 }
      |\zeta-1|^{ -\hat\gamma_3(+0) \over \pi }   |d\zeta|
$$
$$  \leq   {\max_\zeta} \vert { g^{(3)}(\zeta) \over
                    (\zeta-z)(\zeta+1)^{ -\hat\gamma_3(+0) \over \pi } } \vert \
           \varepsilon^{ 1 - |{\hat\gamma_3(+0) \over \pi}| } \to 0,
          \quad \varepsilon \downarrow 0,
$$
because $| { - \hat\gamma_3(+0) \over \pi }| <1$, which complete the
proof of the lemma.
 \hfill\rule{0.5em}{0.5em}

\medskip

{\bf Lemma 3.7.} {\it Almost everywhere on the unit circle
$|\xi|=1$ the following equality holds
$$ \xi^{2(k+1)} (g^+(k,\xi)-g^-(k,\xi))
   =  - \hat r(\xi) (g^+(k,\xi^{-1})+g^-(k,\xi^{-1})) \,, \eqno(3.31)
$$
where
$$  \hat r(\xi) = - { 1 \over 1 + q(\xi) } \cdot
     { n^+(\xi) - n^-(\xi)  \over  \xi -\xi^{-1} }
    R^+(\xi) R^+(\xi^{-1}) \,,           \eqno(3.25')
$$
$$ q(\xi) = 2\pi { \sqrt{ {\rho_R}^\prime (\xi+\xi^{-1}) {\rho_L}^\prime (\xi+\xi^{-1})} \over
                       |n^+(\xi) - n^-(\xi^{-1})| } \,.   \eqno(3.25'')
$$
Moreover, }
$$ |\hat r(\xi)| < C < 1. $$

\medskip

   P r o o f.\  The functions
$\psi^-(k,\xi)$, $\psi^+(k,\xi)$ (the limit values of the Weyl solutions
on the unit circle) are the solutions (with respect to $k$) of
finite-difference equation (0.5). Hence, they are linear combinations
of the functions $\psi^+(k,\xi)$,
$\psi^+(k,\xi^{-1}) = \overline{\psi^+(k,\xi)}$ and
$\psi^-(k,\xi)$, $\psi^-(k,\xi^{-1}) = \overline{\psi^-(k,\xi)}$,
respectively:
$$\psi^-(k,\xi) = a(\xi) \psi^+(k,\xi) + b(\xi) \psi^+(k,\xi^{-1}) \,, \eqno(3.34)
$$
$$\psi^+(k,\xi) = c(\xi) \psi^-(k,\xi) + d(\xi) \psi^-(k,\xi^{-1}) \,. \eqno(3.35)
$$

1. Let us first find the expression for the functions
$a(\xi)$, $b(\xi)$ и $c(\xi)$, $d(\xi)$ and show that
$\vert {b(\xi) \over a(\xi)} \vert \leq 1$,
$\vert {d(\xi) \over c(\xi)} \vert \leq 1$. We remind that
$$ \psi(k,z)=n(z) P_k(z+z^{-1}) + Q_k(z+z^{-1}) \,.
$$
Substituting this expression in (3.34) and (3.35) for $k=0, k=-1$,
and taking into account the initial data on $Q_k(z)$ and $P_k(z)$,
we obtain two systems of linear equation. Solving these linear systems, we
find
$$ a(\xi) = { n^+(\xi^{-1}) - n^-(\xi) \over n^+(\xi^{-1}) - n^+(\xi) } \,,
   \qquad
   b(\xi) = { n^+(\xi) - n^-(\xi) \over n^+(\xi) - n^+(\xi^{-1}) } \,, \eqno(3.36)
$$
$$ c(\xi) = { n^-(\xi^{-1}) - n^+(\xi) \over n^-(\xi^{-1}) - n^-(\xi) }\,,
  \qquad
   d(\xi) = { n^-(\xi) - n^+(\xi) \over n^-(\xi) - n^-(\xi^{-1}) }\,.  \eqno(3.37)
$$

   Further,
$$ { b(\xi) \over a(\xi) } =
   { n^+(\xi) - n^-(\xi) \over n^+(\xi) - n^+(\xi^{-1}) } \cdot
   { n^+(\xi^{-1}) - n^+(\xi) \over n^+(\xi^{-1}) - n^-(\xi) } =
   - { n^+(\xi) - n^-(\xi) \over n^+(\xi^{-1}) - n^-(\xi) } =
   - { n^+(\xi) - n^-(\xi) \over \overline{n^+(\xi)} - n^-(\xi) } \,.
$$
The real parts of the numerator and the denominator in the last expression
are equal. Thus, in order to verify that
$\vert {b(\xi) \over a(\xi)} \vert \leq 1$, it is sufficient to show that
$ \vert {\rm Im}\, ( n^+(\xi) - n^-(\xi) ) \vert  \leq
  \vert {\rm Im}\, ( \overline{n^+(\xi)} - n^-(\xi) ) \vert $.
But it is an evident consequence of the fact that the signs of the
imaginary parts of $n^+(\xi)$ and $n^-(\xi)$, \ $|\xi|=1$, are the same.
In the same way we prove that the absolute value of the expression
$$ {d(\xi) \over c(\xi)} =
    - { n^-(\xi) - n^+(\xi) \over n^-(\xi^{-1}) - n^+(\xi) } \,
$$
is less or equal than $1$.

\bigskip

   2. It is seen, that the coefficients $a(\xi)$, $b(\xi)$, $c(\xi)$, $d(\xi)$
are defined by the limit values of the function $n(z)$ at the points $\xi$,
$\xi^{-1}$ and are not independent. To clarify the connection between the
coefficients in equalities (3.34) and (3.35), we introduce the notation:
$$ v^\pm(\xi) = {\rm Im}\, n^\pm(\xi)
   = {n^\pm(\xi)-n^\pm(\xi^{-1}) \over 2i} \,.
$$
Taking into account that the signs of ${\rm Im}\, n^+(\xi)$
and ${\rm Im}\, n^-(\xi)$ are the same, we can rewrite the equality
(3.28) in the form
$$  v^+(\xi) R^+(\xi) R^+(\xi^{-1}) =
    v^-(\xi) R^-(\xi) R^-(\xi^{-1}) \,.
$$
In these notations
$$ {1 \over a(\xi)} =
   { n^+(\xi) - n^+(\xi^{-1}) \over n^-(\xi) - n^+(\xi^{-1}) } =
   { 2i \, v^+(\xi) \over N^-(\xi) } =
   { 2i  \over \xi-\xi^{-1} } v^+(\xi) \cdot
   { \xi-\xi^{-1} \over N^-(\xi) } \,,
$$
$$ {1 \over c(\xi)} =
   { n^-(\xi) - n^-(\xi^{-1}) \over n^+(\xi) - n^-(\xi^{-1}) } =
   { 2i \, v^-(\xi) \over N^+(\xi) } =
   { 2i  \over \xi-\xi^{-1} } v^-(\xi) \cdot
   { \xi-\xi^{-1} \over N^+(\xi) } \,,
$$

   In view of our factorization theorem the last expression can be rewritten
in the form
$$ {1 \over a(\xi)} =
   { 2i \over \xi-\xi^{-1} } v^+(\xi) \cdot
   R^-(\xi) R^+(\xi^{-1}) =
   { 2i \over \xi-\xi^{-1} } v^+(\xi) \cdot
   R^+(\xi) R^+(\xi^{-1}) {R^-(\xi) \over R^+(\xi)} =
   q(\xi) {R^-(\xi) \over R^+(\xi)} \,,
$$
$$ {1 \over c(\xi)} =
   { 2i \over \xi-\xi^{-1} } v^-(\xi) \cdot
   R^+(\xi) R^-(\xi^{-1}) =
   { 2i \over \xi-\xi^{-1} } v^-(\xi) \cdot
   R^+(\xi) R^-(\xi^{-1}) {R^+(\xi) \over R^-(\xi)} =
   q(\xi) {R^+(\xi) \over R^-(\xi)} \,,
$$
with
$$ q(\xi) =
   { 2i \over \xi-\xi^{-1} } v^+(\xi) \cdot
   R^+(\xi) R^+(\xi^{-1}) =
   { 2i \over \xi-\xi^{-1} } v^-(\xi) \cdot
   R^+(\xi) R^-(\xi^{-1}) \,.
$$

   Let us analyze the coefficient $q(\xi)$. We remind, that, due to the
definition of $R(z)$, we have $R^\pm(\xi^{-1})=\overline{R^\pm(\xi)}$.
Besides, $v^\pm(\xi)>0, \ {\rm Im} \, \xi >0$, and
$v^\pm(\xi)<0, \ {\rm Im} \, \xi <0$.
That is why for all $\xi\ne \pm1$
$$ q(\xi)>0.$$
We also remind that $R(z)= R_{0123}(z)R_\mu(z)$, where
$R^+_{0123}(\xi) = R^-_{0123}(\xi)$, and also
$R_\mu^+(\xi)=\sqrt{v^-(\xi)\over v^+(\xi)} R_\mu^-(\xi)$. Thus,
$$ q(\xi)
   = 2 \, \vert { v^+(\xi) \over \xi-\xi^{-1} } \vert \,
     R_{0123}^+(\xi) R_\mu^+(\xi) R_{0123}^+(\xi^{-1}) R_\mu^+(\xi^{-1})
$$
$$   = 2 \sqrt{v^+(\xi)v^-(\xi)} \,
        \vert { R_{0123}^-(\xi) R_\mu^-(\xi) R_{0123}^+(\xi^{-1}) R_\mu^+(\xi^{-1})
         \over \xi-\xi^{-1} } \vert
$$
$$ =  2 { \sqrt{ | {\rm Im} \, n^+(\xi) \ {\rm Im} \, {n^-(\xi)} | } \over
                       |N^-(\xi)| } |
   =  2 \pi { \sqrt{\rho_+'(\xi+\xi^{-1}) \rho_-'(\xi+\xi^{-1})} \over
                       |n^+(\xi) - n^-(\xi^{-1})| }. \eqno(3.38)
$$

3. We obtain from equality (3.36)
$$ {1\over a(\xi)} \psi^-(k,\xi) = \psi^+(k,\xi) + {b(\xi)\over a(\xi)} \psi^+(k,\xi^{-1}) \,,
$$
which is equivalent, as we have proved, to
$$ q(\xi) {R^-(\xi) \over R^+(\xi)} \psi^-(k,\xi)
       = \psi^+(k,\xi) + {b(\xi)\over a(\xi)} \psi^+(k,\xi^{-1}) \,,
$$
or
$$ q(\xi) R^-(\xi) \psi^-(k,\xi) = R^+(\xi) \psi^+(k,\xi) + {b(\xi)\over a(\xi)}
   {R^+(\xi) \over R^+(\xi^{-1})} R^+(\xi^{-1}) \psi^+(k,\xi^{-1}) \,. \eqno(3.39)
$$
In the same way, we have from (3.35) that
$$ q(\xi) R^+(\xi)  \psi^+(k,\xi) = R^-(\xi) \psi^-(k,\xi) + {d(\xi) \over c(\xi)}
   {R^-(\xi) \over R^-(\xi^{-1})} R^-(\xi^{-1}) \psi^-(k,\xi^{-1}) \,. \eqno(3.40)
$$

   To reduce the two last equalities to the necessary form let us verify
the following relation for the coefficients in the right-hand side:
$$ {b(\xi) \over a(\xi)} \cdot {R^+(\xi) \over R^+(\xi^{-1})} =
   - {d(\xi) \over c(\xi)} \cdot {R^-(\xi) \over R^-(\xi^{-1})} \,. \eqno(3.41)
$$
In fact,
$$ {b(\xi) \over a(\xi)} \cdot {R^+(\xi) \over R^+(\xi^{-1})} +
    {d(\xi) \over c(\xi)} \cdot {R^-(\xi) \over R^-(\xi^{-1})}
= - { n^+(\xi) - n^-(\xi) \over n^+(\xi^{-1}) - n^-(\xi) }
                          \cdot {R^+(\xi) \over R^+(\xi^{-1})}
  - { n^+(\xi) - n^-(\xi) \over n^+(\xi) - n^-(\xi^{-1}) }
                          \cdot {R^-(\xi) \over R^-(\xi^{-1})} =
$$
$$ = { - \Bigl( n^+(\xi) - n^-(\xi) \Bigr) \over \xi-\xi^{-1} }
   \biggl\{ { \xi - \xi^{-1} \over N^-(\xi) }
                          \cdot {R^+(\xi) \over R^+(\xi^{-1})}
            - { \xi - \xi^{-1} \over N^+(\xi) }
                          \cdot {R^-(\xi) \over R^-(\xi^{-1})} \biggr\} = 0.
$$
Having verified (3.41), substract (3.39) from (3.40). We obtain
on the unit circle the equality
$$ (1+q(\xi)) (\psi^+(k,\xi)R^+(\xi)-\psi^-(k,\xi)R^-(\xi))
$$
$$ = - {b(\xi)\over a(\xi)} \cdot {R^+(\xi) \over R^+(\xi^{-1})}
      \Bigl( R^+(\xi^{-1}) \psi^+(k,\xi^{-1}) +
             R^-(\xi^{-1}) \psi^-(k,\xi^{-1}) \Bigr) \,.
$$
Thus, for the function $g(k,z)$, defined in (3.1),
we have on the unit circle an equality, which establish connection
between its limit values at the points $\xi$ and $\xi^{-1}$:
$$ \xi^{2(k+1)} (1+q(\xi)) (g^+(k,\xi)-g^-(k,\xi))
   =  - {b(\xi) \over a(\xi)} \cdot {R^+(\xi) \over R^+(\xi^{-1})}
      (g^+(k,\xi^{-1})+g^-(k,\xi^{-1})) \,. \eqno(3.42)
$$

   It remains to reduce the expression
$$ - {b(\xi) \over a(\xi)} \cdot {R^+(\xi) \over R^+(\xi^{-1})}
  = - { n^+(\xi) - n^-(\xi) \over N^-(\xi) }
                          \cdot {R^+(\xi) \over R^+(\xi^{-1})}
$$
$$ = - { n^+(\xi) - n^-(\xi) \over \xi - \xi^{-1} }
                          \cdot {R^-(\xi) R^+(\xi^{-1})}
                          \cdot {R^+(\xi) \over R^+(\xi^{-1})}
   = - { n^+(\xi) - n^-(\xi) \over \xi - \xi^{-1} }
                          \cdot {R^+(\xi) R^-(\xi)} .
$$
Thus, equality (3.42) can be rewritten as (3.31) with the coefficient
$$  \hat r(\xi) =
       { -1 \over 1 + q(\xi) } \cdot
       { n^+(\xi) - n^-(\xi) \over \xi - \xi^{-1} } \cdot
         R^+(\xi) R^-(\xi) \,,
$$
as was to be proved.

\medskip

   4. To complete the proof of the lemma we have to make sure that
$\max_{|\xi|=1} |\hat r(\xi)| < 1$. For this let us write $\hat r(\xi)$
out in the form
$$  \hat r(\xi) =  { - 1 \over 1 + q(\xi) } \cdot
       \Bigl\{ { n^+(\xi) - n^-(\xi) \over \xi - \xi^{-1} } \cdot
         R^+(\xi) R^+(\xi^{-1}) \Bigr\} \,.
$$

   As it was already said, the absolute value of the factor in braces
is less or equal than $1$. (In fact, $|{b(\xi)\over a(\xi)}|<1$ and
$|{R^+(\xi)\over R^+(\xi^{-1})}| =
 |{R_{0123}^+(\xi)\over R_{0123}^+(\xi^{-1})}| \cdot
 |{R_\mu^+(\xi)\over R_\mu^+(\xi^{-1})}| =1 $.)
Hence, it is sufficient to demonstrate that
$$  \max_{|\xi|=1} \, { 1 \over 1 + q(\xi) } < 1 \,,
$$
i.e.~that the positive coefficient $q(\xi)$ is bounded below from zero.
But, from the definition (3.38),
$$ q(\xi) = 2\pi { \sqrt{\rho_R'(\xi+\xi^{-1}) \rho_L'(\xi+\xi^{-1})} \over
                       |n^+(\xi) - n^-(\xi^{-1})| }
   =  { \sqrt{ | {\rm Im} \, m^R(\tau+i0) \
                        {\rm Im} \, {-1 \over m^L(\tau+i0)} | } \over
                       |M(\tau+i0)| } |_{\tau=\xi+\xi^{-1}}
$$
$$ =  { \sqrt{ | m^R(\tau+i0) | \sin \eta^R(\tau) \
                      | {-1 \over m^L(\tau+i0)} | \sin \eta^L(\tau) } \over
                       |M(\tau+i0)| } |_{\tau=\xi+\xi^{-1}} \,.
$$
We remind, that, to within a constant number,
$$ b_{-1} m^R(\lambda) = P(\lambda,\eta^R), \qquad
   {-1 \over b_{-1} m^L(\lambda)} = P(\lambda,\eta^L),
$$
$$ M(\lambda) = b_{-1} m^R(\lambda) - {1 \over b_{-1} m^L(\lambda)}
    = P(\lambda,\eta) \,,
$$
where the functions $\eta^R(\tau)$, $\eta^L(\tau)$ and $\eta(\tau)$ satisfy
the H\"older condition on the interval $(-2,2)$ and have, according to the
condition E), the same jumps at the points $-2$ and $2$. This means that
$$ \eta^R(\tau) = \eta_1^R(\tau) + \delta_1(\tau) \,,
$$
$$ \eta^L(\tau) = \eta_1^L(\tau) + \delta_1(\tau) \,,
$$
$$ \eta(\tau) = \eta_1(\tau) + \delta_1(\tau) \,,
$$
where the functions $\eta_1^R(\tau)$, $\eta_1^L(\tau)$, $\eta_1(\tau)$ satisfy
the H\"older condition in a certain neighborhood of the segment $[-2,2]$,
and $\delta_1(\tau)$ is the function, which has jumps at the points
$-2$ and $2$, and is constant at the other points. Hence,
$$ q(\xi) =  { \sqrt{ | P(\tau+i0,\eta^R) P(\tau+i0,\delta) | \sin \eta^R(\tau) \
                      | P(\tau+i0,\eta^-) P(\tau+i0,\delta) | \sin \eta^-(\tau) } \over
                      | P(\tau+i0,\eta) P(\tau+i0,\delta) | } |_{\tau=\xi+\xi^{-1}}
$$
$$ =  { \sqrt{ | P(\tau+i0,\eta^R) | \sin \eta^R(\tau) \
               | P(\tau+i0,\eta^-) | \sin \eta^-(\tau) } \over
               | P(\tau+i0,\eta) | } |_{\tau=\xi+\xi^{-1}} \,.
$$
The functions $| P(\tau+i0,\eta^R) |$, $| P(\tau+i0,\eta^-) |$ and
$| P(\tau+i0,\eta) |$ are bounded and bounded away from zero, because
$\eta_1^R(\tau)$, $\eta_1^-(\tau)$, and $\eta_1(\tau)$ satisfy the H\"older
condition. Further, $\sin \eta^R(\tau)$ and $\sin \eta^-(\tau)$
are bounded away from zero according to condition B). Thus,
$$ \min_{|\xi|=1} \, q(\xi) > 0 \,,
$$
which complete the proof of the lemma.
\hfill\rule{0.5em}{0.5em}\medskip

\bigskip

\quad{\bf Deduction of the equation}

\medskip

   Let us denote by $\chi_1(\beta)$, $\chi_2^s(\beta)$
and $\chi_2^a(\beta)$ the indicators of the sets $\Phi\cup\Omega_1$,
$(-1,1)\cap\Omega_2^s$ and $\Omega_2^a$ and define the measure
$d\sigma(\beta)$ and the function $u(k,\beta)$ on the union
$$ \Omega_0=\Phi\cup \Omega_1 \cup ((-1,1)\cap\Omega_2^s)\cup\Omega_2^a
$$
of these sets and also on the unit circle by the equalities
$$ d\sigma(\beta) = \chi_1(\beta) {1 \over |\beta|}
   { |\beta-\beta^{-1}| \over |R(\beta)|^2 }
   d\rho_1(\beta)
 + \chi_2^s(\beta) {1 \over |\beta|}
   { |R(\beta)|^2 \over |\beta-\beta^{-1}| }
   d\rho_2(\beta)
 + \chi_2^a(\beta) {1 \over |\beta|}
     { p(\beta) \over 2\pi q(\beta) } d\beta\,,
    \eqno(3.43)
$$

$$ u(k,\beta) = (\chi_1(\beta)+\chi_2^s(\beta))
                 \beta^{-2(k+1)} g(k,\beta^{-1})
$$
$$ - \chi_2^a(\beta) \, i \, p(\beta)^{-1} q(\beta)
     { \beta-\beta^{-1} \over |\beta-\beta^{-1}| }
     (g^+(k,\beta) - g^-(k,\beta))
$$
$$ + \chi_{окр}(\beta) \, {1\over2} \, \hat r(\beta)^{-1}
     (g^+(k,\beta) - g^-(k,\beta)) \,,
   \eqno(3.44)
$$
where the measures $d\rho_1(\tau)$, $d\rho_2(\tau)$ and the functions
$p(\beta)$, $q(\beta)$, $\hat r(\beta)$ are defined by equalities
(3.3), (3.4) and (3.21), (3.22), (3.32)--(3.33), respectively.

   According to the definition, the support of the measure
$d\sigma(\beta)$ is contained in $\Omega_0$.

   In view of (3.2) and lemmas 3.1--3.7
$$ g(k,z) = 1 + { 1 \over 2 \pi i }
   \int\limits_{\Gamma} {g(k,\xi) \over \xi-z } d\xi
$$
$$ = 1 + \int\limits_{\Phi \cup \Omega_1}
   {  -\alpha^{-2(k+1)} g(k,\alpha^{-1})
      { \alpha-\alpha^{-1} \over |\alpha-\alpha^{-1}| } |\alpha| \over
      \alpha - z  } \cdot
   { |\alpha-\alpha^{-1}| \over |R(\alpha)|^2 } \cdot
   {1 \over |\alpha|} d\rho_1(\alpha)
$$
$$ + \int\limits_{\Omega_2^s \cap (-1,1)}
   {  -\alpha^{-2(k+1)} g(k,\alpha^{-1})
      { \alpha-\alpha^{-1} \over |\alpha-\alpha^{-1}| } |\alpha| \over
      \alpha - z  } \cdot { |R(\alpha)|^2 \over |\alpha-\alpha^{-1}| }
      {1 \over |\alpha|} d\rho_2(\alpha)
$$
$$ + \int\limits_{\Omega_2^a}
     { (g^+(k,\alpha) - g^-(k,\alpha)) q(\alpha) |\alpha|  \over
        i \, p(\alpha) (\alpha-z)  }  \cdot
     {  { \alpha-\alpha^{-1} \over |\alpha-\alpha^{-1}| } \over
        { \alpha-\alpha^{-1} \over |\alpha-\alpha^{-1}| }          }  \cdot
     {1 \over |\alpha|} \cdot { p(\alpha) \over 2 \pi q(\alpha) } d\alpha
$$
$$ + {1 \over \pi } \int\limits_{|\xi|=1}
     { {1\over2} \hat r(\xi)^{-1}
       (g^+(k,\xi) - g^-(k,\xi)) \over
        \xi - z  }  \cdot
      {1\over i} \hat r(\xi) d\xi
$$

$$ = 1 + \int\limits_{\Omega_0}
     {  u(k,\alpha) { \alpha-\alpha^{-1} \over |\alpha-\alpha^{-1}| } \over
      \alpha - z  } |\alpha|  d \sigma(\alpha)
   + {1 \over \pi } \int\limits_{-\pi}^\pi
     { u(k,e^{i\theta}) \over 1 - z e^{ - i \theta } } \hat r(e^{i\theta})
     d\theta
$$

$$ = 1
   + \int\limits_{\Omega_0}
     { u(k,\alpha) s(\alpha) \over 1 - z \alpha^{-1} }
     d \sigma (\alpha)
   + {1\over\pi} \int\limits_{-\pi}^\pi
     { \hat r(e^{i\theta}) u (k,e^{i\theta})
       \over 1 - z e^{-i\theta} } d \theta
$$
where, as it was already defined,
$   s (\alpha) = { \alpha \over |\alpha| } \cdot
                    { \alpha-\alpha^{-1} \over |\alpha-\alpha^{-1}| }
     = \cases { -1, \quad \alpha \in (-1,1),  \cr
                 1, \quad \alpha \in (-\infty,-1) \cup (1,+\infty).  \cr  }
$ \quad
Hence, when
$\beta\in\Omega_0$ or $|\beta| = 1$, we have
$$ { g^+(k,\beta^{-1}) + g^-(k,\beta^{-1}) \over 2 } = 1 +
    \int\limits_{\Omega_0 \backslash \Omega_2^a}
     {  u(k,\alpha) { \alpha-\alpha^{-1} \over |\alpha-\alpha^{-1}| } |\alpha| \over
      \alpha - \beta^{-1}  } d \sigma(\alpha)
    + {\rm v.p.} \int\limits_{\Omega_2^a}
     {  u(k,\alpha) { \alpha-\alpha^{-1} \over |\alpha-\alpha^{-1}| }
                                                     |\alpha|      \over
      \alpha - \beta^{-1}  } d \sigma(\alpha)
$$
$$  + {1 \over \pi } {\rm v.p.} \int\limits_{-\pi}^\pi
     { \hat r(e^{i\theta}) u(k,e^{i\theta}) \over 1 - \beta^{-1} e^{-i\theta} }
     d\theta \,,   \eqno(3.45)
$$
where ${\rm v.p.} {\displaystyle \int} $ denotes the principal value
of the integral. Let us denote, for the sake of brevety, the sum of
the integrals in the right-hand side of (3.45) by
${\rm v.p.} {\displaystyle \int}$ (although not all of
them are singular for $ \beta \not\in \Omega_0\backslash\Omega_2^a$).
If $\beta\in \Omega_0\backslash\Omega_2^a$, then according to (3.44),
$$ { g^+(k,\beta^{-1}) + g^-(k,\beta^{-1}) \over 2 } =
   - (\chi_1(\beta) + \chi_2^s(\beta)) \beta^{2(k+1)} u(k,\beta),
$$
and if $\beta\in \Omega_2^a$, then, in view of (3.44) and lemma 3.5,
$$ { g^+(k,\beta^{-1}) + g^-(k,\beta^{-1}) \over 2 } =
      { \beta-\beta^{-1} \over |\beta-\beta^{-1}| }
     { i \, q(\beta) \over p(\beta) } \beta^{2(k+1)}
     (g^+(k,\beta) - g^-(k,\beta))
$$
$$   + { \beta \over |\beta| }
       { \beta-\beta^{-1} \over |\beta-\beta^{-1}| }
     { i \, m(\beta) \over 2 \, p(\beta) }
     (g^+(k,\beta^{-1}) - g^-(k,\beta^{-1}))
   = - \beta^{2(k+1)} u(k,\beta)
     + { \beta \over |\beta| } { m(\beta) \over 2 q(\beta^{-1}) }
        u(k,\beta^{-1})
$$
(we remind that $p(\beta)=p(\beta^{-1})$).
Hence, for all $\beta\in\Omega_0$
$$ { g^+(k,\beta^{-1}) + g^-(k,\beta^{-1}) \over 2 } =
   - \chi_{\Omega_0}(\beta) \beta^{2(k+1)} u(k,\beta)
   + \chi_{\Omega_2^a}(\beta) { \beta \over |\beta| }
     { m(\beta) \over 2 q(\beta^{-1}) } u(k,\beta^{-1}) \,.
$$
Finally, if $|\beta|=1$, then
$$ { g^+(k,\beta^{-1}) + g^-(k,\beta^{-1}) \over 2 } =
   - {1\over2} \hat r(\beta)^{-1} \beta^{2(k+1)}
     (g^+(k,\beta) - g^-(k,\beta))  = - \beta^{2(k+1)} u(k,\beta) \,.
$$
Thus, taking into account these equalities and (3.45), we conclude that
$u(k,\beta)$ satisfy the equation
$$ \beta^{2(k+1)} u(k,\beta) -
   \chi_{\Omega_2^a}(\beta) { \beta \over |\beta| }
   { m(\beta) \over 2 q(\beta^{-1}) } u(k,\beta^{-1})
$$
$$ + \, {\rm v.p.}\int\limits_{\Omega_0}
     { u(k,\alpha) \over 1 - \beta^{-1} \alpha^{-1}  }
     s(\alpha) d \sigma (\alpha)
   + {1\over\pi} \, {\rm v.p.} \int\limits_{-\pi}^\pi
     { \hat r(e^{i\theta}) u(k,e^{i\theta})
       \over 1 - \beta^{-1} e^{-i\theta}  } d \theta
   +  1 =0 \,,   \eqno(3.46)
$$
where $ \chi_{\Omega_0}(\beta) + \chi_{окр.}(\beta) = 1 $
on the set on which the equation is considered.

In this section  we proved

\bigskip

{\bf Theorem 3. } \ { \it The function $g(k,z)$ can be represented in
the form
$$ g(k,z) = 1
   + \int\limits_{\Omega_0}
     { u(k,\alpha) s(\alpha) \over 1 - z \alpha^{-1} }
     d \sigma (\alpha)
   + {1\over\pi} \int\limits_{-\pi}^\pi
     { \hat r(e^{i\theta}) u(k,e^{i\theta})
       \over 1 - z e^{-i\theta} } d \theta.  \eqno(3.47)
$$
for $ z \not\in \Omega_0 \cup {\bf T}$, where the function $u(k,\alpha)$
is the solution of integral equation $(3.46)$,
and the function $s(\alpha)$, the coefficients $q(\beta)$, $m(\beta)$,
$\hat r(e^{i\theta})$ and the measure $d\sigma(\alpha)$ are defined by
equalities (2.31), (3.22), (3.23), (3.32)--(3.33) and (3.44), respectively.}

\medskip

{\bf Remark.} \ Equation (3.46), as it is seen from its form, is considered
on the unit circle $\bf T$ and on the support of the measure $d\sigma(\alpha)$,
which is a subset of the set $\Omega_0$ (let us denote this support
by $\Omega_0'$). In fact, the set $\Omega_0' \cup {\bf T}$
{\it exactly} corresponds (in the parametrization $z+z^{-1}$) to the spectrum
of the matrix $J$. We outline a short explanation of this.
(For the sake of convenience we take our considerations in $z$-plane,
so in this remark by the word "spectrum" we mean the image in the plane
$z$ of the honest spectrum in the plane $\lambda=z+z^{-1}$.)

   Since $J$ is a finite-dimensional perturbation of the direct sum
of the operators $J^R$ and $J^L$, the continuous spectrum $J$
is the same as the union of the continuous  spectrum $J^R$ and $J^L$,
i.e.~the nonisolated points of $\Omega_1\cup\Omega_2^a\cup{\bf T}$.
Besides, the spectrum of $J$ includes its isolated eigenvalues. According
to (2.1) and (2.2), these are zeroes of the function $N(z)$ (a part of the
set $\Phi$) and the common poles of the functions $n(z)$ and $n(z^{-1})$
(i.e.~the set $\Omega_2^s$). Hence, the spectrum of $J$ belongs to
$\Phi\cup \Omega_1 \cup \Omega_2^s \cup\Omega_2^a \cup {\bf T} =
 \Omega_0 \cup {\bf T} $.

   But what we need to prove is to prove that the spectrum of $J$
is exactly $\Omega_0' \cup {\bf T}$, that is to prove that

(i) the parts of $\Omega_0$ that do not belong to the support of the measure
$d\sigma(\alpha)$, do not belong to the spectrum.

(ii) the whole of the set $\Omega'_0 \cup {\bf T} $ (i.e.~the circle ${\bf T}$ and
the support $\Omega'_0$ of the measure $d\sigma(\alpha)$)
belongs to the spectrum.

   Let us begin with the inverse inclusion (ii). Since all the points
$\Omega_2^a\cup{\bf T}$ are nonisolated ($\Omega_2^a$ entirely belongs
to the support of the measure $d\sigma(\alpha)$), and $\rho_R(\tau)$, $\rho_L(\tau)$
are absolutely continuous on $\Omega_2^a\cup T$, all this set
belongs to the continuous spectrum of $J^R$ and $J^L$, from where it follows
that $\Omega_2^a\cup{\bf T}$ {\it completely} belongs to the spectrum
of $J$ (more precisely, to the absolutely continuous one).
Further, according to (3.43) and (3.7), the support of the measure
$d\sigma(\alpha)$ contains the set $(-1,1)\cap\Omega_2^s$. But, according
to (2.2), this set completely belongs to the spectrum, because it consists
of common poles of the functions $n(z)$ and $n(z^{-1})$ and, hence, these are
poles of the function ${n(z)n(z^{-1})\over b_{-1} N(z)}$. To finish the proof
of (ii) we have to show that the spectrum of $J$ include
the whole intersection of the support of $d\sigma(\alpha)$ with the set
$\Phi\cup\Omega_1$. But it is evident from the definition of $d\sigma(\alpha)$
in terms of the measure $d\rho_1(\alpha)$ on this set (formula (3.43))
and from the definition of the measure $d\rho_1(\alpha)$ (formula (3.4))
by the element of the resolvent $N(z)^{-1}$ (formula (2.1));
we remind that all the points of increase of $\rho_1(t)$ are singularities
of $N(z)^{-1}$, and all the singularities of the resolvent belong to the
spectrum.

   The last speculation also prove (i). In fact, $\Omega_2^a$
and $(-1,1)\cap\Omega_2^s$ completely belong to the support of the measure
$d\sigma(\alpha)$. That is, in the set $\Omega_0$ only
a certain sublet $(\Phi\cup\Omega_1)\backslash \Omega'_0$ of the set
$\Phi\cup\Omega_1$ may not belong to the support of the measure
$d\sigma(\alpha)$; more precisely, some isolated points of this set.
(In other terms,
$\Omega_0\backslash \Omega'_0 = (\Phi\cup\Omega_1)\backslash \Omega'_0$.)
But, according to (3.43), if
$\sigma\Bigl((\Phi\cup\Omega_1)\backslash \Omega'_0\Bigr)=0$, then also
$\rho_1\Bigl((\Phi\cup\Omega_1)\backslash \Omega'_0\Bigr)=0$,
from where, according to (3.4), it follows that the points of the set
$(\Phi\cup\Omega_1)\backslash \Omega'_0$ are points of holomorphy
of $N(z)^{-1}$. Hence, according to (2.1), $R(-1,-1,z+z^{-1})$
is regular on this set. The function $R(0,0,z+z^{-1})$ is also regular
at these points. (Since the set $\Omega_1$ is asymetric, either $n(z)$,
or $n(z^{-1})$ is regular at each point, so, in view of (2.2), every pole
of $n(z)$ of $n(z^{-1})$ is reduced with the same pole of $N(z)$).
Thus, the resolvent is regular on the set
$(\Phi\cup\Omega_1)\backslash \Omega'_0$. Therefore, this set does not belong
to the spectrum.

\bigskip

{\bf Conclusion. } \  {\it Thus for the reconstruction of the infinite
Jacobi matrix $J$ by its spectral data $n(z)$ one have to find the functions
$\hat r(e^{i\theta})$, $m(\beta)$, $q(\beta)$ and the measure
$d\sigma(\alpha)$, corresponding to the Weyl function $n(z)$, then to solve
the equation $(3.46)$ for all $k \in {\bf Z}$ and define from $u(k,\beta)$
the function $g(k,z)$ and the entries $a_k$, $b_k$ of the Jacobi matrix
by (3.47) and (0.29) ($b_k$ are reconstructed up to their sign).
Here the points of increase of the function $\sigma$, together with
the points of the unit circle are exactly the set of such $z$ for which
$z+z^{-1}$ belongs to the spectrum of matrix $J$, and also both points
$z\in {\bf R}$ and $z^{-1}\in {\bf R}$ are the points of increase of the
function $\sigma$ if and only if $z+z^{-1}$ belongs to the absolutely
continuous spectrum of multiplicity $2$ of the matrix $J$.

 }

\bigskip
\bigskip

\centerline {\bf 4. Solvability of the fundamental equation}

\bigskip

In this section we will prove the solvability of the equation of more general
form than (3.46), namely
$$ e^{\beta t} \beta^{2(k+1)} u(k,\beta) -
   e^{\beta^{-1} t} \chi_{\Omega_2^a}(\beta) { \beta \over |\beta| }
   {m(\beta) \over 2 q(\beta^{-1}) } u(k,\beta^{-1})
$$
$$ + e^{\beta^{-1} t}\, {\rm v.p.}\int\limits_{\Omega_0}
     { u(k,\alpha) \over 1 - \beta^{-1} \alpha^{-1}  }
     s(\alpha) d \sigma (\alpha)
   + e^{\beta^{-1} t}{1\over\pi} \, {\rm v.p.} \int\limits_{-\pi}^\pi
     { \hat r(e^{i\theta}) u(k,e^{i\theta})
       \over 1 - \beta^{-1} e^{-i\theta}  } d \theta
    = -1 \,,   \eqno(4.1)
$$
where $t\ge 0$ is a real parameter. It is evident that, when $t=0$,
(4.1) coincide with the equation of the inverse problem (3.46).
After dividing (4.1) by $e^{\beta^{-1} t}$ and denoting
$$ v(k,\beta) = e^{(\beta - \beta^{-1}) t} \beta^{2(k+1)} u(k,\beta) \,.
$$
we rewrite the equation for $v(k,\beta)$ in the form
$$ v(k,\beta) -  \chi_{\Omega_2^a}(\beta) { \beta \over |\beta| }
   {m(\beta) \over 2 q(\beta^{-1}) }
   { v(k,\beta^{-1}) \over e^{-(\beta - \beta^{-1}) t} \beta^{-2(k+1)} }
  + \, {\rm v.p.}\int\limits_{\Omega_0}
     { v(k,\alpha) \over 1 - \beta^{-1} \alpha^{-1}  } \cdot
     { s(\alpha) \over e^{(\alpha - \alpha^{-1}) t} \alpha^{2(k+1)} }
     d \sigma (\alpha)
$$
$$ + {1\over\pi} \, {\rm v.p.} \int\limits_{-\pi}^\pi
     { \hat r(e^{i\theta}) v(k,e^{i\theta})        1
       \over 1 - \beta^{-1} e^{-i\theta}  } \cdot
     { 1 \over exp((e^{i\theta} - e^{-i\theta}) t) e^{2i\theta(k+1)} }
     d \theta
   = - { 1 \over e^{\beta^{-1} t} } \,,   \eqno(4.2)
$$
Let us introduce for $|\beta|=1$ the function
$$ \tilde r(\beta) = { \hat r(\beta) \over e^{(\beta - \beta^{-1}) t} \beta^{2(k+1)} }
$$
and the measure on the real line
$$ d \tilde\sigma(\alpha) = { d \sigma(\alpha) \over
         e^{(\alpha - \alpha^{-1}) t} \alpha^{2(k+1)} } \,. \eqno(4.3)
$$
Since on the circle we have $|e^{(\beta - \beta^{-1}) t} \beta^{2(k+1)}|=1$
\ ($\beta - \beta^{-1}$ is a pure imaginary number), then like it was in
lemma 3.7,
$$ \max_{|\beta|=1} |\tilde r(\beta)| \leq 1 \,.
$$

\medskip

In new notations the equation (4.2) is replaced by
$$ v(k,\beta) -  \chi_{\Omega_2^a}(\beta) { \beta \over |\beta| }
   {m(\beta) \over 2 q(\beta^{-1}) } e^{(\beta - \beta^{-1}) t} \beta^{2(k+1)}
    v(k,\beta^{-1})
$$
$$ + \, {\rm v.p.}\int\limits_{\Omega_0}
     { v(k,\alpha) s(\alpha) \over 1 - \beta^{-1} \alpha^{-1}  }
     d \tilde\sigma (\alpha)
   + {1\over\pi} \, {\rm v.p.} \int\limits_{-\pi}^\pi
     { \tilde r(e^{i\theta}) v(k,e^{i\theta})
       \over 1 - \beta^{-1} e^{-i\theta}  } d \theta
   = - { 1 \over e^{\beta^{-1} t} } \,.   \eqno(4.4)
$$

We will seek the solution of this equation in the class
$L^2({\bf R}_{\tilde\sigma} \cup {\bf T}_{\hat m\over \pi})$ \
(${\bf T}$ is the unit circle and ${\hat m\over \pi}$ is the Lebesque
measure on it, divided by $\pi$).

   So, we will prove the invertibility of the operator $L$:
$$ (L v) (\beta) = v(\beta) -
   \chi_{\Omega_2^a}(\beta) { \beta \over |\beta| }
   {m(\beta) \over 2 q(\beta^{-1}) } e^{(\beta - \beta^{-1}) t} \beta^{2(k+1)}
    v(\beta^{-1})
$$
$$ + \, {\rm v.p.}\int\limits_{\Omega_0}
     { v(\alpha) s(\alpha) \over 1 - \beta^{-1} \alpha^{-1}  }
     d \tilde\sigma (\alpha)
   + {1\over\pi} \, {\rm v.p.} \int\limits_{-\pi}^\pi
     { \tilde r(e^{i\theta}) v (e^{i\theta})
       \over 1 - \beta^{-1} e^{-i\theta}  } d \theta,
$$
which is defined in the introduced space (it is evident that the right-hand
side $ - e^{-\beta^{-1} t}$ of the equation belongs to this space).
For this we will use the following lemma

\medskip

{\bf Lemma 4.1.} /see~[20], p.~107/
{\it Let $A \in B({\bf H})$ ($B({\bf H})$ denotes the algebra of the bounded
operators in the Hilbert space ${\bf H}$) and let there exist such positive
number $d$ that ${\rm Re}\, (A \vec f, \vec f) \ge d \|\vec f\|^2$ for all
$\vec f \in {\bf H}$.  Then the operator $A$ is invertible and
$\|A^{-1}\| \leq d^{-1}$. }

\medskip

It is seen from the last lemma that we need to calculate the real part of the
expression

$$ (L v,v) = \int\limits_{\Omega_0} d\tilde\sigma(\beta) \overline{v(\beta)}
\Biggl\{ v(\beta) - \chi_{\Omega_2^a}(\beta) { \beta \over |\beta| }
   {m(\beta) \over 2 q(\beta^{-1}) } e^{(\beta - \beta^{-1}) t} \beta^{2(k+1)}
    v(\beta^{-1})
$$
$$ + \, {\rm v.p.}\int\limits_{\Omega_0}
     { v(\alpha) s(\alpha) \over 1 - \beta^{-1} \alpha^{-1}  }
     d \tilde\sigma (\alpha)
   + {1\over\pi} \, {\rm v.p.} \int\limits_{-\pi}^\pi
     { \tilde r(e^{i\theta}) v (e^{i\theta})
       \over 1 - \beta^{-1} e^{-i\theta}  } d \theta \Biggr\}
$$
$$ + {1\over\pi} \int\limits_{-\pi}^\pi d\phi \ \overline{v(e^{i\phi})}
     \Biggl\{ v(e^{i\phi})
   + \, {\rm v.p.}\int\limits_{\Omega_0}
     { v(\alpha) s(\alpha) \over 1 - \alpha^{-1} e^{-i\phi} }
     d \tilde\sigma (\alpha)
   + {1\over\pi} \, {\rm v.p.} \int\limits_{-\pi}^\pi
     { \tilde r(e^{i\theta}) v (e^{i\theta})
       \over 1 - e^{-i\phi} e^{-i\theta}  } d \theta \Biggr\}
$$
$$ = (I) + (II) +(III) + (IV) + (VI) + (VII) $$
(we decomposed $(Lv,v)$ in the sum of terms $(I)$ -- $(VII)$, which
are obtained after removing the parentheses). Let us transform each one
of this terms. First, we remark that the second term can be omitted from
our consideration because it is a pure imaginary number whereas we only
need to estimate ${\rm Re}\, (L v, v)$. In fact, taking into account that,
according to (4.3) and (3.33), on $\Omega_2^a$ the measure
$d\tilde\sigma(\beta)
  =  { d\sigma(\beta) \over e^{(\beta - \beta^{-1}) t} \beta^{2(k+1)}  }
  =  { 1 \over e^{(\beta - \beta^{-1}) t} \beta^{2(k+1)}  } \cdot
     {p(\beta) \over 2\pi |\beta| q(\beta) } d\beta$,
and that, according to (3.21) and (3.23), we have
$p(\alpha)=p(\alpha^{-1})$, $m(\alpha)=-m(\alpha^{-1})$,
we obtain that
$$ (II) = - \int\limits_{\Omega_2^a} d\tilde\sigma(\beta) \overline{v(\beta)}
    \cdot { \beta \over |\beta| } \cdot
   {m(\beta) \over 2 q(\beta^{-1}) } e^{(\beta - \beta^{-1}) t} \beta^{2(k+1)}
    v(\beta^{-1})
$$
$$ = - \int\limits_{\Omega_2^a}
    d\beta {p(\beta) \over 2\pi} \cdot
 { \overline{v(\beta)} \over
   e^{(\beta - \beta^{-1}) t} \beta^{2(k+1)} |\beta| q(\beta) } \cdot
   { \beta \over |\beta| } \cdot
   e^{(\beta - \beta^{-1}) t} \beta^{2(k+1)}
   {m(\beta) \over 2 q(\beta^{-1}) } v(\beta^{-1})
$$
$$ = - \int\limits_{\Omega_2^a}
    d\beta {p(\beta) \over 2\pi}
   { \overline{v(\beta)} v(\beta^{-1})  \over |\beta| q(\beta) } \cdot
   { \beta \over |\beta| } \cdot
   {m(\beta) \over 2 q(\beta^{-1}) }
$$
$$ = \Bigl[ \alpha = \beta^{-1} ] = - \int\limits_{\Omega_0}
   {d\alpha \over \alpha^2} {p(\alpha) \over 2\pi}
   \overline {
               { \overline{v(\alpha)} v(\alpha^{-1})
               \over |\alpha|^{-1} q(\alpha^{-1}) }
              }                                       \cdot
   { \alpha\over |\alpha| } \cdot
   {m(\alpha) \over 2 q(\alpha) }
$$
$$ = - \int\limits_{\Omega_0}
    d\alpha {p(\alpha) \over 2\pi}
   \overline {
               { \overline{v(\alpha)} v(\alpha^{-1}) }
                              \over |\alpha| q(\alpha)
             }                                          \cdot
   { \alpha \over |\alpha| } \cdot
   {m(\alpha) \over 2 q(\alpha^{-1}) } \ = \ - \overline{(II)}.
$$
Thus, ${\rm Re}\, (II) =0$, and later on we will not take $(II)$
into account.

   Let us now extract the real part from the third term. We denote by
$A$ the integral operator
$$  (Av)(\beta)= \, {\rm v.p.}\int\limits_{\Omega_0}
     { v(\alpha) s(\alpha) \over 1 - \beta^{-1} \alpha^{-1}  }
     d \tilde\sigma (\alpha)  \,,
$$
defined in $L^2_{\tilde\sigma}({\bf R})$. For this operator
$$ A= { A + A^* \over 2} + { A - A^* \over 2} \,.
$$
The second term ${ A - A^* \over 2}$ is an asymetric operator in
$L^2_{\tilde\sigma}({\bf R})$. Hence,
$({ A - A^* \over 2} v,v)_{\tilde\sigma}$ is a pure imaginary number.
Further, ${ A + A^* \over 2}$ is a symmetric operator of the form
$$  ({ A + A^* \over 2} v)(\beta)= \,
    {\rm v.p.}\int\limits_{\Omega_0}
     { S(\alpha,\beta) \over 1 - \beta^{-1} \alpha^{-1}  }
     v(\alpha) d \tilde\sigma (\alpha)  \,,
$$
where $ S(\alpha,\beta) = { s(\alpha) + s(\beta) \over 2 } $.
Hence,
$$ {\rm Re} \, (III) = {\rm Re} \, ( A v,v)_{\tilde\sigma}
   = ({ A + A^* \over 2} v,v)_{\tilde\sigma}
   = \int\limits_{\Omega_0}
     d \tilde\sigma (\beta) \overline{ v(\beta) }
     \int\limits_{\Omega_0}
     { S(\alpha,\beta) \over 1 - \beta^{-1} \alpha^{-1}  }
     v(\alpha) d \tilde\sigma (\alpha)  \,.
$$
Let us remark that the internal integral is no more singular because
its singularities were the points for which $\alpha=\beta^{-1}$,
and now for all these points $ S(\alpha,\beta) = 0 $.

\medskip

Thus,
$$ {\rm Re} \, (Lv,v) = \int\limits_{\Omega_0}
     d \tilde\sigma (\alpha) \overline{ v(\beta) } v(\beta)
$$
$$ + \int\limits_{\Omega_0}
     d \tilde\sigma (\beta) \overline{ v(\beta) }
     \int\limits_{\Omega_0}
     { S(\alpha,\beta) \over 1 - \beta^{-1} \alpha^{-1}  }
     v(\alpha) d \tilde\sigma (\alpha)
   + {\rm Re} \, \int\limits_{\Omega_0} d\tilde\sigma(\beta) \overline{v(\beta)}
     \cdot {1\over\pi} \, \int\limits_{-\pi}^\pi
     { \tilde r(e^{i\theta}) v (e^{i\theta})
       \over 1 - \beta^{-1} e^{-i\theta}  } d \theta
$$
$$ + {1\over\pi} \int\limits_{-\pi}^\pi d\phi \
     \overline{v(e^{i\phi})} v(e^{i\phi})
   + {\rm Re} \, {1\over\pi} \int\limits_{-\pi}^\pi d\phi \
     \overline{v(e^{i\phi})}
     \, \int\limits_{\Omega_0}
     { v(\alpha) s(\alpha) \over 1 - \alpha^{-1} e^{-i\phi} }
     d \tilde\sigma (\alpha)
$$
$$ + {\rm Re} \, {1\over\pi} \int\limits_{-\pi}^\pi d\phi \
     \overline{v(e^{i\phi})} \cdot
     {1\over\pi} \, {\rm v.p.} \int\limits_{-\pi}^\pi
     { \tilde r(e^{i\theta}) v (e^{i\theta})
       \over 1 - e^{-i\phi} e^{-i\theta}  } d \theta \,.
$$

   The introduced $L^2$-space is expanded
into the direct sum of two spaces: $L^2$ on the real line and $L^2$
on the circle:
$$ L^2({\bf R}_{\tilde\sigma} \cup {\bf T}_{\hat m\over \pi})
   = L^2({\bf R}_{\tilde\sigma}) \dot + L^2({\bf T}_{\hat m\over \pi}) \,.
$$
Let us expand the function $v$ into the sum of its projections in
$\bf R$ and $\bf T$:
$$ v = v_r + c_t \,.
$$
Then
$$ {\rm Re} \, (Lv,v) = \|v_r\|^2
   + \int\limits_{\Omega_0\cap(-1,1)}
     d \tilde\sigma (\beta) \overline{ v(\beta) }
     \int\limits_{\Omega_0\cap(-1,1)}
     { -1 \over 1 - \beta^{-1} \alpha^{-1}  }
     v(\alpha) d \tilde\sigma (\alpha)
$$
$$   + \int\limits_{\Omega_0\cap{\bf R}\backslash (-1,1)}
     d \tilde\sigma (\beta) \overline{ v(\beta) }
     \int\limits_{\Omega_0\cap{\bf R}\backslash (-1,1)}
     { 1 \over 1 - \beta^{-1} \alpha^{-1}  }
     v(\alpha) d \tilde\sigma (\alpha)
$$
$$ + {\rm Re} \,
      \Bigl\{ \int\limits_{\Omega_0\cap(-1,1)} +
              \int\limits_{\Omega_0\cap{\bf R}\backslash (-1,1)} \Bigr\}
      d\tilde\sigma(\beta) \overline{v(\beta)} \cdot
     {1\over\pi} \, \int\limits_{-\pi}^\pi
     { \tilde r(e^{i\theta}) v (e^{i\theta})
       \over 1 - \beta^{-1} e^{-i\theta}  } d \theta
$$
$$ + \|v_t\|^2
   + {\rm Re} \, {1\over\pi} \int\limits_{-\pi}^\pi d\phi \
     \overline{v(e^{i\phi})}
     \, \Bigl\{ - \int\limits_{\Omega_0\cap(-1,1)} +
                  \int\limits_{\Omega_0\cap{\bf R}\backslash (-1,1)}
        \Bigr\}
     { v(\alpha) \over 1 - \alpha^{-1} e^{-i\phi} }
     d \tilde\sigma (\alpha)
$$
$$ + {\rm Re} \, {1\over\pi} \int\limits_{-\pi}^\pi d\phi \
     \overline{v(e^{i\phi})} \cdot
     {1\over\pi} \, {\rm v.p.} \int\limits_{-\pi}^\pi
     { \tilde r(e^{i\theta}) v (e^{i\theta})
       \over 1 - e^{-i\phi} e^{-i\theta}  } d \theta \,. \eqno(4.5)
$$
In the first six of the seven integrals obtained
we decompose the kernels of the form
${ 1 \over 1 - \beta^{-1} \alpha^{-1}  }$,
${ 1 \over 1 - \beta^{-1} e^{-i\theta} }$,
${ 1 \over 1 - \alpha^{-1} e^{-i\phi} }$
into sums of geometric progression:

(1) $\alpha \in (-1,1), \ \beta \in (-1,1)$:

$${ 1 \over 1 - \beta^{-1} \alpha^{-1}  }
   = - \alpha\beta \cdot { 1 \over 1 - \alpha \beta }
   = - \alpha\beta \sum_{n=0}^{+\infty} \alpha^n \beta^n
   = - \sum_{n=1}^{+\infty} \alpha^n \beta^n \,;
$$

(2) $\alpha \in {\bf R}\backslash (-1,1), \
      \beta \in {\bf R}\backslash (-1,1)$:

$${ 1 \over 1 - \beta^{-1} \alpha^{-1}  }
   = \sum_{n=0}^{+\infty} {(\alpha^{-1} \beta^{-1})}^n
   = \sum_{n=0}^{+\infty} {1 \over \alpha^n \beta^n } \,;
$$

(3) $\theta \in (-\pi,\pi), \ \beta \in (-1,1)$:

$${ 1 \over 1 - \beta^{-1} e^{-i\theta}  }
   = - \beta e^{i\theta} \cdot { 1 \over 1 - \beta e^{i\theta} }
   = - \beta e^{i\theta} \sum_{n=0}^{+\infty} \beta^n e^{in\theta}
   = - \sum_{n=1}^{+\infty} \beta^n e^{in\theta} \,;
$$

(4) $\theta \in (-\pi,\pi), \ \beta \in {\bf R}\backslash (-1,1)$:

$${ 1 \over 1 - \beta^{-1} e^{-i\theta}  }
   = \sum_{n=0}^{+\infty} {(\beta^{-1} e^{-i\theta})}^n
   = \sum_{n=0}^{+\infty} { 1 \over \beta^n e^{in\theta} } \,;
$$

(5) $\phi \in (-\pi,\pi), \ \alpha \in (-1,1)$  (аналогично (3)):

$${ 1 \over 1 - \alpha^{-1} e^{-i\phi}  }
   = - \sum_{n=1}^{+\infty} \alpha^n e^{in\phi} \,;
$$

(6) $\phi \in (-\pi,\pi), \ \alpha \in {\bf R}\backslash (-1,1)$
    (analogously to (4)):

$${ 1 \over 1 - \alpha^{-1} e^{-i\phi}  }
   = \sum_{n=0}^{+\infty} { 1 \over \alpha^n e^{in\phi} } \,.
$$

Let us now calculate the seven integral in (4.5). First,
$$   {1\over\pi} \, {\rm v.p.} \int\limits_{-\pi}^\pi
     { \tilde r(e^{i\theta}) v (e^{i\theta})
       \over 1 - e^{-i\phi} e^{-i\theta}  } d \theta
  =  {1\over2\pi} \lim_{\varepsilon \downarrow 0}
      \Bigl\{ \int\limits_{-\pi}^\pi
     { \tilde r(e^{i\theta}) v (e^{i\theta})
       \over 1 - e^{-i\phi} e^{-i\theta} (1-\varepsilon) } d \theta
      + \int\limits_{-\pi}^\pi
     { \tilde r(e^{i\theta}) v (e^{i\theta})
       \over 1 - e^{-i\phi} e^{-i\theta} (1+\varepsilon) } d \theta \Bigr\} \,.
$$
so we once more decompose the kernels of the two last integrals
in the sum of geometric progression:

(1) $ |e^{-i\phi} e^{-i\theta} (1-\varepsilon)|<1 $

$$ { 1 \over 1 - e^{-i\phi} e^{-i\theta} (1-\varepsilon) }
   = \sum_{n=0}^{+\infty} e^{-in\phi} e^{-in\theta} (1-\varepsilon)^n \,;
$$

(2) $ |e^{-i\phi} e^{-i\theta} (1+\varepsilon)|>1 $

$$ { 1 \over 1 - e^{-i\phi} e^{-i\theta} (1+\varepsilon) }
   = - e^{i\phi} e^{i\theta} (1+\varepsilon) \cdot
   { 1 \over 1 - e^{i\phi} e^{i\theta} (1+\varepsilon)^{-1} }
   = - \sum_{n=1}^{+\infty} e^{in\phi} e^{in\theta} {(1+\varepsilon)}^{-n} \,.
$$

Hence,
$$   {1\over\pi} \, {\rm v.p.} \int\limits_{-\pi}^\pi
     { \tilde r(e^{i\theta}) v (e^{i\theta})
       \over 1 - e^{-i\phi} e^{-i\theta}  } d \theta
$$
$$ =  {1\over2\pi} \lim_{\varepsilon \downarrow 0}
      \Bigl\{ \sum_{n=0}^{+\infty} e^{-in\phi} (1-\varepsilon)^n
       \int\limits_{-\pi}^\pi \tilde r(e^{i\theta}) v (e^{i\theta})
           e^{-in\theta} d \theta
      - \sum_{n=1}^{+\infty} e^{in\phi} (1+\varepsilon)^{-n}
       \int\limits_{-\pi}^\pi \tilde r(e^{i\theta}) v (e^{i\theta})
           e^{in\theta} d \theta \,.
$$
Thus, the last integral in expression (4.5) is equal to
$$  {1\over\pi} \int\limits_{-\pi}^\pi d\phi \
     \overline{v(e^{i\phi})} \cdot
     {1\over\pi} \, {\rm v.p.} \int\limits_{-\pi}^\pi
     { \tilde r(e^{i\theta}) v (e^{i\theta})
       \over 1 - e^{-i\phi} e^{-i\theta}  } d \theta
$$
$$  = {1\over2} \sum_{n=0}^{+\infty} {1\over\pi} \int\limits_{-\pi}^\pi
       e^{-in\phi} \overline{v(e^{i\phi})} d\phi \cdot
       {1\over\pi} \int\limits_{-\pi}^\pi
        \tilde r(e^{i\theta}) v (e^{i\theta}) e^{-in\theta} d \theta
$$
$$  - {1\over2} \sum_{n=1}^{+\infty} {1\over\pi} \int\limits_{-\pi}^\pi
       e^{in\phi} \overline{v(e^{i\phi})} d\phi \cdot
       {1\over\pi} \int\limits_{-\pi}^\pi
        \tilde r(e^{i\theta}) v (e^{i\theta}) e^{in\theta} d \theta \,.
$$

   After substituting the obtained expression into sum (4.5)
and taking into account that the support of the measure
$d \tilde\sigma (\alpha)$ lies in $\Omega_0$, we obtain that
$$ {\rm Re} \, (Lv,v) = \|v_r\|^2
   + \int\limits_{(-1,1)}
     d \tilde\sigma (\beta) \overline{ v(\beta) }
     \int\limits_{(-1,1)}
     \Bigl[ \sum_{n=1}^{+\infty} \alpha^n \beta^n \Bigr]
     v(\alpha) d \tilde\sigma (\alpha)
$$
$$   + \int\limits_{{\bf R}\backslash (-1,1)}
     d \tilde\sigma (\beta) \overline{ v(\beta) }
     \int\limits_{{\bf R}\backslash (-1,1)}
     \Bigl[ \sum_{n=0}^{+\infty} { 1 \over \alpha^n \beta^n }
     \Bigr]      v(\alpha) d \tilde\sigma (\alpha)
$$
$$ - {\rm Re} \, \int\limits_{(-1,1)}
     d\tilde\sigma(\beta) \overline{v(\beta)} \cdot
     {1\over\pi} \, \int\limits_{-\pi}^\pi
     \tilde r(e^{i\theta})
     \Bigl[ \sum_{n=1}^{+\infty} \beta^n e^{in\theta} \Bigr]
     v (e^{i\theta}) d \theta
$$
$$ + {\rm Re} \,
      \int\limits_{{\bf R}\backslash (-1,1)}
      d\tilde\sigma(\beta) \overline{v(\beta)} \cdot
     {1\over\pi} \, \int\limits_{-\pi}^\pi
     \tilde r(e^{i\theta})
     \Bigl[ \sum_{n=0}^{+\infty} { 1 \over \beta^n e^{in\theta} } \Bigr]
     v (e^{i\theta}) d \theta
$$
$$ + \|v_t\|^2
   + {\rm Re} \, {1\over\pi} \int\limits_{-\pi}^\pi d\phi \
     \overline{v(e^{i\phi})}
     \int\limits_{(-1,1)}
     \Bigl[ \sum_{n=1}^{+\infty} \alpha^n e^{in\phi} \Bigr]
     v(\alpha) d \tilde\sigma (\alpha)
$$
$$   + {\rm Re} \, {1\over\pi} \int\limits_{-\pi}^\pi d\phi \
     \overline{v(e^{i\phi})}
     \int\limits_{{\bf R}\backslash (-1,1)}
     \Bigl[ \sum_{n=0}^{+\infty} { 1 \over \alpha^n e^{in\phi} } \Bigr]
     v(\alpha) d \tilde\sigma (\alpha)
$$
$$  + {\rm Re} \, {1\over2} \sum_{n=0}^{+\infty} {1\over\pi} \int\limits_{-\pi}^\pi
       e^{-in\phi} \overline{v(e^{i\phi})} d\phi \cdot
       {1\over\pi} \int\limits_{-\pi}^\pi
        \tilde r(e^{i\theta}) v (e^{i\theta}) e^{-in\theta} d \theta
$$
$$  - {\rm Re} \,  {1\over2} \sum_{n=1}^{+\infty} {1\over\pi} \int\limits_{-\pi}^\pi
       e^{in\phi} \overline{v(e^{i\phi})} d\phi \cdot
       {1\over\pi} \int\limits_{-\pi}^\pi
        \tilde r(e^{i\theta}) v (e^{i\theta}) e^{in\theta} d \theta
$$

$$ = \|v_r\|^2  +  \sum_{n=1}^{+\infty}  \int\limits_{(-1,1)}
     \beta^n \overline{ v(\beta) } d \tilde\sigma (\beta)
     \int\limits_{(-1,1)} \alpha^n v(\alpha) d \tilde\sigma (\alpha)
$$
$$   + \sum_{n=0}^{+\infty} \int\limits_{{\bf R}\backslash (-1,1)}
     \beta^{-n} \overline{ v(\beta) } d \tilde\sigma (\beta)
     \int\limits_{{\bf R}\backslash (-1,1)}
     \alpha^{-n} v(\alpha) d \tilde\sigma (\alpha)
$$
$$ - {\rm Re} \, \sum_{n=1}^{+\infty} \int\limits_{(-1,1)}
     \beta^n \overline{v(\beta)} d\tilde\sigma(\beta) \cdot
     {1\over\pi} \, \int\limits_{-\pi}^\pi
     \tilde r(e^{i\theta}) v(e^{i\theta}) e^{in\theta} d \theta
$$
$$ + {\rm Re} \, \sum_{n=0}^{+\infty}
      \int\limits_{{\bf R}\backslash (-1,1)}
      \beta^{-n} \overline{v(\beta)} d\tilde\sigma(\beta) \cdot
     {1\over\pi} \, \int\limits_{-\pi}^\pi
     \tilde r(e^{i\theta}) v(e^{i\theta}) e^{-in\theta} d \theta
$$
$$ + \|v_t\|^2  +  {\rm Re} \, \sum_{n=1}^{+\infty}
     {1\over\pi} \int\limits_{-\pi}^\pi
     e^{in\phi} \overline{v(e^{i\phi})} d\phi \
     \int\limits_{(-1,1)} \alpha^n v(\alpha) d \tilde\sigma (\alpha) +
$$
$$ + {\rm Re} \, \sum_{n=0}^{+\infty}
     {1\over\pi} \int\limits_{-\pi}^\pi
     e^{-in\phi} \overline{v(e^{i\phi})} d\phi \
     \int\limits_{{\bf R}\backslash (-1,1)}
     \alpha^{-n} v(\alpha) d \tilde\sigma (\alpha)
$$
$$  + {\rm Re} \, {1\over2} \sum_{n=0}^{+\infty} {1\over\pi} \int\limits_{-\pi}^\pi
       e^{-in\phi} \overline{v(e^{i\phi})} d\phi \cdot
       {1\over\pi} \int\limits_{-\pi}^\pi
        \tilde r(e^{i\theta}) v (e^{i\theta}) e^{-in\theta} d \theta
$$
$$  - {\rm Re} \,  {1\over2} \sum_{n=1}^{+\infty} {1\over\pi} \int\limits_{-\pi}^\pi
       e^{in\phi} \overline{v(e^{i\phi})} d\phi \cdot
       {1\over\pi} \int\limits_{-\pi}^\pi
        \tilde r(e^{i\theta}) v (e^{i\theta}) e^{in\theta} d \theta \,
$$

$$  = {\rm Re} \, \Bigl\{ (\vec a,\vec a) + (\vec a,\vec b)
                        + (\vec c,\vec a) + {1\over2} (\vec b,\vec c) \Bigr\}
      + \|v_r\|^2 + \|v_t\|^2,
$$
that is
$$ {\rm Re} \, (Lv,v)
      = {\rm Re} \, \Bigl\{ (\vec a,\vec a) + (\vec a,\vec b)
                        + (\vec c,\vec a) + {1\over2} (\vec b,\vec c) \Bigr\}
      + \|v_r\|^2 + \|v_t\|^2
$$
where the vectors $\vec a, \vec b, \vec c \in l^2({\bf Z})$, whose
coordinates are defined by the formulas
$$  a_n = \int\limits_{{\bf R}\backslash (-1,1)}
     \alpha^n v(\alpha) d \tilde\sigma (\alpha),  \quad n\leq 0,
$$
$$  a_n = \int\limits_{(-1,1)}
     \alpha^n v(\alpha) d \tilde\sigma (\alpha),  \quad n > 0,
$$
$$  b_n = {1\over\pi} \, \int\limits_{-\pi}^\pi
     \tilde r(e^{i\theta}) v(e^{i\theta}) e^{in\theta} d \theta,
     \quad n \leq 0,
$$
$$  b_n = - {1\over\pi} \, \int\limits_{-\pi}^\pi
     \tilde r(e^{i\theta}) v(e^{i\theta}) e^{in\theta} d \theta,
     \quad n > 0,
$$
$$  c_n = {1\over\pi} \int\limits_{-\pi}^\pi
     \overline{v(e^{i\phi})} e^{-in\phi} d\phi,  \quad   n \in {\bf Z}.
$$
Here $b_n$ and $c_n$ are the Fourier coefficients in the expansion of the
functions $\sqrt{2} \tilde r(e^{i\theta}) v(e^{i\theta})$ and
$\sqrt{2} v(e^{i\phi})$ resp., defined under $(-\pi,\pi)$, by orthogonal systems
${ \Bigl\{  - {n \over |n|} { e^{in\theta} \over \sqrt{2} }
                             \Bigr\}_{n=-\infty}^{\infty}   }$ and
${\Bigl\{   { e^{in\phi} \over \sqrt{2} }
                             \Bigr\}_{n=-\infty}^{\infty}   }$.
That means,
$$ \|\vec b\|^2 \leq 2 \, {\rm max} \, |\tilde r(e^{i\theta})|^2 \cdot \|v_t\|^2 \,,
   \qquad
   \|\vec c\|^2 = 2 |v_t\|^2 \,.
$$
Now we use the identity
$$  {\rm Re} \, \Bigl\{ (\vec a,\vec a) + (\vec a,\vec b)
                        + (\vec c,\vec a) + {1\over2} (\vec b,\vec c) \Bigr\}
  = \| \vec a + {\vec b + \vec c  \over  2} \|^2
      - {1\over4} \Bigl( \|\vec b\|^2 + \|\vec c\|^2 \Bigr),
$$
which is true for any scalar product (can be simply verified starting from
the right-hand side). Applying this identity to our expression, we obtain
$$ {\rm Re} \, (Lv,v)
    = {\rm Re} \, \Bigl\{ (\vec a,\vec a) + (\vec a,\vec b)
                        + (\vec c,\vec a) + {1\over2} (\vec b,\vec c) \Bigr\}
      + \|v_r\|^2 + \|v_t\|^2 =
$$
$$  = \| \vec a + {\vec b + \vec c  \over  2} \|^2
      - {1\over4} \Bigl( \|\vec b\|^2 + \|\vec c\|^2 \Bigr)
      + \|v_r\|^2 + \|v_t\|^2
$$
$$  \ge \|v_r\|^2 + \|v_t\|^2 -
    {1\over4} \Bigl( 2 \, {\rm max} \, |\tilde r(e^{i\theta})|^2 \|v_t^2\|
       + 2 |v_t\|^2 \Bigr)
    \ge d \|v\|^2
$$
where $d = {1\over2} (1 - {\rm max} \, |\tilde r(e^{i\theta})|^2 )  >  0$.
Comparing this result with lemma 4.1, we see that  the operator $L$
is invertible and, moreover,
$$ \|L^{-1}\| \leq d^{-1}.
$$
Thus, the solvability of the equation (4.4), as well as (4.1), is proved.
So, we deduce

\bigskip

{\bf Theorem 4.} \ {\it Equations (3.46) and (4.1) are uniquely solvable
for every $k$, i.e.~the operator in the left-hand side of (4.1) is invertible
in the space $L^2({\bf R}_{\sigma} \cup {\bf T}_{\hat m\over \pi})$
and the inverse operator is bounded. }

\bigskip

  The collection
$\{ \hat r(e^{i\theta}),{m(\alpha)\over q(\alpha^{-1})},d\sigma(\alpha) \}$
consisting of the functions
$\hat r(e^{i\theta}),{m(\alpha)\over q(\alpha^{-1})}$ and the measure
$d\sigma(\alpha) $, which determine integral equation (3.46) and
representation (3.47), is called the {\it reduced spectral data} of
the Jacobi matrix $J$. The map
$$ J \mapsto n(z) \mapsto \{ \hat r(e^{i\theta}),
             {m(\alpha)\over q(\alpha^{-1})},d\sigma(\alpha) \} \,,
$$
described in the previous sections, is the solution of the direct spectral
problem, and the map
$$\{ \hat r(e^{i\theta}),{m(\alpha)\over q(\alpha^{-1})},d\sigma(\alpha) \}
   \mapsto J
$$
solves the inverse spectral problem.

   Summing up the results obtained in the preceding sections, we arrive
to the following theorem:

\bigskip

{\bf Theorem 5 \ (the main theorem).} \ {\it The infinite Jacobi matrix $J$
whose Weyl functions satisfy the conditions A)--E), is uniquely defined by
its reduced spectral data. The equations (3.46), reconstructed
according to these data, have a unique solution in the space
$L^2({\bf R}_{\sigma} \cup {\bf T}_{\hat m\over \pi})$ for all
$k\in{\bf Z}$. To reconstruct the Jacobi matrix $J$ according
to given spectral data it is necessary to solve equation (3.46) and then
use formulas (3.47) and (0.29).}

\bigskip
\bigskip

Acknowledgments. \ {\it This work is supported by
INTAS 2000-272 and by the State Foundation of Fundamental
Researches (Ukraine) №~1.4/20.}

\newpage

\centerline{\LARGE REFERENCES}

\bigskip

\begin{flushleft}
1. {\it S.P.~Novikov ed.,} Solyton theory. Inverse problem method.
Nauka, Moscow (1980).\\
2. {\it Teschl,} \, Jacobi operators and completely Integrable nonlinear lattices.
Math. Surv and Monographs, vol. 72. (2000)\\
3. {\it M.~Kac and P.~van Moerbeke,} \, Adv. in Math. (1975), v.~16, №~2,
p.~160--169.\\
4. {\it Yu.~Berezanski,} \, The integration of semi-infinite Toda chain
by means of inverse spectral problem. Reports on mathematical physycs
(1986). v.~24, p.~21--46.\\
5. {\it N.~Zhernakov,} \,  The integration of the Toda lattice in the
class if Hilbert-Schmidt operators. Ukrainskiy matematicheskiy zhurnal
(1987), v.~39, №~5, p.~645--646.\\
6. {\it N.~Akhiezer,} \, Classical moments problem. "FIZMATGIZ", Moskow
(1961).\\
7. {\it Yu.~Березанский,} \, Eigenfunction expansion of the self-adjoint
operators. "Naukova dumka", Kiev (1965).\\
8. {\it M.Kudryavtsev,} \, Resolution of the Cauchy problem for the Toda
lattice with non-stabilized initial data. In preparation.\\
9. {\it M.Kudryavtsev,} \, Inverse problems for for finite-difference
operators of mathematical physics and their application. PhD thesis,
The Institute for Low Temperature Physics and Engineering, Kharkov,
Ukraine (2000). \\
10. {\it M.~Kudryavtsev,} \,  Solution of the Cauchy problem for a Toda chain
with initial data that are not stabilized. Dokl. NAN Ukrainy
(2001) №~3, p.~14-19. \\
11. {\it L.~Maslov,} \, Spectral properties of the reflectionless Jacobi
matrices. Teoriya funktsiy, funktsionalnyy analiz i
ih prilozheniya (1992), №~57, p.~25--46.\\
12. {\it L.~Maslov,} \, Description of the set of strong limits of the
reflectionless Jacobi matrices. Teoriya funktsiy, funktsionalnyy
analiz i ih prilozheniya (1993). №~58. p.~54--60.\\
13. {\it A.~Boutet de Monvel and V.~Marchenko,} \,
The Cauchy problem for nonlinear Shr\"odinger equation with bounded
initial data. Mat. fiz., an., geom. (1997), v.4, №~1/2, p.~3--45. \\
14. {\it N.~Akhiezer and I.~Glazman,} \, Theory of linear operators in
Hilbert spaces. "Nauka", Moskow (1966).\\
15. {\it M.A.~Kudryavtsev,} \, The Riemann problem with additional
singularities (in Russian). Mat. fiz., an., geom. (2000). v.~7, №~2, p.~196--208.
\ English translation: http://arXiv.org/abc/math.SP/0110242 \\
16. {\it M.G.~Krein M.G.,} The Markov moments problem and extremal
problems. "Nauka", Moscow (1973).\\
17. {\it Muskhelishvili,} Singular integral equations. "Fizmatgiz", Moscow (1962).\\
18. {\it J.~Garnett,} \, Bounded analytic functions. Academic press, New
York, London (1981).\\
19. {\it P.~Koosis,} \, Introduction in theory of spaces $H^p$. "Mir", Moscow
(1984).\\
20. {\it V.A.~Marchenko,} \, Nonlinear equations and operator algebras.
"Naukova dumka", Kiev (1986).\\
\end{flushleft}

\end{large}
\end{document}